\documentclass[12pt,leqno]{article}

=155mm\textheight =220mm

\usepackage{amsmath}
\usepackage{amscd}
\usepackage{amsmath}
\usepackage{amssymb}
\usepackage{latexsym}
\makeatletter

\@addtoreset{equation}{section}
\makeatother

\newtheorem{thm}{Theorem}[subsection]
\newtheorem{pr}[thm]{Proposition}
\newtheorem{df}[thm]{Definition}
\newtheorem{lm}[thm]{Lemma}
\newtheorem{cor}[thm]{Corollary}

\newtheorem{eg}[thm]{Example}

\newcommand{\qed}{\hfill{\rule{7pt}{7pt}}}

\input xy
\xyoption{all}
\CompileMatrices

\begin{document}

\title{Wild ramification
and the characteristic cycle\\
of an $\ell$-adic sheaf}
\author{Takeshi Saito}
\maketitle

\begin{abstract}
We propose a geometric method to
measure the wild ramification
of a smooth \'etale sheaf
along the boundary.
Using the method,
we study
the graded quotients
of the logarithmic 
ramification groups
of a local field of characteristic $p>0$
with arbitrary residue field.
We also define the characteristic
cycle of an $\ell$-adic sheaf,
satisfying certain conditions,
as a cycle on the logarithmic cotangent
bundle and prove that
the intersection with the 0-section
computes the characteristic class,
and hence the Euler number.
\end{abstract}

Let $X$ be a separated
scheme of finite type over
a perfect field $k$
of characteristic $p>0$.
We consider 
a smooth $\ell$-adic \'etale sheaf
${\cal F}$
on a smooth dense open subscheme
$U\subset X$ for 
a prime
$\ell\neq p$.
The ramification of
${\cal F}$ along the boundary
$X\setminus U$
has been studied
traditionally by
using a finite \'etale covering
of $U$ trivializing
${\cal F}$ modulo $\ell$.
In this paper,
we propose a new geometric method,
inspired by the definition
of the ramification groups
\cite{AS}, \cite{AS2} and \cite{aml}.

The basic geometric construction
used in this paper
is the blowing-up
at the ramification
divisor embedded diagonally
in the self log product.
A precise definition will be given at
the beginning of \S2.3.
We will consider two types of blow-up.
The preliminary one,
called the log blow-up,
is the blow-up
$(X\times X)'\to X\times X$
at every $D_i\times D_i$
where $D_i$ denotes
an irreducible component
of a divisor
$D=X\setminus U$ with simple normal crossings
in a smooth scheme $X$ over $k$.
The second one
is the blow-up
$(X\times X)^{(R)}\to (X\times X)'$
at $R=\sum_ir_iD_i$,
with some rational multiplicities $r_i\ge0$,
embedded in the log
diagonal $X\to 
(X\times X)'$.
This construction
globalizes
that used in the definition
of the ramification
groups in \cite{AS} and \cite{AS2}
recalled in \S1.

Inspired by \cite{ksc},
we consider
the ramification along the boundary
of the smooth sheaf
${\cal H}={\cal H}om(
{\rm pr}_2^*{\cal F},
{\rm pr}_1^*{\cal F})$
on the dense open subscheme
$U\times U\subset (X\times X)^{(R)}$.
We introduce a measure of
wild ramification by
using the extension property
of the identity
regarded as a section
of the restriction 
on the diagonal of
the sheaf ${\cal H}$,
in Definition \ref{dfbdR}.

Let $j^{(R)}:
U\times U\to (X\times X)^{(R)}$
denote the open immersion.
A key property of the sheaf
${\cal H}$ established in Propositions
\ref{prD} and \ref{prR}
is that the restriction
of $j^{(R)}_*{\cal H}$
on the complement
$(X\times X)^{(R)}\setminus U\times U$
admits a description
by the Artin-Schreier sheaves
defined by certain linear forms.
This fact is derived
from a groupoid structure
of $(X\times X)^{(R)}$
inherited from the natural
one on $X\times X$.
We prove in Theorem \ref{thmgr} that
this property at the generic point
of an irreducible component
implies the following
properties
of the ramification groups
conjectured in \cite{aml} Conjecture 9.4:
The graded pieces
of the ramification groups
are killed by $p$
and
their character groups
are described by
differential forms.

The definition of the measure
of the wild ramification in this paper
is closely related
to that of the characteristic 
class in \cite{cc}.
In Definition \ref{dfCC},
we propose a definition
of the characteristic cycle
of an $\ell$-adic sheaf
as a cycle of the
logarithmic cotangent bundle,
under the conditions (R) and (C)
stated in \S3.2. 
Roughly speaking, 
the conditions mean
that the ramification is 
controlled at the generic points
of the irreducible components
of the ramification divisor.
Consequently, the characteristic cycle
in this case
does not have components supported
on subvarieties of codimension 
at least two.
We show that
its intersection product
with the 0-section
computes the characteristic class,
in Theorem \ref{thmCC}.
This is a generalization
of Kato's formula in the rank one case
\cite{Kato2}.

One expects that
the same construction
works for ${\cal D}$-modules
with irregular singularities.
It should give another evidence
for the analogy between 
the wild ramification of
$\ell$-adic sheaves
and irregular singularities
of ${\cal D}$-modules.

The author would like
to express his sincere gratitude
to Ahmed Abbes for
stimulating discussions.
He is pleased
to acknowledge that 
a large part of this paper
is based on a colaboration
with him.
The author also thanks to him
for the information on
the reference \cite{Oe}
on Epp's theorem \cite{epp}.
The author thanks an anonymous referee
for pointing a mistake in
the definition \ref{dfadd}
in an earlier version.

This research is partly supported by
JSPS Grant-in-Aid for Scientific Research (B) 14340002 and
18340002.

\tableofcontents

\section*{Notation}
$k$ denotes a perfect field
of characteristic $p>0$.
A scheme over $k$
is assumed to be
separated of finite type
over $k$.
For a locally free ${\cal O}_X$-module
${\cal E}$ of finite rank 
on a scheme $X$,
$E={\mathbf V}({\cal E})$
denotes the contravariant
vector bundle
defined by the quasi-coherent
${\cal O}_X$-algebra
$S^{\bullet}{\cal E}$.
Similarly,
${\mathbf P}({\cal E})$
denotes the projective space
bundle
${\cal P}roj\ S^{\bullet}{\cal E}$.
The dual of ${\cal E}$
is denoted by ${\cal E}^\vee$.
For a closed subscheme
$X\subset Y$ defined
by the ideal ${\cal I}_X\subset {\cal O}_Y$,
the conormal sheaf
${\cal I}_X/{\cal I}_X^2$ is denoted
by ${\cal N}_{X/Y}$.

$\ell$ denotes a prime number
invertible in $k$ and
$\Lambda$ denotes
a finite local
${\mathbb Z}_\ell$-algebra.

\setcounter{section}0

\section{Ramification groups}

The theory of logarithmic ramification groups 
of a local field with imperfect residue field
as developed in \cite{AS} and \cite{AS2}
relies on rigid geometry 
and on log geometry in
an essential way.
In \S\S1.1-1.3,
we give
some interpretations 
purely in terms
of schemes,
without using rigid geometry
or log geometry.
In \S1.3,
we state the main result
Theorem \ref{thmgr}
on the structure of
the graded quotients.
We prove it
in \S1.4 by computing
the nearby cycles.

In this section,
$K$ denotes a discrete valuation field,
${\cal O}_K$ denotes the valuation ring,
and ${\mathfrak m}_K$ denotes the maximal ideal.
The residue field
${\cal O}_K/{\mathfrak m}_K$ is denoted by $F$
and $v_K:K \to {\mathbb Z}\cup\{\infty\}$
denotes the discrete valuation
normalized by $v_K(\pi)=1$
for a prime element $\pi$.
We put $S={\rm Spec}\ {\cal O}_K$.

\subsection{Basic constructions}
Let $A$ be a finite flat ${\cal O}_K$-algebra.
We put $T={\rm Spec}\ A$.
We consider a closed immersion
$T\to P$
to a smooth scheme $P$ over $S$.
Let ${\cal I}_T={\rm Ker}({\cal O}_P\to {\cal O}_T)$
be the ideal sheaf defining
the closed subscheme $T$ in $P$.

For a pair $(m,n)$
of integers $m\ge 0$ and $n>0$,
let $Q=P_T^{[m/n]}\to P$ 
be the blow-up at the ideal 
${\cal I}_T^n+{\mathfrak m}_K^m{\cal O}_P$
and $P_T^{(m/n)}\subset P_T^{[m/n]}$
be the complement
of the support of
$({\cal I}_T^n{\cal O}_Q+{\mathfrak m}_K^m{\cal O}_Q)/
{\mathfrak m}_K^m{\cal O}_Q$.
The morphism
$P_T^{(m/n)}\to P$
is affine and 
$P_T^{(m/n)}$
is defined by the
quasi-coherent sub ${\cal O}_P$-algebra
${\cal O}_P[{\mathfrak m}_K^{-m}{\cal I}_T^n]
\subset K\otimes {\cal O}_P$.
The maps 
$P_T^{(m/n)}\to P_T^{[m/n]}\to P$
induce isomorphisms
$P_{T,K}^{(m/n)}\to P_{T,K}^{[m/n]}\to P_K$
on the generic fibers.
For $m=0$, we have
$P_T^{(0/n)}= P_T^{[0/n]}= P$.
The immersion $T\to P$
is uniquely lifted to an immersion
$T\to P^{(m/n)}$.

Let $d>0,m'\ge 0,n'>0$
be integers such that
$m'\le dm$ and $n'=dn$.
Then the inclusion
$({\mathfrak m}_K^{-m}{\cal I}_T^n)^d\supset 
{\mathfrak m}_K^{-m'}{\cal I}_T^{n'}$
induces a canonical map
$P_T^{(m/n)}\to P_T^{(m'/n')}$
that is an isomorphism on
the generic fibers.
If $(m',n')=(dm,dn)$,
the canonical map
$P_T^{(m/n)}\to P_T^{(m'/n')}$
is finite.

For a rational number
$r=m/n\ge 0$,
let $P^{(r)}_T$ be the normalization
of $P_T^{(m/n)}$.
For $r>0$,
let $\widehat {P^{(r)}_T}$
be the formal completion
of $P^{(r)}_T$ along the closed fiber
$P^{(r)}_{T,F}$.
For $r'\le r$,
the canonical maps
$P_T^{(r)}\to P_T^{(r')}$
of schemes and
$\widehat {P^{(r)}_T}\to
\widehat {P^{(r')}_T}$
of affine formal schemes
are induced.

We compare the construction above
to those in \cite{AS} and \cite{AS2}.

\begin{eg}
Assume $K$ is complete.

1. Let 
$Z=(z_1,\ldots,z_n)$
be a system of generators of 
a finite flat ${\cal O}_K$-algebra $A$
and consider the closed immersion
$T={\rm Spec}\ A\to 
P={\rm Spec}\ {\cal O}_K[X_1,\ldots,X_n]$
defined by $Z$.
Then, the affinoid variety $X_Z^a$
in \cite{AS} 3.1 
is defined by the formal ${\cal O}_K$-scheme
$\widehat {P^{(r)}_T}$
for $a=r$.

2. Let $T\to P$
be a closed immersion of a finite flat
${\cal O}_K$-scheme $T$ to a smooth scheme $P$
and let ${\rm Spf}\ {\mathbf A}$ be 
the formal completion $\widehat P|_T$
of $P$ along the closed subscheme $T$.
Then the 
affinoid variety $X^j({\mathbf A}\to A)$
in \cite{AS2} Definition 1.5
is defined by the formal ${\cal O}_K$-scheme
$\widehat {P^{(r)}_T}$ for $j=r$.
\end{eg}

\begin{lm}\label{lmhtp}
Let $T$ be a finite flat
scheme over $S$ and
$T\to P$ and $T\to Q$
be closed immersions
to smooth schemes over $S$.
Let $P\to Q$ be
a smooth morphism over $S$
such that the diagram
$$\xymatrix{
T\ar[r]\ar[rd]&P\ar[d]\\
&Q}$$
is commutative.
Then, for a positive integer $r>0$,
the map $P\to Q$
induces a smooth map
$P^{(r)}_T\to Q^{(r)}_T$
and an isomorphism
$$\begin{CD}
P^{(r)}_{T,F}
@>>>
Q^{(r)}_{T,F}\times_{T_F}
{\mathbf V}(
{\mathfrak m}_K^{-r}
\otimes_{{\cal O}_K}
\Omega^1_{P/Q}\otimes_{{\cal O}_P}{\cal O}_{T_F}).
\end{CD}$$
\end{lm}
{\it Proof.}
We show the map
$P^{(r)}_T\to Q^{(r)}_T$ is smooth.
Let $t\in T$
be a closed point
and $d$ be the relative dimension
of $P\to Q$ at the image of $t$.
Since
$T\to P$
defines a section
of the smooth morphism
$P\times_QT\to T$,
there exist an neighborhood
$V\subset P$ of the image of $t$
and an \'etale morphism
$V\to {\mathbf A}^d_Q=
Q[X_1,\ldots,X_d]$
inducing an open immersion
$T\cap V
\to 
T\subset
T\times_Q{\mathbf A}^d_Q$
to the 0-section.
Then,
$P^{(r)}_T
\times_PV$
is 
isomorphic to
$V\times_{{\mathbf A}^d_Q}
Q^{(r)}_T
[X_1/\pi^r,\ldots,X_d/\pi^r]$
and is smooth over
$Q^{(r)}_T$.

Since
the map
$Q^{(r)}_{T,F}
\to Q_{T,F}$
factors through
the closed immersion
$T_F\to Q_{T,F}$,
the isomorphism
$P^{(r)}_T
\times_PV
\to
V\times_{{\mathbf A}^d_Q}
Q^{(r)}_T
[X_1/\pi^r,\ldots,X_d/\pi^r]$
above induces
an open immersion
$P^{(r)}_{T,F}
\times_PV
\to
Q^{(r)}_{T,F}
[X_1/\pi^r,\ldots,X_d/\pi^r]$.
Since
$Q^{(r)}_{T,F}
[X_1/\pi^r,\ldots,X_d/\pi^r]$
is canonically identified
with 
$Q^{(r)}_{T,F}\times_{T_F}
{\mathbf V}(
{\mathfrak m}_K^{-r}
\otimes_{{\cal O}_K}
\Omega^1_{P/Q}\otimes_{{\cal O}_P}{\cal O}_{T_F})$,
the assertion follows.
\qed

Let $\overline T$ denote
the normalization of $T$.
For positive integers $m,n>0$
and for $r=m/n$,
the immersion $T\to P$
induces an immersion $T\to P^{(m/n)}_T$
and hence a finite maps
$\overline T
\to P^{(r)}_T$
and
$\overline T
\to \widehat{P^{(r)}_T}$.

\begin{lm}\label{lmspl}
Let $T={\rm Spec}\ A$ be a finite flat
scheme over $S$ and
$T\to P$ be
a closed immersion 
to a smooth scheme over $S$.
Assume that $T_K$
is isomorphic to
the disjoint union of finitely many copies of
${\rm Spec}\ K$.
Then there exists an integer
$r>0$
such that 
the map $\overline T
\to P^{(r)}_T$ is a closed immersion.
\end{lm}
{\it Proof.}
By the assumption on $T_K$,
the semi-local ring $A$
is the product of finitely many
local rings.
Hence, we may assume $P={\rm Spec}\ R$
is affine and hence
$P^{(r)}_T={\rm Spec}\ R^{(r)}$
is also affine.
Further by the assumption
on $T_K$, the normalization
of $A$ is generated by
the idempotents in
$A\otimes_{{\cal O}_K}K$.
Hence, it is sufficient to show
that, for every idempotent
$e\in A\otimes_{{\cal O}_K}K$,
there exists an integer $r>0$
such that 
$e$ is in the image of $R^{(r)}\to
A\otimes_{{\cal O}_K}K$.
Take a non-zero element $a
\in {\mathfrak m}_K$ such that
$ae\in A$.
We show that $r=2v_K(a)$
satisfies the condition.

Take a lifting $f\in R$
of $ae\in A$.
Then $g=f^2-af$ is in the
kernel $I={\rm Ker}(R\to A)$.
Since $g/a^2$ is in
$R^{(r/1)}$, the solution
$f/a\in R^{(r/1)}\otimes_{{\cal O}_K}K$
of the equation $X^2-X=g/a^2$
lies in $R^{(r)}$.
\qed

We study the relation
of the basic construction
with a base change
of discrete valuation rings.
Let $T\to P$ be a closed immersion
of a finite scheme to
a smooth scheme
over $S={\rm Spec}\ {\cal O}_K$ as above.
Let $S'={\rm Spec}\ {\cal O}_{K'}\to S$
be a surjection of spectra
of discrete valuation rings
of ramification index $e$.
Then, the base change
$T'=T\times_SS'\to P'=P\times_SS'$
is a closed immersion
of a finite flat scheme to
a smooth scheme
over $S'$.
For integers $m,n>0$,
the induced map
$P^{\prime [em/n]}_{T'}
\to 
P^{[m/n]}_T\times_S
S'$ is an isomorphism.
Hence, for $r=m/n$,
the scheme
$P^{\prime (er)}_{T'}$
is the normalization
of $P^{(r)}_T\times_SS'$
and the formal scheme
$\widehat{P^{\prime (er)}_{T'}}$
is the normalization of
$\widehat{P^{(r)}_T}
\times_{\widehat S}\widehat S'$.

We prepare some facts on
the properties 
(S$_k$) and (R$_k$)
of locally noetherian schemes
\cite{EGA} Chap.\ IV \S\S 5.7, 5,8.

\begin{lm}\label{lmGR}
Let $f:X\to S$ be a flat scheme of
finite type over a regular noetherian
scheme $S$.
For a point $s\in S$,
we put $c(s)=\dim {\cal O}_{S,s}$.
Let $k\ge 0$ be an integer.

1. The following conditions
are equivalent.

{\rm (1)} 
For every point $s\in S$,
the fiber $X_s$ satisfies
the condition
{\rm (S}$_{k-c(s)}${\rm )}.

{\rm (2)} 
$X$ satisfies 
the condition 
{\rm (S}$_k${\rm )}.

2. The condition {\rm (1)} 
implies the condition {\rm (2)}.

{\rm (1)} 
For every point $s\in S$,
the fiber $X_s$ satisfies
the condition {\rm (R}$_{k-c(s)}${\rm )}.

{\rm (2)} 
$X$ satisfies 
the condition 
{\rm (R}$_k${\rm )}.
\end{lm}
{\it Proof.}
1.
Let $x\in X$ be a point
and put $s=f(x)\in S$.
Let $t_1,\ldots,t_c\in {\mathfrak m}_s$
be a regular system of parameters
where $c=c(s)$.
Since $f:X\to S$
is flat,
$f^*t_1,\ldots,f^*t_c\in {\mathfrak m}_x$
is a regular sequence of ${\cal O}_{X,x}$
and ${\cal O}_{X_s,x}={\cal O}_{X,x}/(f^*t_1,\ldots,f^*t_c)$.
Hence, we have equalities
$\dim {\cal O}_{X_s,x}=\dim {\cal O}_{X,x}-c(s)$
(\cite{EGA} Chap.\ 0 Proposition (16.3.4)) and
${\rm prof}\ {\cal O}_{X_s,x}={\rm prof}\ {\cal O}_{X,x}-c(s)$
(loc.\ cit. Proposition (16.4.6) (ii))
and the assertion follows..

2. Further,
${\cal O}_{X,x}$ is regular
if ${\cal O}_{X_s,x}$ is regular
(loc.\ cit.\ Proposition (17.3.3 (ii))).
\qed

\begin{cor}\label{corEpp}
Let $S={\rm Spec}\ {\cal O}_K$ 
be the spectrum of a discrete valuation
ring and $f:X\to S$ be
a normal scheme of finite type
with smooth generic fiber.

1. 
There exists a surjection of spectra
$S'={\rm Spec}\ {\cal O}_{K'}\to S$
of discrete valuation rings
such that $K'$ is a finite 
extension of $K$
and that the normalization 
$X'$ of $X\times_SS'$
has geometrically reduced fibers
over $S'$.

2. Assume
$X\to S$
has geometrically reduced fibers.
Then, for any
surjection $S'\to S$ of spectra
of discrete valuation rings,
the base change
$X\times_SS'$
is normal.
\end{cor}
{\it Proof.}
1. 
We apply a variant,
l'appendice Th\'eor\`eme 2 \cite{Oe},
of Epp's theorem \cite{epp}
corrected in \cite{corr}
to the generic points
of the irreducible 
components of
the closed fiber of $X\to S$.
Then, we fine
a surjection $S'={\rm Spec}\ {\cal O}_{K'}\to S$
of spectra of discrete valuation rings
and an open subscheme
$U$ of the normalization 
$X'$ of the base change 
$X\times_{S}S'$ such that
$K'$ is a finite extension
of $K$ and that
$U$ is smooth 
over $S'$ and contains
the generic point of 
every irreducible component
of the closed fiber.

We show that $X'$ has geometrically reduced
fibers.
Since the generic fiber is smooth,
it suffices to show 
that the geometric closed fiber is reduced.
Since $X'$ is normal,
it satisfies the condition ${\rm (S_2)}$.
By Lemma \ref{lmGR}.1, the closed fiber 
satisfies ${\rm (S_1)}$ and
hence the geometric closed fiber also
satisfies ${\rm (S_1)}$.
Since the geometric closed fiber has
a dense open subscheme smooth over the base field,
it also satisfies the condition ${\rm (R_0)}$.
Hence the geometric closed fiber is reduced.

2.
Since the closed fiber of
$X_{S'}$ is reduced,
it satisfies 
the conditions ${\rm (R_0)}$
and ${\rm (S_1)}$.
Since the generic fiber is regular,
$X_{S'}$ satisfies
the conditions ${\rm (R_1)}$
and ${\rm (S_2)}$
by Lemma \ref{lmGR}.
\qed

Let $X$ be a normal scheme 
of finite type over $S={\rm Spec}\ {\cal O}_K$.
Assume that
the generic fiber of $X$ is
smooth
and that the closed geometric fiber
is reduced.
Then, the formal completion
$\widehat X$ along the closed fiber
is the stable integral model
of the affinoid variety
defined by $\widehat X$ itself.
Thus, Corollary \ref{corEpp}.1
implies the finiteness theorem
of Grauert-Remmert
for algebraizable formal schemes.

Applying Corollary
\ref{corEpp} to
$P^{(r)}_T\to S$,
we obtain the following.

\begin{cor}\label{corGR}
Let $T\to P$ be a closed immersion
of a finite scheme to
a smooth scheme
over $S={\rm Spec}\ {\cal O}_K$
and $r>0$ be a rational number.

1. There exists a
surjection of spectra
$S'={\rm Spec}\ {\cal O}_{K'}\to S$
of discrete valuation rings
of ramification index $e$
such that $K'$ is a finite 
extension of $K$
and that 
$P^{\prime (er)}_{T'}
\to S'$
has geometrically reduced fibers.

2. Assume
$P^{(r)}_T\to S$
has geometrically reduced fibers.
Then, for any
surjection $S'\to S$ of spectra
of discrete valuation rings
of ramification index $e$,
the canonical map
$P^{\prime (er)}_{T'}
\to
P^{(r)}_T\times_SS'$
is an isomorphism.
\end{cor}

\begin{df}\label{dfstb}
Let $T$ be a finite flat
scheme over $S$ and
$T\to P$ 
be a closed immersion
to a smooth scheme over $S$.
Let $r>0$ be a rational number
and $S'\to S$
be a surjection of spectra
of discrete valuation rings
of ramification index $e$.

We say $P^{\prime (er)}_{T'}
\to S'$ is a stable model
of $P^{(r)}_T$ 
if its geometric fibers
are reduced.
If $P^{\prime (er)}_{T'}
\to S'$ is a stable model,
we call 
$P^{\prime (er)}_{T'}
\times_{S'}{\rm Spec}\ {\overline F}$
the stable
closed fiber
and write it by
$\overline P^{(r)}_{T,\overline F}$.
\end{df}

By Corollary \ref{corGR}.1,
there exists an $S'$
such that $P^{\prime (er)}_{T'}
\to S'$ is a stable model.
By Corollary \ref{corGR}.2,
the stable
closed fiber
$\overline P^{(r)}_{T,\overline F}$
is independent of
the choice of such $S'$.
The finite map
$P^{\prime (er)}_{T'}
\to
P^{(r)}_T\times_SS'$
induces a finite map
$\overline P^{(r)}_{T,\overline F}
\to
P^{(r)}_{T,\overline F}$.

Similarly
as in the definition of
the stable closed fiber
$\overline P^{(r)}_{T,\overline F}$, 
if $T_K$ is \'etale over $K$,
for $S'={\rm Spec}\ {\cal O}_{K'}\to S$ such that
$T\times_S{\rm Spec}\ K'$ is the
disjoint union of finitely many copies
of ${\rm Spec}\ {K'}$,
the geometric fiber
$\overline {T\times_S{S'}}\times_{S'}
\overline F$ of the normalization
is independent of
the choice of $S'$.
We write it by
$\overline T_{\overline F}$.

\begin{df}\label{dfram}
Let $T$ be a finite flat
scheme over $S$ such that
$T_K$ is \'etale over $K$.

1. Let $r>0$ be a rational number.
Let $T\to P$ be a closed immersion
to a smooth scheme over $S$
and $S'={\rm Spec}\ {\cal O}_{K'}\to S$
be a surjection of spectra
of discrete valuation rings
of ramification index $e$
satisfying the following conditions:
The \'etale covering $T_K
\to {\rm Spec}\ K$
splits over $K'$
and hence the normalization
$\overline T_{S'}$
of $T\times_SS'$
is isomorphic to
the disjoint union of finitely many
copies of $S'$.
The product $e r$ is an integer
and the geometric fibers
of $P^{(er)}_{T,S'}\to {S'}$
are reduced. 

We say the ramification of
$T$ over $S$
is bounded by $r$ if,
the map
$\overline T_{S'}
\to P^{(er)}_{T,{S'}}$
induces an injection
$$\overline T_{\overline F}
\to \pi_0(\overline P^{(r)}_{T,\overline F})$$
of finite sets.

2. 
Let $r\ge 0$ be a rational number.
We say the ramification of
$T$ over $S$
is bounded by $r+$ if
the ramification of
$T$ is bounded by every rational number
$s>r$.
\end{df}

By Lemma \ref{lmhtp},
the map $\overline T_{\overline F}
\to \pi_0(\overline P^{(r)}_{T,\overline F})$
is independent of $P$.
Let $T$ be a finite flat
scheme over $S$
and $S'\to S$
be a surjection
of spectra of
discrete valuation rings
of ramification $e$.
Then, it is clear from
the definition that
the ramification of
$T$ over $S$ is bounded by
$r$ if and only if
the ramification of
$T\times_SS'$ over $S'$ is bounded by
$er$.

We will see later that Definition
\ref{dfram}
is equivalent to the definition
in \cite{AS} Definition 6.3
for finite flat ${\cal O}_K$-algebra
locally of complete intersection.

\begin{lm}\label{lmlci}
For a finite flat
scheme $T$ over $S$,
the following conditions are
equivalent.

{\rm (1)} $T$ is locally of 
complete intersection.

{\rm (2)} There exists
a cartesian diagram
\begin{equation}\begin{CD}
T@>>> Q\\
@VVV @VVV\\
S@>>> P
\end{CD}
\label{eqTS}
\end{equation}
of schemes over $S$
satisfying the following condition:

{\rm (CI)}
The vertical arrows
are quasi-finite flat and
the horizontal arrows
are closed immersions.
The schemes $P$ and $Q$ are smooth over $S$.
\end{lm}
{\it Proof.}
(1) $\Rightarrow$ (2)
Take a surjection
${\cal O}_K[X_1,\ldots,X_d]
\to A$.
Then the closed immersion
$T\to Q={\mathbf A}^d_S$
is regular of
codimension $d$.
By Nakayama's lemma,
there exists elements
$f_1,\ldots,f_d
\in I={\rm Ker}({\cal O}_K[X_1,\ldots,X_d]
\to A)$
such that
$(f_1,\ldots,f_d)=I$
on a neighborhood of $T$.
We define a map
$Q\to P={\mathbf A}^d_S$
by $f_1,\ldots, f_d$.
Then, shrinking $Q$ if necessary,
the diagram 
(\ref{eqTS}) is cartesian
and the map $Q\to P$
is quasi-finite 
and flat.

(2) $\Rightarrow$ (1)
Since the immersion $S\to P$
is regular,
the immersion $T\to Q$
is also regular and 
$T$ is
locally of complete intersection over $S$.
\qed

We compute the scheme
$P^{(r)}_T$
explicitly in the case where $T=S
\to P$ is a section
of a smooth scheme $P\to S$
of relative
dimension $d$.
The conormal sheaf
${\cal N}_{S/P}={\cal I}_S/{\cal I}_S^2$ is canonically
identified with the 
free ${\cal O}_K$-module
$\Omega^1_{P/S}\otimes_{{\cal O}_P}{\cal O}_S$
of rank $d$.

\begin{lm}\label{lmmn}
Let $S\to P$ be a section
of a smooth scheme $P\to S$ and
$r>0$ be a rational number.
Let $j:P_K=
P\times_SK\to P$
be the open immersion
and ${\cal I}_S
\subset {\cal O}_P$ be the ideal sheaf
of $S$ regarded as a subscheme
of $P$ by the section $s:S\to P$.

1. The affine $P$-scheme
$P^{(r)}_S$ is defined 
by the quasi-coherent
${\cal O}_P$-algebra
\begin{equation}
\sum_{l\ge 0}{\mathfrak m}_K^{-[lr]}\cdot {\cal I}_S^l
\subset j_*{\cal O}_{P_K}
\label{eqPr}
\end{equation}
where $[lr]$ denotes the integral part.

2. Assume $r$
is an integer.
Then $P^{(r)}_S=
P^{(r/1)}_S$ is smooth
over ${\cal O}_K$ and the closed fiber
$P^{(r)}_{S,F}$
is identified with the
$F$-vector space
$\Omega^1_{P/S}\otimes_{{\cal O}_P}
F\otimes_F
{\mathfrak m}_K^{-r}/
{\mathfrak m}_K^{-r+1}$.
\end{lm}
{\it Proof.}
1. Let $n\ge 1$ be an integer
such that $m=nr$ is an 
integer. Then,
$P^{(r)}_S$ is
defined by the normalization $A$
of the quasi-coherent ${\cal O}_P$-algebra
${\cal O}_P[{\mathfrak m}_K^{-mr}\cdot {\cal I}_S^n]
\subset j_*{\cal O}_{P_K}$.
Since
$\sum_{l\ge 0}{\mathfrak m}_K^{-[lr]}\cdot {\cal I}_S^l$
is integral over
${\cal O}_P[{\mathfrak m}_K^{-mr}\cdot {\cal I}_S^n]$,
we have an inclusion
$\sum_{l\ge 0}{\mathfrak m}_K^{-[lr]}\cdot {\cal I}_S^l
\subset A$.

We show the inclusion is
an equality.
Since the question is \'etale
local on $P$,
we may assume $P={\mathbf V}(M)=
{\rm Spec}\ S^\bullet(M)$
for a free ${\cal O}_K$-module of finite rank
and $S\to P$
is the 0-section of the vector bundle
$P$ over $S$.
Then, we have
$P^{(r)}_S=
{\rm Spec}\ \bigoplus_{l\ge 0}
{\mathfrak m}_K^{-[lr]}S^l(M)$.

2. We show $P^{(r)}_S=
P^{(r/1)}_S$ is smooth
over ${\cal O}_K$.
Since the question is \'etale
local on $P$,
we may assume $P={\mathbf V}(M)=
{\rm Spec}\ S^\bullet(M)$
as above.
Then, 
$P^{(r)}_S=
P^{(r/1)}_S=
{\rm Spec}\ 
S^\bullet ({\mathfrak m}_K^{-r}M)$
is smooth over ${\cal O}_K$.

We show that the closed immersion
$P^{(r)}_{S,F}=
P^{(r/1)}_{S,F}\to
{\rm Spec}\ 
S^\bullet ({\mathfrak m}_K^{-r}\otimes_{{\cal O}_K}
{\cal N}_{S/P}\otimes_{{\cal O}_K}F)$
is an isomorphism.
We conclude
by reducing to the case
$P={\mathbf V}(M)=
{\rm Spec}\ S^\bullet(M)$
as above.
\qed

Let $v:{\overline K}\to {\mathbb Q}
\cup \{\infty\}$
be the extension
of the normalized discrete valuation
$v:K\to {\mathbb Z}
\cup \{\infty\}$
to a separable closure.
For a rational number $r$,
we put
${\mathfrak m}^r_{\overline K}=
\{a\in {\overline K}|v(a)\ge r\}$
and
${\mathfrak m}^{r+}_{\overline K}=
\{a\in {\overline K}|v(a)> r\}$.

\begin{cor}\label{cormn}
Let $m,n>0$ be positive integers 
such that $r=m/n$ and $(m,n)=1$.
Then for the reduced closed fiber
$(P^{(r)}_{S,F})_{\rm red}$,
we have a commutative diagram
\begin{equation}
\begin{CD}
\overline P^{(r)}_{S,\overline F}
@>>>
{\mathbf V}(
\Omega^1_{P/S}\otimes_{{\cal O}_P}
\overline F
\otimes_{\overline F}
{\mathfrak m}_{\overline K}^{(-r)}/
{\mathfrak m}_{\overline K}^{(-r)+}
)\\
@.=
\displaystyle{{\rm Spec}\ 
S^\bullet
(\Omega^1_{P/S}\otimes_{{\cal O}_P}
\overline F
\otimes_{\overline F}
{\mathfrak m}_{\overline K}^{(-r)}/
{\mathfrak m}_{\overline K}^{(-r)+}
)}\\
@VVV @VVV\\
(P^{(r)}_{S,F})_{\rm red}
@>>>
\displaystyle{{\rm Spec}\ 
\bigoplus_{l\ge 0}
S^{nl}\Omega^1_{P/S}\otimes_{{\cal O}_P}F
\otimes_F
{\mathfrak m}_K^{-ml}/{\mathfrak m}_K^{-ml+1}}.
\end{CD}
\label{eqPrd}
\end{equation}
The horizontal arrows
are isomorphisms
induced by {\rm (\ref{eqPr})}
and
the right vertical arrows
is induced by the natural inclusion.
\end{cor}
{\it Proof.}
The lower horizontal isomorphism
is defined
by ${\cal O}_{P^{(r)}_S}
=
\sum_{l\ge 0}{\mathfrak m}_K^{-[lr]}\cdot {\cal I}_S^l$
in Lemma \ref{lmmn}.1.
The upper isomorphism
is the lower one for the base change
$S'={\rm Spec}\ {\cal O}_{K'}\to S$ such that 
$e_{K'/K}r$ is an integer.
The commutativity of the diagram is clear
from the construction.
\qed

We consider a cartesian 
diagram (\ref{eqTS})
satisfying the condition {\rm (CI)}
in Lemma \ref{lmlci}.
For a positive integers
$m,n>0$,
the diagram
\begin{equation}\begin{CD}
Q^{(m/n)}_{T}
@>>> Q\\
@VVV @VVV\\
P^{(m/n)}_S@>>> P
\end{CD}
\label{eqTSmn}
\end{equation}
is cartesian.
Hence the canonical map
$Q^{(m/n)}_{T}
\to P^{(m/n)}_S$
is also quasi-finite and flat
and induces
a finite map
$Q^{(m/n)}_{T,F}
\to P^{(m/n)}_{S,F}$
on the closed fibers.
For $r=m/n$, we have a quasi-finite morphism
$Q^{(r)}_{T}\to P^{(r)}_S$ of schemes
and a finite morphism of $\widehat {Q^{(r)}_{T}}
\to \widehat {P^{(r)}_S}$
of affine formal schemes
over $\widehat S={\rm Spf}\ \widehat {\cal O}_K$.
If $Q\to P$ is \'etale,
the diagram (\ref{eqTSmn})
with $(m/n)$ replaced by $(r)$
is also cartesian.

A diagram (\ref{eqTS})
satisfying the condition {\rm (CI)}
in Lemma \ref{lmlci}
naturally arises in the following ways.

\begin{eg}
1. Let $A$ be a finite flat 
${\cal O}_K$-algebra locally of complete intersection
and ${\cal O}_K[T_1,\ldots,T_n]/
(f_1,\ldots,f_n)\to A$
be an isomorphism over ${\cal O}_K$.
We define a closed immersion 
$T={\rm Spec}\ A\to 
Q={\rm Spec}\ 
{\cal O}_K[T_1,\ldots,T_n]$
by the surjection
${\cal O}_K[T_1,\ldots,T_n]\to A$.
We also define
a section
$S\to P={\rm Spec}\ 
{\cal O}_K[S_1,\ldots,S_n]$ by
$S_1,\cdots, S_n\mapsto 0$.
Then, by defining $Q\to P$
by $S_i\mapsto f_i$,
we obtain a cartesian diagram
{\rm (\ref{eqTS})}
satisfying the condition {\rm (CI)}
in Lemma \ref{lmlci}
on a neighborhood of $T$.

2. Let $X$ be a smooth scheme
over $k$ and
$D$ be a smooth divisor of $X$.
We consider the 
local ring ${\cal O}_K={\cal O}_{X,\xi}$
at the generic point $\xi$ of $D$. 
Let $f:Y\to X$ be a quasi-finite flat morphism
of smooth schemes over $k$
and assume $V=Y\times_XU\to U=X\setminus D$
is \'etale.
We assume $T=Y\times_XS$
is finite over $S$.
We put $P=X\times_kS$
and $Q=Y\times_kS$
and let $Q\to P$
be $f\times 1_S$.
We consider the immersions
$S\to P$ and $T\to Q$
defined by the natural maps
$S\to X$ and $T\to Y$.
Then we obtain a cartesian diagram
{\rm (\ref{eqTS})}
satisfying the condition {\rm (CI)}
in Lemma \ref{lmlci}.
\end{eg}

\begin{lm}\label{lmbddr}
Let $T$ be a 
finite flat scheme over $S$
of degree $d$
such that $T_K$ is \'etale
over $K$.
We consider
a cartesian diagram
$$
\begin{CD}
T@>>> Q\\
@VVV @VVV\\
S@>>> P
\end{CD}
$$
satisfying the condition 
{\rm (CI)} in Lemma {\rm \ref{lmlci}}.
We consider the following conditions:

{\rm (1)} The ramification of $T$
is bounded by $r$.

{\rm (2)} 
The number of connected components
of the scheme $Q^{(r)}_{T,\overline F}$
is $d$.

{\rm (3)} 
The scheme $Q^{(r)}_{T,\overline F}$
over $P^{(r)}_{S,\overline F}$
is isomorphic to 
the disjoint union of $d$ copies
of $P^{(r)}_{S,\overline F}$.

{\rm (4)} 
The map $Q^{(r)}_{T,\overline F}\to 
P^{(r)}_{S,\overline F}$
is finite and \'etale.

{\rm (5)} 
The induced map
$\overline T\to 
Q^{(r)}_T$
is a closed immersion.

{\rm (6)} The ramification of $T$
is bounded by $r+$.

\noindent
Then, we have implications
$(1)\Leftrightarrow(2)
\Leftrightarrow(3)
\Rightarrow(4)\Rightarrow(5)
\Rightarrow(6)$.
If $Q_K\to P_K$ is finite \'etale,
we have
$(4)\Leftrightarrow(5)$.
\end{lm}
{\it Proof.}
(1)$\Rightarrow$(3)
We may assume that the map $Q_T^{(r)}
\to P_S^{(r)}$
is finite flat of degree $d$
on the generic fiber.
Assume the ramification of $T$ is bounded by $r$.
For each $t\in \overline T_{\overline F}$,
let $Q_{T,\overline F}^{(r),t}$
denote the connected component
containing the image of $t$
by the map $\overline T\to 
Q_T^{(r)}$.
Then,
we have an open and closed immersion
$\coprod_{t\in \overline T_{\overline F}}
Q_{T,\overline F}^{(r),t}\to
Q_{T,\overline F}^{(r)}$.
Since the number of the points
in every geometric fiber of the map
$Q_T^{(r)}
\to P_S^{(r)}$
is at most $d$,
we obtain
an equality
$\coprod_{t\in \overline T_{\overline F}}
Q_{T,\overline F}^{(r),t}=
Q_{T,\overline F}^{(r)}$.

(3)$\Rightarrow$(2)
Clear.

(2)$\Rightarrow$(1)
By \cite{AC} Chap.\ V \S2.4 Theorem 3,
the image of every connected component
of $\overline Q_{T,\overline F}^{(r)}$
is $\overline P_{S,\overline F}^{(r)}$.
Hence the map
$\overline T_{\overline F}
\to \pi_0(\overline P_{S,\overline F}^{(r)})$
is surjective.
Since the cardinalities are the same,
it is a bijection.

(3)$\Rightarrow$(4)
Clear.

(4)$\Rightarrow$(5)
We may assume that the map
$Q_T^{(r)}
\to P_S^{(r)}$
is finite \'etale.
Then, the diagram 
\begin{equation}
\begin{CD}
\overline T@>>> Q_T^{(r)}\\
@VVV @VVV\\
S@>>> P_S^{(r)}
\end{CD}\label{eqcart}
\end{equation}
is cartesian and the upper horizontal
arrow is a closed immersion.

(5)$\Rightarrow$(6)
Let $s>r$ be a rational number.
Then, we have a commutative diagram
$$\begin{CD}
\overline T_{\overline F}@>>> 
(Q_{T,\overline F}^{(s)})_{\rm red}
@>>>
\overline T_{\overline F}@>>> Q_T^{(r)}\\
@VVV @VVV@VVV @VVV\\
{\rm Spec}\ \overline F
@>>> P_{S,{\overline F}}^{(r)}@>>>
{\rm Spec}\ \overline F
@>>> P_S^{(r)}.
\end{CD}$$
Since the composition of
the left two upper horizontal arrows
is the identity,
the map
$\overline T_{\overline F}\to
\pi_0(Q_{T,\overline F}^{(s)})$
is an injection.
\qed

The equivalence 
(1)$\Leftrightarrow$(2) means that
Definition \ref{dfram}.1
is equivalent to that in
\cite{AS} Definition 6.3 if $A$
is locally of complete intersection.
Under the assumption
that $Q_K\to P_K$ is finite \'etale,
we have an equivalence
$(4)\Leftrightarrow(5)
\Leftrightarrow(6)$
(cf.\ \cite{AS2} Corollary 4.12).
The author does not know
how to prove the implication
$(6)\Rightarrow(4)$
without using rigid geometry.

\begin{cor}
Let $T$ be a finite flat
scheme locally of
complete intersection over $S$ and
$T\to P$
a closed immersion 
to a smooth scheme over $S$.
Assume $T_K$ is \'etale over $K$.
Then, there exists
a positive rational number
$r$
such that 
the ramification of
$T$ over $S$ is bounded by $r$.
\end{cor}
{\it Proof.}
By Lemma \ref{lmbddr}
(5)$\Rightarrow$(6),
it is a consequence of
Lemma \ref{lmspl}.
\qed

\subsection{Logarithmic variant}

We keep the notation
in the previous subsection.
We consider a logarithmic variant
of the constructions in
the previous section,
without using log geometry.
We work with Cartier divisors
to replace log structures.

Let $D_S
\subset S={\rm Spec}\ {\cal O}_K$ 
be the Cartier divisor
${\rm Spec}\ F$.
Let $T$ be a flat scheme
of finite type
over $S$ and $D_T$ be a Cartier divisor of
$T$ satisfying the following condition:
\begin{itemize}
\item[(D)] For each $t\in T$, 
there exists an integer
$e_t\ge 1$ such that
the pull-back of $D_S$ is equal to
$e_tD_T$ on a neighborhood of $t$.
\end{itemize}
The condition (D) implies
that the complement $T\setminus D_T$
is equal to the generic fiber $T_K$.
If $P$ is a regular flat scheme
of finite type over $S$
and if the reduced
closed fiber $D_P=
(P\times_SD_S)_{\rm red}$
is regular,
then the Cartier divisor
$D_P$ satisfies the condition (D).
For $(T,D_T)$ satisfying the condition
(D),
let $e_T$ denote
the least common multiple
${\rm lcm}_{t\in T}e_t$.
The condition 
$e_T=1$ is equivalent to
that $D_T$ is the pull-back of $D_S$.

Let $T$ be a flat scheme of finite type
over $S$
and $D_T$ be a Cartier divisor of
$T$ satisfying the condition (D).
For a surjection 
$S'={\rm Spec}\ {\cal O}_{K'}\to S$
of the spectra of discrete valuation rings
of ramification index $e'
=e_{K'/K}$,
we define the log
base change or the log product
$T'=T\times_S^{\log}S'$ as follows.
First, we consider the case
where we have $e_t=e$
for every $t\in T$ and
there exists a generator 
$f$ of the ideal
of $D_T$.
Let $\pi'$ be a prime element of $K'$.
We define 
$u\in \Gamma(T,{\cal O}_T^\times)$
and $v\in {\cal O}_{K'}^\times$ by
$\pi=uf^e$ and $\pi=v\pi^{\prime e'}$
and a morphism
$T\times_SS'\to
{\rm Spec}\ {\mathbb Z}
[X,Y,U^{\pm1},V^{\pm1}]/(UX^e-VY^{e'})$
by $X\mapsto f,Y\mapsto \pi',
U\mapsto u,V\mapsto v$.
Let $d=(e,e')$ be the maximum common divisor
and put $e=de_1$ and $e'=de'_1$.
Let $a$ and $b$ be integers satisfying
$d=ae+be'$.
We define
\begin{equation}
\begin{array}{rcl}
T\times^{\log}_SS'&=&
(T\times_SS')\times_{
{\rm Spec}\ {\mathbb Z}
[X,Y,U^{\pm1},V^{\pm1}]/(UX^e-VY^{e'})}
{\rm Spec}\ {\mathbb Z}
[Z,W^{\pm1},U^{\pm1}]\\
&&\\
&=&(T\times_SS')
[Z,W^{\pm1}]/
(f-Z^{e'_1}W^a,
\pi'-Z^{e_1}W^{-b},
v-uW^d)
\end{array}
\label{eqlp}
\end{equation}
where
${\mathbb Z}
[X,Y,U^{\pm1},V^{\pm1}]/(UX^e-VY^{e'})\to
{\mathbb Z}
[Z,W^{\pm1},U^{\pm1}]$
is defined by
$X\mapsto Z^{e'_1}W^a,$
$Y\mapsto Z^{e_1}W^{-b},
V\mapsto UW^d$.
This is independent of
the choices and is well-defined.
In the general case,
we define 
$T\times_S^{\log}S'$
by patching.

The canonical map
$T'=T\times_S^{\log}S'
\to T\times_SS'$
is finite.
If $e_T=1$,
the canonical map
$T'=T\times_S^{\log}S'
\to T\times_SS'$
is an isomorphism.

If $T'$ is flat over $S'$,
we define a Cartier divisor $D_{T'}$
locally to be that defined by
$Z$ in (\ref{eqlp}).
Then, the divisor $D_{T'}$
satisfies the condition (D)
by putting 
$e_{t'}=e_t/{\rm gcd}(e_t,e_{K'/K})$
for $t'\in T'$ above $t\in T$.
We have $e_{T'}=
e_T/{\rm gcd}(e_T,e_{K'/K})$.
In particular, if $e_{K'/K}$
is divisible by $e_T$,
we have $e_{T'}=1$ and
the divisor $D_{T'}$
is the pull-back of $D_{S'}$.

\begin{df}\label{dflog}
Let $K$ be a discrete valuation field
and let $D_S$ be the Cartier divisor
${\rm Spec}\ F$
of $S={\rm Spec}\ {\cal O}_K$.

1. Let $T$ be a flat scheme 
of finite type over $S$
and $D_T$ be a Cartier divisor
of $T$
satisfying the condition
{\rm (D)}.
We say $(T,D_T)$ is
log flat over $S$,
if, for an arbitrary surjection 
$S'={\rm Spec}\ {\cal O}_{K'}\to S$
of the spectra of discrete valuation rings,
the log base change
$T'=T\times_S^{\log}S'\to S'$ is flat.

2. 
Let $P$ be a regular flat scheme
of finite type over $S$
such that the reduced
closed fiber $D_P=(P\times_SD_S)_{\rm red}$
is regular.
We say $P$ is
log smooth over $S$,
if \'etale locally on $P$, 
there exists a smooth map
$P\to P_e$ for some $e\ge 1$
where $$P_e=
\begin{cases}{\rm Spec}\ {\cal O}_K[t]/(t^e-\pi)
&\quad\text{ if }e\in {\cal O}_K^\times,\\
{\rm Spec}\ {\cal O}_K[t,u^{\pm1}]/(ut^e-\pi)
&\quad\text{ if otherwise.}
\end{cases}$$
and $\pi$ is a prime element of $K$.

3. Let $T\to P$ be 
a closed immersion of
flat schemes
over $S$ and $D_T$ 
and $D_P$ be Cartier divisors
satisfying the condition {\rm (D)}.
If $D_T=D_P\times_PT$,
we say $T\to P$ is
an exact closed immersion.
\end{df}

\begin{lm}\label{lmlsm}
Let $P$ be a regular flat scheme
of finite type over $S$
such that $D_P=
(P\times_SD_S)_{\rm red}$ is regular
and that $P$ is log smooth
over $S$.
Let
$S'={\rm Spec}\ {\cal O}_{K'}\to S$
be a surjection 
of the spectra of discrete valuation rings.
We put $P'=P\times_S^{\log}S'$.

1. The scheme $P'$ is regular and flat over $S'$,
$D_{P'}=(P'\times_{S'}D_{S'})_{\rm red}$
is regular and 
$P'$ is log smooth over $S'$.
If the ramification index
$e'=e_{K'/K}$ is divisible
by $e_P$,
the map $P'\to S'$ is smooth.

2. Let $T$ be a finite flat scheme over
$S$ and $T\to P$
be a regular exact closed immersion.
Then, $T$ is log flat
and $T'=T\times_S^{\log}S'\to P'$
is also a regular exact closed immersion.
\end{lm}
{\it Proof.}
1. It suffice to prove the case where
$P=P_e$ for an integer $e\ge 1$.
If $e$ is invertible in ${\cal O}_K$,
in the notation of (\ref{eqlp}),
the log product
$P_e\times^{\log}_SS'$
is given by
\begin{eqnarray*}
&&{\rm Spec}\ {\cal O}_{K'}[t]/(t^e-\pi)
[Z,W^{\pm1}]/(t-Z^{e'_1}W^a,
\pi'-Z^{e_1}W^{-b},
v-W^d)\\
&=&
{\rm Spec}\ {\cal O}_{K'}
[W,Z]/(W^d-v,Z^{e_1}-W^b\pi').
\end{eqnarray*}
Since
$W^b\pi'$ is a prime element
of the unramifed extension
${\cal O}_{K'}
[W]/(W^d-v)$,
the assertion follows.
Assume $e$ is not
invertible in ${\cal O}_K$.
Then, in the notation of (\ref{eqlp}),
$P_e\times^{\log}_SS'$
is given by
\begin{eqnarray*}
&&
{\rm Spec}\ {\cal O}_{K'}[t,u^{\pm1}]/(ut^e-\pi)
[Z,W^{\pm1}]/(t-Z^{e'_1}W^a,
\pi'-Z^{e_1}W^{-b},
v-uW^d)\\
&=&{\rm Spec}\ {\cal O}_{K'}
[Z,W^{\pm1}]/(W^{-b}Z^{e_1}-\pi').
\end{eqnarray*}
First, we consider
the case where $e_1$ is invertible in ${\cal O}_K$.
In this case, the \'etale covering
$P_e\times^{\log}_SS'[V]/(V^{e_1}-W)
=
{\rm Spec}\ {\cal O}_{K'}
[Z,V^{\pm1}]/((V^{-b}Z)^{e_1}-\pi')$
of
$P_e\times^{\log}_SS'$
is smooth over
$P'_{e_1}={\cal O}_{K'}
[T]/(T^{e_1}-\pi')$.
Assume $e_1$ is not
invertible in ${\cal O}_K$.
Then, by the definition of $b$,
we have $(b,e_1)=1$ and
$b$ is invertible in ${\cal O}_K$.
Hence 
$P_e\times^{\log}_SS'=
{\rm Spec}\ {\cal O}_{K'}
[Z,W^{\pm1}]/(W^{-b}Z^{e_1}-\pi')$
is \'etale over
$P'_{e_1}={\rm Spec}\ {\cal O}_{K'}
[Z,V^{\pm1}]/(VZ^{e_1}-\pi')$.

If $e'$ divides $e$,
we have $e_1=1$
and $P_e\times^{\log}_SS'$
is smooth over $P'_1=S'$.

2.
By the definition of the base change,
the map $T'\to P'$ is a closed immersion
and $T'$ is finite over $S'$.
Since the ideal ${\cal I}_{T'}\subset
{\cal O}_{P'}$ is locally generated
by $d$ elements where
$d$ is the relative dimension of $P'$
over $S'$,
the immersion $T'\to P'$
is regular and $T'$ is flat over $S'$.
By the definition of $D_{T'}$,
the immersion $T'\to P'$
is regular.
\qed

\begin{cor}\label{cortb}
Let $P$ be a regular flat scheme 
of finite type over $S$
such that $D_P=
(P\times_SD_S)_{\rm red}$ is 
irreducible and regular
and that $P$ is log smooth
over $S$.
Let 
$S'={\rm Spec}\ {\cal O}_{K'}
\to P$ be the localization
at the generic point $\xi$ of $D_P$.

Let $L$ be
a finite separable extension
of $K$,
$T={\rm Spec}\ {\cal O}_L$
and $D_T=
(T\times_SD_S)_{\rm red}$.
Then,
$T'=T\times_S^{\log}S'$
is equal to
${\rm Spec}\ {\cal O}_{L\otimes_KK'}$
and we have
$D_{T'}=
(T'\times_{S'}D_{S'})_{\rm red}$.
\end{cor}
{\it Proof.}
It suffices to show
that $T'=T\times_S^{\log}S'$
is regular
and that
$D_{T'}$ is defined by
the reduced closed point
at each closed point $\xi'\in T'$.
Let $t\in T$ be the image
of $\xi'$ and
$T_t$ be the localization at $t$.
Then, the localization
of $T'$ at $\xi'$
is equal to
a localization of
$P\times_S^{\log}T_t$
and the assertion follows
from Lemma \ref{lmlsm}.1
\qed

For the convenience of a reader
familiar with 
the terminologies on log geometry
as in \cite{KS0} \S4,
we include a lemma,
not used in the sequel,
showing that 
the Definition \ref{dflog} above is
a special case of the standard definitions.

\begin{lm}\label{lmsyn}
We consider
$S={\rm Spec}\ {\cal O}_K$
as a log scheme with 
the log structure
defined by $D_S$.

1. 
Let $T$ be a
flat scheme of finite type
over $S$ and $D_T$
be a Cartier divisor
satisfying the condition {\rm (D)}.
Then, the following conditions
are equivalent:

{\rm (1)} The log scheme $T$
with the log structure defined
by $D_T$ is log flat
over $S$.

{\rm (2)} $(T,D_T)$
is log flat over $S$
in the sense of Definition {\rm \ref{dflog}.1}.

2. 
Let $P$ be a regular flat scheme
of finite type over $S$
such that the reduced
closed fiber $D_P=
(P\times_SD_S)_{\rm red}$
is regular.
Then, the following conditions
are equivalent.

{\rm (1)}  The log scheme $P$
with the log structure defined
by $D_P$ is log smooth
over $S$.

{\rm (2)} $(P,D_P)$
is log smooth over $S$
in the sense of Definition {\rm \ref{dflog}.2}.
\end{lm}
{\it Proof.}
1.
$(1)\Rightarrow(2)$
Let
$S'={\rm Spec}\ {\cal O}_{K'}\to S$
be a surjection 
of the spectra of discrete valuation rings
and we show that
the base change
$T'=T\times_S^{\log}S'\to S'$ is flat
at each closed point
$t'\in T'$.
We put $e'=e_{t'}$.
Let $S'_1$ be the localization
of $P_{e'}$ over $S'$
and consider the cartesian diagram
\begin{equation}
\begin{CD}
T'@<<< T'_1=T'\times_{S'}^{\log}S'_1\\
@VVV @VVV\\
S'@<<< S'_1.
\end{CD}
\label{eqfl}
\end{equation}
Since $e'=e_{t'}$,
the map $T'_1\to S'_1$ is strict
on a neighborhood $V'_1$
of the inverse image of $t'$.
Since $T'\to S'$ is log flat,
the map $V'_1\to S'_1$ is log flat
and strict ane hence is flat.
Since $S'_1\to S$ is log flat,
the map $V'_1\to T'$ is also log flat
and strict ane hence is flat.
Hence the map $T'\to S'$ is flat.

$(2)\Rightarrow(1)$
Let $t\in T$ be a closed point
and put $e=e_t$.
Let $S_1$ be the localization
of $P_e$
and consider the cartesian diagram
(\ref{eqfl})
with $'$ removed everywhere.
Then, as above,
there exists an open neighborhood
$V_1\subset T_1$ of the inverse image
of $t$ such that
$V_1\to T$ and $V_1\to S_1$
are flat.
Hence by \cite{KS0} Proposition 4.3.10,
the map $T\to S$ is log flat.

2. $(2)\Rightarrow(1)$
Since $P_e$ is log smooth over $S$,
the assertion follows.

$(1)\Rightarrow(2)$
We consider the ring homomorphism
${\mathbb Z}[{\mathbb N}]
\to {\cal O}_K$
sending 
$1\in {\mathbb N}$
to a prime element $\pi$.
The question is \'etale local.
Hence, we may assume
that $P={\rm Spec}\
{\cal O}_K\otimes_{{\mathbb Z}[{\mathbb N}]}
{\mathbb Z}[M]$
for a morphism 
${\mathbb N}\to M$ of fs-monoids
such that
the map ${\mathbb Z}\to M^{\rm gp}$
is an injection and that
the order of the torsion part
of the cokernel
is invertible in ${\cal O}_K$.
Further 
$\overline M=M/M^{\times}$
is isomorphic to ${\mathbb N}$.
We may assume $M^{\rm gp}$
is torsion free.

If the order $e$ of the cokernel
of ${\mathbb Z}\to 
\overline M^{\rm gp}$
is invertible in ${\cal O}_K$,
we may assume
$M={\mathbb N}$.
In this case,
we have $P=P_e$.
Assume $e$ is not invertible.
In this case,
we may assume
$M={\mathbb N}\times 
{\mathbb Z}$
and the map
${\mathbb N}\to M$
sends $1$ to $(e,1)$.
Then, we also have
$P=P_e$.
\qed

In this subsection,
$T$ denotes a finite flat 
scheme over $S$
and $D_T$ denotes
a Cartier divisor of $T$
satisfying the condition (D)
such that $(T,D_T)$
is log flat over $S$.
Recall that $e_T$
denotes the least common multiple
of the integers $e_t\ge 1$ for
closed points $t\in T$.

\begin{df}\label{dflram}
Let $T$ be a finite flat 
scheme over $S$
such that $T_K$ is \'etale over $K$
and let $D_T$ denote
a Cartier divisor of $T$
satisfying the condition {\rm (D)}
such that $(T,D_T)$
is log flat over $S$.

1. For a rational number $r>0$,
we say that the log ramification of $(T,D_T)$
over $S$
is bounded by $r$
if, for a surjection
$S'={\rm Spec}\ {\cal O}_{K'}\to S$
of spectra of discrete valuation rings
such that $e=e_{K'/K}$
is divisible by $e_T$,
the ramification of
the finite flat scheme
$T\times_S^{\log}S'$ over $S'$ is bounded by $er$.

2. For a rational number $r\ge 0$,
we say that the log ramification of $(T,D_T)$
over $S$
is bounded by $r+$
if the log ramification of
$(T,D_T)$ over $S$ is bounded by $s$
for every rational number $s>r$.

3.
Let $L$ be a finite \'etale $K$-algebra,
$T={\rm Spec}\ {\cal O}_L$
and $D_T=
{\rm Spec}\ ({\cal O}_L\otimes
_{{\cal O}_K}F)_{\rm red}$.
Then, we say that the log ramification of $L$
over $K$ is bounded by $r$
(resp.\ by $r+$)
if the log ramification of $(T,D_T)$
is bounded by $r$
(resp.\ by $r+$).
\end{df}

Let $(T,D_T)$ be as in
Definition \ref{dflram}
and $S'\to S$
be a surjection
of spectra of
discrete valuation rings
of ramification index $e$.
Then, it is clear from
the definition that
the log ramification of
$T$ over $S$ is bounded by
$r$ if and only if
the ramification of
$T\times^{\log}_SS'$ over $S'$ is bounded by
$er$.

\begin{lm}\label{lmtb}
Let $P$ be a regular flat scheme 
of finite type over $S$
such that $D_P=
(P\times_SD_S)_{\rm red}$ is 
irreducible and regular
and that $P$ is log smooth
over $S$
and let $\xi$ be
the generic point of $D_P$.
We put
${\cal O}_{K'}=
{\cal O}_{P,\xi}$
and consider
the surjection
$S'={\rm Spec}\ {\cal O}_{K'}\to S$
of ramification index $e$.

Then, 
for a finite separable extension $L$
of $K$,
the log ramification of
$L$ over $K$ is bounded by
$r$ if and only 
the log ramification of
$L\otimes_KK'$ over $K'$ is bounded by
$er$.
\end{lm}
{\it Proof.}
Clear from Corollary \ref{cortb}
and the above remark.
\qed

Let $T$ be a finite flat
scheme over $S$ and
$D_T$ be a Cartier divisor
satisfying the condition (D).
We consider an exact closed immersion
$T\to P$
to a log smooth scheme $P$ over $S$.
Let $S'={\rm Spec}\ {\cal O}_{K'}\to S$
be a surjection of spectra
of discrete valuation rings
of ramification index $e$.
Then, the base change
$T'=T\times^{\log}_SS'\to P'=P\times^{\log}_SS'$
is an exact closed immersion
to a log smooth scheme
over $S'$.
Assume $e$ is divisible by the integer
$e_P$.
Then, the map $P'\to S'$
is smooth.
Thus for positive integers
$m,n>0$ and $r=m/n$,
we apply the construction in \S1.1
to define
$P^{\prime [em/n]}_{T'},
P^{\prime (em/n)}_{T'},
P^{\prime [er]}_{T'},
P^{\prime (er)}_{T'},
\widehat{P^{\prime (er)}_{T'}}$
etc.

\begin{eg}
Assume $K$ is complete.

1.
Let $L$ be a finite separable extension of $K$
and ${\cal O}_K[X_1,\ldots,X_n]
/(f_1,\ldots,f_n)\to {\cal O}_L$
be an isomorphism.
Let $m\le n$ be an integer
such that the images $z_1,\ldots,z_m$ 
of $X_1,\ldots,X_m$
are non-zero and that
$z_i$ is a prime element of $L$
for some $1\le i\le m$.
We define a map
${\mathbb N}^{m+1}\to {\mathbb N}$
by sending the standard basis
of ${\mathbb N}^{m+1}$ 
to $e_{L/K},v_L(z_1),\cdots,$
$v_L(z_m)$.
Let $M$ be the inverse image
of ${\mathbb N}$ by the induced map
${\mathbb Z}^{m+1}\to {\mathbb Z}$.
We define ${\mathbb N}^{m+1}
\to {\cal O}_K[X_1,\ldots,X_n]$
by sending the standard basis
to $\pi,X_1,\ldots,X_m$
where $\pi$ is a prime element of $K$.

We put
$P={\rm Spec}\ {\cal O}_K[X_1,\ldots,X_n]
\otimes_{{\mathbb Z}[{\mathbb N}^{m+1}]}
{\mathbb Z}[M]$.
Then, $P$ is regular,
the reduced closed fiber
of $P$ is regular
and $P$ is log smooth over $S$.
Further the isomorphism
${\cal O}_K[X_1,\ldots,X_n]
/(f_1,\ldots,f_n)\to {\cal O}_L$
induces an exact closed immersion
$T={\rm Spec}\ {\cal O}_L\to P$.
For a finite extension 
$K'$ over $K$ with ramification
index $e$ divisible by $e_{L/K}$,
the affinoid variety over $K'$
defined by the formal ${\cal O}_{K'}$-scheme
$\widehat {P^{\prime (er)}_{T'}}$
is the 
affinoid variety $Y_{Z,I,P}^a$
defined in \cite{AS} 3.1 for $a=r$
and $I=\{1,\ldots,n\}
\supset P=\{1,\ldots,m\}$.

2.
Assume ${\rm Spf}\ {\mathbf A}$ is the completion
of $P$ at $T={\rm Spec}\ A$.
For a finite extension 
$K'$ over $K$ with ramification
index $e$ divisible by $e_{L/K}$,
the affinoid $K'$-variety
defined by the formal ${\cal O}_{K'}$-scheme
$\widehat {P^{\prime (er)}_{T'}}$
is the 
affinoid variety $X^j_{\log}({\mathbf A}\to A)_K'$
defined in \cite{AS2} \S4.2 for $j=r$.
\end{eg}

We consider a cartesian diagram
\begin{equation}\begin{CD}
T@>>> Q\\
@VVV @VVV\\
S@>>> P
\end{CD}
\label{eqlTS}
\end{equation}
of schemes over $S$
satisfying the following condition:

(LCI)
The vertical arrows
are quasi-finite and flat
and 
the horizontal arrows
are closed immersions.
The scheme $P$ is smooth over $S$,
$Q$ is regular flat
over $S$,
$D_Q=(Q\times_SD_S)_{\rm red}$
is smooth over $F$
and $Q$ is log smooth over $S$.

\noindent
We consider the Cartier divisor
$D_T=D_Q\times_QT$.
Then, by Lemma \ref{lmlsm}.2,
the pair $(T,D_T)$ is log flat over $S$.

Let $S'\to S$ be a surjection
of the spectra of
discrete valuation rings
of ramification index $e$.
We assume that
$e_Q$ divides $e$.
Then by Lemma \ref{lmlsm},
the log product
$Q'=Q\times_S^{\log}S'$
is smooth 
and
the immersion
$T'=T\times_S^{\log}S'
\to Q'$
is a regular immersion.
Hence the cartersian diagram
\begin{equation}\begin{CD}
T'@>>> Q'\\
@VVV @VVV\\
S'@>>> P'&=P\times_SS'
\end{CD}
\end{equation}
satisfies the condition
(CI) in Lemma  \ref{lmlci}.

\subsection{Logarithmic ramification groups}

In \cite{AS} Definitions 3.4 and 3.12,
we introduced two filtrations,
the non-logarithmic one and
the logrithmic one,
by ramification groups
of the absolute Galois group.
In this paper, we will be only interested
in the logarithmic filtration.

Assume $K$ is a henselian discrete valuation field.
Let $\overline K$
be a separable closure
and 
$G_K={\rm Gal}(\overline K/K)$
be the absolute Galois group.
In \cite{AS} Definition 3.12,
we define a decreasing filtration
by logarithmic ramification
groups $G_{K,\log}^r\subset G_K$
indexed by positive rational numbers
$r\ge 0$.
We put
$G_{K,\log}^{r+}=
\overline{
\bigcup_{q>r}G_{K,\log}^q}$
and
${\rm Gr}^r_{\log}G_K=
G_{K,\log}^r/
G_{K,\log}^{r+}$.
For $r=0$,
$G_{K,\log}^{0+}
\subset G_{K,\log}^0$
are equal to
the inertia subgroup and
its pro-$p$
Sylow subgroup $P\subset I$.

We consider the 
opposite category $({\rm FE}/K)$ of
finite \'etale $K$-algebras.
We identify
the category $({\rm FE}/K)$
with that of finite discrete sets with
continuous action of 
the absolute Galois group $G_K$
by the fiber functor
$X\mapsto X(\overline K)$.
For a rational number $r\ge 0$,
the \'etale $K$-algebras $L$
such that the log ramification
is bounded by $r+$ form
a Galois subcategory $({\rm FE}/K)^{r+}$
of $({\rm FE}/K)$
corresponding to
a normal closed subgroup $G_{K,\log}^{r+}
\subset G_K={\rm Gal}(\overline K/K)$.
For an extension of
discrete valuation field
$K'$ over $K$ of ramification
index $e$,
the natural map
$G_{K_1}\to G_K$
sends $G_{K',\log}^{er}$
into $G_{K,\log}^r$.

In the rest of this section, we 
assume that $K$
satisfies the following condition:

\smallskip
(Geom) There exist
a smooth scheme $X$ over $k$,
an irreducible divisor $D$ smooth over $k$
with the generic point $\xi$
and an isomorphism $S\to {\rm Spec}\ O^h_{X,\xi}$
to the henselization of the local ring.

\smallskip
Let $\Omega^1_F(\log)$
denote the $F$-vector space
$\Omega^1_{X/k}(\log D)_{\xi}\otimes_{{\cal O}_{X,\xi}}F$.
It fits in an exact sequence
$0\to \Omega^1_{F/k}\to
\Omega^1_F(\log)\overset{\rm res}
\to F\to 0$.
We extend the normalized discrete valuation
$v_K:K \to {\mathbb Z}\cup \{\infty\}$ to
$v_K:\overline K \to {\mathbb Q}\cup \{\infty\}$.
Let $r>0$ be a rational number.
We put ${\mathfrak m}_{\overline K}^r=
\{a\in \overline K|v_K(a)\ge r\}$
and ${\mathfrak m}_{\overline K}^{r+}=
\{a\in \overline K|v_K(a)> r\}$.
Let $\Theta^{(r)}_{\log}
=\Theta^{(r)}_{F,\log}$
denote the $\overline F$-vector space
${\mathbf V}(
\Omega^1_F(\log)
\otimes_F
{\mathfrak m}_{\overline K}^{(-r)}
/{\mathfrak m}_{\overline K}^{(-r)+}
)$.

Let $P'=(X\times_kS)'$
be the blow-up of
$X\times_kS$
at $D\times_kD_S$
and define the
log product 
$P=(X\times_kS)^\sim
\subset P'$
to be the complement
of the proper transforms
of $D\times_kS$
and of $X\times_kD_S$.
Then, $P$ is smooth over $S$
and the canonical
map $S\to X$ induces
a section $S\to P$.
Thus, for a rational number
$r>0$,
applying the construction in \S1.1,
we define the schemes
$P_S^{(r)},
\overline P_{S,{\overline F}}^{(r)}$,
etc.
Since ${\cal N}_{S/P}
=\Omega^1_{X/k}(\log D)_{\xi}$,
we have a canonical isomorphism
\begin{equation}
\begin{CD}
\overline P_{S,{\overline F}}^{(r)}
@>>>
\Theta_{\log}^{(r)}
\end{CD}
\label{eqtheta}
\end{equation}
by Lemma \ref{lmmn}.

Under the condition (Geom),
a canonical surjection
$\pi_1^{\rm ab}
(\Theta^{(r)}_{\log})
\to {\rm Gr}_{\log}^rG_K$
is defined in \cite{AS2} (5.12.1).
We recall the construction.
Let $L$ be a finite \'etale algebra over $K$.
After replacing $X$ by an \'etale neighborhood of $\xi$
if necessary,
there exists a finite flat morphism
$f:Y\to X$ of smooth schemes over $k$
such that $V=Y\times_XU\to U=X\setminus D$
is \'etale and that $Y\times_XS=T={\rm Spec}\ {\cal O}_L$.
We also assume that
$V\subset Y$ is the complement of
a smooth divisor $E$.

Similarly
as the construction of
$P=(X\times_kS)^\sim$,
let $Q'=(Y\times_kS)'$
be the blow-up of
$Y\times_kS$
at $E\times_kD_S$
and $Q=(Y\times_kS)^\sim
\subset Q'$
be the complement
of the proper transforms
of $E\times_kS$
and of $Y\times_kD_S$.
We consider the immersions
$S\to P$ and $T\to Q$
defined by the natural maps
$S\to X$ and $T\to Y$.
Then we obtain a cartesian diagram
{\rm (\ref{eqlTS})}
satisfying the condition
(LCI).

Let $K'$ be
a finite extension
such that the ramification
index $e'$ is divisible by $e_{L/K}$.
We put
$S'={\rm Spec}\ {\cal O}_{K'}$
and consider the diagram
$$\begin{CD}
T'&=T\times_S^{\log}S'
@>>>
Q'&=Q\times_S^{\log}S'\\
@VVV& @VVV\\
S'&
@>>>
P'&=P\times_SS'
\end{CD}$$
satisfying the condition
(CI).
Assume that the log ramification
of $L$ over $K$
is bounded by $r+$.
Then, 
the conditions (4) and (6)
in Lemma \ref{lmbddr}
are equivalent
in this case and
the induced map
\begin{equation}
\begin{CD}
\overline Q_{T',\overline F}^{\prime(er)}
@>>>
\overline P_{S',\overline F}^{\prime(er)}
=\Theta^{(r)}_{F,\log}
\end{CD}
\label{eqcvL}
\end{equation}
is finite \'etale.
This construction defines a functor
$({\rm FE}/K)^{r+}
\to ({\rm FE}/\Theta^{(r)}_{\log})$
to the category
of finite \'etale schemes
over $\Theta^{(r)}_{\log}$
and hence a morphism
$\pi_1(\Theta^{(r)}_{\log})
\to G_K/G_{K,\log}^{r+}$.
In \cite{AS2} Theorem 5.12.1,
it is proved that
it factors through the abelian quotient
and induces a surjection
\begin{equation}
\begin{CD}
\pi_1^{\rm ab}
(\Theta^{(r)}_{\log})
@>>> {\rm Gr}_{\log}^rG_K.
\end{CD}
\label{eqGr}
\end{equation}

We give a compatibility
of the map (\ref{eqGr})
with a log smooth base change.
Let $S\to X$ be as above.
Let $t$ be a uniformizer of $D
\subset X$
and $e_1\ge 1$ be an integer.
Let $X_1$ be a scheme
smooth over
$$Z_{e_1}=
\begin{cases}
X[T]/(T^{e_1}-t)
&\quad\text{ if $e_1$
is invertible in $k$,}\\
X[T,U^{\pm 1}]/(UT^{e_1}-t)
&\quad\text{ if $e_1$
is $0$ in $k$}
\end{cases}$$
and assume
$D_1=(D\times_XX_1)_{\rm red}$
is irreducible.
Let ${\cal O}_{K_1}$
be the henselization
${\cal O}_{X_1,\xi_1}^h$
at the generic point $\xi_1$
of $D_1$.

\begin{lm}\label{lmr}
Let $S_1={\rm Spec}\ {\cal O}_{K_1}
\to S={\rm Spec}\ {\cal O}_K$
be the surjection
of the spectra of discrete valuation
rings of ramification index $e_1$
above and 
let $r>0$ be a rational number.
Let $F_1$ denote the residue field of $K_1$
and let $\pi:\Theta_{F_1,\log}^{(e_1r)}
\to 
\Theta_{F,\log}^{(r)}$
be the map induced by
$F_1\otimes_F \Omega_F(\log)
\to \Omega_{F_1}(\log)$.

Then, the induced map 
${\rm Gr}_{\log}^{er}G_{K_1}
\to 
{\rm Gr}_{\log}^rG_K$
is a surjection
and the diagram
\begin{equation}
\begin{CD}
\pi_1^{\rm ab}
(\Theta_{F_1,\log}^{(e_1r)})
@>{\pi_*}>>
\pi_1^{\rm ab}
(\Theta_{F,\log}^{(r)})\\
@VVV @VVV\\
{\rm Gr}_{\log}^{e_1r}G_{K_1}
@>>>
{\rm Gr}_{\log}^rG_K
\end{CD}
\label{eqtb}
\end{equation}
is commutative.
\end{lm}
{\it Proof.}
The natural map
$F_1\otimes_F \Omega_F(\log)
\to \Omega_{F_1}(\log)$
is injective
and hence
$\pi:\Theta_{F_1,\log}^{(e_1r)}
\to 
\Theta_{F,\log}^{(r)}
\times_{\overline F}
\overline F_1$
is a
surjection of
vector spaces.
Thus,
$\pi_*:
\pi_1^{\rm ab}
(\Theta_{F_1,\log}^{(e_1r)})
\to
\pi_1^{\rm ab}
(\Theta_{F,\log}^{(r)})$
is a surjection.
Hence,
it suffices to show
the commutativity of
the diagram (\ref{eqtb}).

Let $Y\to X$
and $Q=(Y\times S)^\sim
\to 
P=(X\times S)^\sim$
be finite coverings
and $S'\to S$
be a finite surjection
appeared in 
the construction of 
the map (\ref{eqcvL}).
The normalization
$Y_1$ of the fiber product
$Y\times_XX_1$
is smooth over $k$
and $V_1=V\times_UU_1$
is the complement
of a smooth divisor
$D_{Y_1}\subset Y_1$.
By Corollary \ref{cortb},
the log product
$S'\times_S^{\log}S_1$
is normal and is
a finite disjoint union
of spectra of discrete 
valuation rings.
Let $S'_1={\rm Spec}\
{\cal O}_{K'_1}$ be a connected
component
and $e'$ be the
ramification index $e_{K'_1/K}$.
Applying the construction of 
the map (\ref{eqcvL})
to $Y_1\to X_1$
and $S'_1\to S_1$,
we obtain
a finite \'etale covering
\begin{equation}
\begin{CD}
\overline Q_{T'_1,\overline F_1}^{\prime(e'r)}
@>>>
\overline P_{S'_1,\overline F_1}^{\prime(e'r)}
=\Theta^{(e'r)}_{F_1,\log}.
\end{CD}
\label{eqcvL1}
\end{equation}
It suffices to show
that the diagram
$$\begin{CD}
\overline Q_{T'_1,\overline F_1}^{\prime(e'r)}
@>>>
\overline P_{S'_1,\overline F_1}^{\prime(e'r)}
\\
@VVV @VVV\\
\overline Q_{T',\overline F}^{\prime(er)}
@>>>
\overline P_{S',\overline F}^{\prime(er)}
\end{CD}$$
is cartesian.

By the construction,
it suffices to show
that the map
$P_1\to 
P\times_SS_1$
is smooth.
Since
$P_1=(X_1\times S_1)^\sim
\to (Z_{e_1}\times S_1)^\sim$
is smooth,
it is reduced
to the case
where $X_1=Z_{e_1}$.
First, we consider
the case
$e_1$ is invertible in $k$.
Let $\pi\in {\cal O}_K$ be
the image of $t$
and $\pi_1\in {\cal O}_{K_1}$ be
the image of $T$.
Then, 
$P_1
=(X[T]/(T^{e_1}-t)\times S_1)^\sim
=X\times S_1[T]/(T^{e_1}-t)
[V^{\pm1}](T-V\pi_1)$
equals $X\times S_1
[V^{\pm1}]/(V^{e_1}\pi-t)$.
This is \'etale over
$P\times_SS_1=
(X\times S)^\sim\times_SS_1
=X\times S_1[W^{\pm1}]/(t-W\pi)$.

We assume
$e_1$ is not invertible in $k$.
Let $\pi\in {\cal O}_K$ be
the image of $t$
and let $\pi_1,u\in {\cal O}_{K_1}$ be
the image of $T,U$.
Then, 
$P_1
=(X[T,U^{\pm1}]/(UT^{e_1}-t)\times S_1)^\sim$
equals
$X\times S_1[T,U^{\pm1}]/(UT^{e_1}-t)
[V^{\pm1}](T-V\pi_1)
=X\times S_1
[U^{\pm1},V^{\pm1}]/(UV^{e_1}\pi_1^{e_1}-t)$.
This is smooth over
$P\times_SS_1
=X\times S_1[W^{\pm1}]/(t-W\pi)$.
\qed

For an $\overline F$-vector
space $V$ of finite dimension,
we introduce a quotient
$\pi_1^{\rm alg}(V)$
of 
$\pi_1^{\rm ab}(V)$
annihilated by $p$.
We regard $V$ as a smooth group scheme 
${\rm Spec}\ S^\bullet V^\vee$
over $\overline F$.
Let $({\rm FE}/V)^{\rm alg}$
be the full subcategory 
of $({\rm FE}/V)$
whose objects are finite \'etale
morphisms $f:X\to V$ such that
there exists a structure of algebraic group scheme
on $X$ and that $f$ is a morphism of algebraic groups.
Let $\pi_1^{\rm alg}(V)$
be the quotient of
$\pi_1^{\rm ab}(V)$
corresponding to the subcategoy
$({\rm FE}/V)^{\rm alg}$.
The pro-finite group
$\pi_1^{\rm alg}(V)$
is the Pontrjagin dual
of the extension group
${\rm Ext}(V,{\mathbb F}_p)$
in the category of
smooth algebraic groups over $
\overline F$.
The map
$V^\vee=
{\rm Hom}_{\overline F}(V,\overline F)
\to {\rm Ext}(V,{\mathbb F}_p)$
sending  
a linear form $f:V\to 
{\mathbb A}^1_{\overline F}$
to the pull-back by $f$ of the 
Artin-Schreier sequence
$0\to {\mathbb F}_p
\to {\mathbb A}^1_{\overline F}
\overset{t\mapsto t^p-t}
\to {\mathbb A}^1_{\overline F}\to 0$
is an isomorphism.
Thus we have defined a canonical isomorphism
\begin{equation}
\begin{CD}
V^\vee
@>>>
{\rm Hom}(\pi_1^{\rm alg}(V),
{\mathbb Q}/{\mathbb Z})
={\rm Ext}(V,
{\mathbb F}_p).
\end{CD}
\label{eqExt}
\end{equation}

\begin{lm}\label{lmisog}
Let $V$ be an $\overline F$-vector space
of finite dimension.
For a continuous character
$\chi:\pi_1^{\rm ab}(V)
\to {\mathbb Q}/{\mathbb Z}$
of finite order,
the following conditions are equivalent.

{\rm (1)} $\chi$ factors through the quotient
$\pi_1^{\rm alg}(V)$.

{\rm (2)} $-^*\chi={\rm pr}_1^*\chi-
{\rm pr}_2^*\chi$
in ${\rm Hom}(\pi_1^{\rm ab}(V\times V)
,{\mathbb Q}/{\mathbb Z})$.
\end{lm}
{\it Proof.}
The implication
(1)$\Rightarrow$(2) is clear.
We show (2)$\Rightarrow$(1).
Let $\chi:\pi_1^{\rm ab}(V)
\to {\mathbb Q}/{\mathbb Z}$
be a character satisfying
$-^*\chi={\rm pr}_1^*\chi-
{\rm pr}_2^*\chi$.
Taking the pull-back
by $i_2:V\to V\times V$,
we obtain
$(-1)^*\chi=-\chi$.
Hence we have 
$+^*\chi={\rm pr}_1^*\chi+
{\rm pr}_2^*\chi$.
By induction on $n$,
we have
$n\cdot \chi=[n]^*\chi$.
Hence, we have
$p\cdot \chi=0$.

Let $f:X\to V$
be the ${\mathbb F}_p$-torsor
corresponding to $\chi$.
By $+^*\chi={\rm pr}_1^*\chi+
{\rm pr}_2^*\chi$, 
we have an isomorphism
$(X\times X)/{\mathbb F}_p
\to X\times_V(V\times V)$
of ${\mathbb F}_p$-torsors
on $V\times V$.
We consider the composition
$\tilde +:X\times X\to
(X\times X)/{\mathbb F}_p
\to X\times_V(V\times V)
\to X$.
Take a point $\tilde 0\in f^{-1}(0)$.
By shifting by
the ${\mathbb F}_p$-action,
we may assume $\tilde 0\tilde +\tilde 0=
\tilde 0$.
Then we can easily verify
that $\tilde +$
defines a group structure on $X$
and the map $f:X\to V$
is compatible with the group structure.
\qed

We will prove the following theorem
in the next subsection.

\begin{thm}\label{thmgr}
Let $K$ be a henselian
discrete valuation field
satisfying the condition
{\rm (Geom)}.
The graded quotient
${\rm Gr}_{\log}^rG_K$
is annihilated by $p$
and the surjection
{\rm (\ref{eqGr})}
induces a surjection
\begin{equation}
\begin{CD}
\pi_1^{\rm alg}
(\Theta^{(r)}_{\log})
@>>> {\rm Gr}_{\log}^rG_K.
\end{CD}
\label{eqisog}
\end{equation}
\end{thm}

By the isomorphism (\ref{eqExt}),
Theorem \ref{thmgr}
has the following corollary.

\begin{cor}\label{corgr}
The dual of the surjection
$\pi_1^{\rm ab}
(\Theta^{(r)}_{\log})
\to {\rm Gr}_{\log}^rG_K$
defines an injection
$$\begin{CD}
{\rm rsw}:{\rm Hom}(
{\rm Gr}_{\log}^rG_K,
{\mathbb F}_p)
@>>>
\Omega^1_F(\log)\otimes_F
{\mathfrak m}^{(-r)}_{\overline K}/
{\mathfrak m}^{(-r)+}_{\overline K}.
\end{CD}$$
\end{cor}

Theorem \ref{thmgr}
implies the prime-to-$p$ part
of the Hasse-Arf theorem.
Let $V$ be an
$\ell$-adic representation $V$
of $G_K$.
Since $P=G_{K,\log}^{0+}$
is a pro-$p$ group,
there exists a unique
direct sum decomposition
$V=\bigoplus_{q\ge0,
q\in {\mathbb Q}}
V^{(q)}$
by sub $G_K$-modules
such that
the $G_{K,\log}^{r+}$-fixed part
is given by
$V^{G_{K,\log}^{r+}}=
\bigoplus_{q\ge r}
V^{(q)}.$
We put ${\rm Sw}_KV=
\sum_rr\cdot {\rm rank}\ V^{(r)}
\in {\mathbb Q}$.

\begin{cor}\label{corHA}
$${\rm Sw}_KV\in {\mathbb Z}[\frac1p].$$
\end{cor}
{\it Proof.}
It suffices to show
that $\dim V\cdot r
\in {\mathbb Z}[\frac1p]$
assuming
$V=V^{(r)}$.
This is equivalent to
that $\dim V$
is divisible by
the prime-to-$p$ part $m$
of the denominator of $r$.
Let $\chi:
{\rm Gr}^r_{\log}G_K
\to \mu_p
\subset \overline{{\mathbb Q}_\ell}^\times$
be a character
appearing in the restriction 
of $V$.
The injection
${\rm Hom}(
{\rm Gr}^r_{\log}G_K
,{\mathbb F}_p)
\to
{\rm Hom}_{\overline F}
({\mathfrak m}_{\overline K}^r/
{\mathfrak m}_{\overline K}^{r+}
,\Omega^1_F(\log)\otimes
\overline F)$
is compatible
with the action of
$I\subset G_K$
and the action
of $I$ is
by the multiplication
through the quotient
$I\to \mu_m$.
Hence there are $m$ conjugates
of $\chi$
appearing 
with the same multiplicities in $V$.
Thus the assertion follows.
\qed

By the same limit argument as in
the proof of Theorem  2.15 \cite{AS2},
Theorem \ref{thmgr} implies
the following.

\begin{cor}
For an arbitrary henselian
discrete valuation field $K$
of characteristic $p>0$,
the pro-finite abelian group
${\rm Gr}^r_{\log}G_K$
is annihilated by $p$.
\end{cor}

\subsection{Nearby cycles}

Let $X$ be a smooth scheme over $k$,
$D$ be a smooth irreducible divisor of $X$
and $U=X\setminus D$
be the complement.
Let $\xi$ be the generic point of $D$
and ${\cal O}_K={\cal O}_{X,\xi}^h$ be the 
henselization of the local ring at $\xi$.
We put $S={\rm Spec}\ {\cal O}_K$
and let $\eta={\rm Spec}\ K$ be
the generic point.
We consider the log product
$P=(X\times_kS)^\sim$ as
in the last subsection
and the section
$S\to P$
induced by the canonical map
$S\to X$.

For a rational number $r\ge 0$,
we consider the cartesian
diagram
$$\begin{CD}
P^{(r)}_{S,F}@>{i^{(r)}}>>
P^{(r)}_S@<{j^{(r)}}<<
P^{(r)}_{S,\eta}&=U\times \eta@>{{\rm pr}_1}>>U\\
@VVV @V{p^{(r)}}VV @VV{{\rm pr}_2}V @.\\
{\rm Spec}\ F
@>i>>
S@<j<< \eta&={\rm Spec}\ K. @.
\end{CD}$$
Let
$s^{(r)}:S\to P_S^{(r)}$
be the section
induced by $S\to P$.
For $r=0$,
we have $P^{(0)}_S=
P=(X\times_kS)^\sim$.
Let $\psi^{(r)}$ be the nearby cycle
functor for $p^{(r)}:P^{(r)}_S\to S$
and $\psi$ be the nearby cycle
functor for the identity $S\to S$.
For a sheaf
${\cal F}_\eta$ on $\eta$,
we identify $\psi({\cal F}_\eta)$
with the $G_K$-module
${\cal F}_{\bar \eta}$.

\begin{df}\label{dfFr}
Let ${\cal F}$ be
a locally constant 
constructible sheaf
of $\Lambda$-modules
on $U=X\setminus D$.
The stalk ${\cal F}_{\bar \eta}$
defines a representation
of the absolute Galois group $G_K$.

For a rational number $r>0$,
we say that the log ramification of ${\cal F}$
at $\xi$ along $D$ is bounded by $r$
if $G_{K,\log}^r$
acts trivially on
${\cal F}_{\bar \eta}$.
Similarly, for a rational number $r\ge 0$,
we say that the log ramification of ${\cal F}$
along $D$ is bounded by $r+$
if $G_{K,\log}^{r+}$
acts trivially on
${\cal F}_{\bar \eta}$.
\end{df}

Since $P=G_{K,\log}^{0+}$
is a pro-$p$ group,
there exists a unique
direct sum decomposition
\begin{equation}
{\cal F}_{\bar \eta}=
\bigoplus_{q\ge0,
q\in {\mathbb Q}}
{\cal F}_{\bar \eta}^{(q)}
\label{eqFq}
\end{equation}
by sub $G_K$-modules
such that
the $G_{K,\log}^{r+}$-fixed part
is given by
$${\cal F}_{\bar \eta}^{G_{K,\log}^{r+}}=
\bigoplus_{q\ge r}
{\cal F}_{\bar \eta}^{(q)}.$$
Replacing $X$ by
an \'etale neighborhood of $\xi$
if necessary,
we may assume that
there exists a
direct sum decomposition
${\cal F}=\bigoplus_{q\ge 0}
{\cal F}^{(q)}$
inducing (\ref{eqFq}).

We identify
the stable closed fiber
$\overline P_{S,{\overline F}}^{(r)}$
with the $\overline F$-vector space
$\Theta^{(r)}_{\log}$
by the isomorphism in
Corollary \ref{cormn}.

\begin{pr}\label{prgr}
Let $r> 0$ be a rational number and
let $\pi^{(r)}:
\overline P_{S,{\overline F}}^{(r)}
\to P_{S,{\overline F}}^{(r)}$
be the canonical map.
Let ${\cal F}$
be a smooth sheaf on $U$.
We assume
${\cal F}_{\bar \eta}=
{\cal F}_{\bar \eta}^{(q)}$
for a rational number
$q\ge 0$.

1. Assume $q=r$.
Let ${\cal F}_{\bar \eta}=
\bigoplus_\chi
{\cal F}_{\bar \eta}^{(\chi)}$
be the decomposition by
characters
$\chi:{\rm Gr}_{\log}^rG_K
\to \Lambda^\times$.
Let ${\cal L}_\chi$
be the smooth sheaf of rank $1$
on $\overline P_{S,{\overline F}}^{(r)}
=\Theta^{(r)}_{\log}$
defined by the composition
$\pi_1(\Theta^{(r)}_{\log})^{\rm ab}
\to {\rm Gr}_{\log}^rG_K
\to \Lambda^\times$.

Then, 
there exists a canonical isomorphism
\begin{equation}
\begin{CD}
\psi^{(r)}({\rm pr}_1^*{\cal F})
@>>>
\bigoplus_\chi
\pi^{(r)}_*{\cal L}_\chi
\otimes p^{(r)*}
{\cal F}_{\bar \eta}^{(\chi)}
\end{CD}
\label{eqps0}
\end{equation}
on $P^{(r)}_{S,\overline F}$.

2. 
If $q<r$, then
there exists a canonical isomorphism
\begin{equation}
\begin{CD}
\psi^{(r)}({\rm pr}_1^*{\cal F})
@>>>
\pi^{(r)}_*\Lambda 
\otimes p^{(r)*}\psi
({\cal F}_{\eta})
\end{CD}
\label{eqps<}
\end{equation}
on $P^{(r)}_{S,\overline F}$.
\end{pr}
{\it Proof.}
Let $V\to U$
be the finite \'etale covering
trivializing ${\cal F}$.
Replacing $X$ by
an \'etale neighborhood of $\xi$,
we may assume
that $V$ is the complement
of a smooth divisor
of the normalization $Y$
of $X$ in $V$
as in the previous subsection.
We consider the diagram (\ref{eqlTS}).
Since $q\le r$,
there exists a finite 
extension $K'$ of $K$ of ramification index $e$
such that $er$ is an integer
and that, 
for the base change by 
$S'={\rm Spec}\ {\cal O}_{K'}\to S$,
the map $Q^{(er)}_{T'}
\to P^{(er)}_{S'}$
is finite \'etale by Lemma \ref{lmbddr}.
If $q<r$,
the finite \'etale covering
$Q^{(er)}_{T',\overline F}
\to P^{(er)}_{S',\overline F}$
is trivial.

The pull-back
of ${\rm pr}_1^*{\cal F}$
to $U\times {\rm Spec}\ K'$
is extended to a smooth sheaf
${\cal G}$ on
$P^{(er)}_{S'}$.
By the definition of
the surjection
$\pi_1(\Theta^{(r)}_{\log})^{\rm ab}
\to {\rm Gr}_{\log}^rG_K$,
the pull-back of
${\cal G}$
to $\overline P_{S,{\overline F}}^{(r)}
=P^{(er)}_{S'}\times_{S'}{\overline F}$
is isomorphic to
$\bigoplus_{\chi}
{\cal L}_{\chi}^{{\rm rank}\ 
{\cal F}^{(\chi)}}$ if $q=r$.
If $q>r$,
the pull-back of
${\cal G}$
to 
$\overline P_{S,{\overline F}}^{(r)}$
is constant.

We consider the nearby cycle
functor $\psi'$ for the smooth map
$p^{\prime (er)}:P^{(er)}_{S'}\to S'$.
Let $s':S'\to P^{(er)}_{S'}$
be the section induced by
$S\to P^{(r)}_S$. Then
$\psi'({\rm pr}_1^*{\cal F})$
is the restriction of
${\cal G}$
on
$P^{(er)}_{S',\overline F}$
and the base change map 
\begin{equation}
\begin{CD}
s^{\prime*}\psi'({\rm pr}_1^*{\cal F})
={\cal G}_0
@>>>
\psi({\cal F})
\end{CD}
\label{eqbc0'}
\end{equation}
is an isomorphism,
where $0\in \Theta_{\log}^{(r)}$ 
denotes the origin.
Thus we obtain a canonical isomorphism
\begin{equation}
\begin{CD}
\psi'({\rm pr}_1^*{\cal F})
@>>>
\begin{cases}
\bigoplus_{\chi}
{\cal L}_{\chi}
\otimes p^{\prime (er)*}{\cal F}_{\bar\eta}^{(\chi)}
&\quad \text{ if }q=r\\
p^{\prime (er)*}\psi({\cal F}_{\eta})
&\quad \text{ if }q>r.
\end{cases}
\end{CD}
\label{egLchi4}
\end{equation}
Since $
\psi^{(r)}({\rm pr}_1^*{\cal F})=
\pi^{(r)}_*
\psi'({\rm pr}_1^*{\cal F})$,
the isomorphism
(\ref{egLchi4})
induces isomorphisms
(\ref{eqps0})
and
(\ref{eqps<}).
\qed

\begin{cor}\label{corgr00}
Let $r\ge 0$ be a rational number.

1.
If the log ramification
of ${\cal F}$ is bounded by $r+$,
then the base change map
\begin{equation}
\begin{CD}
s^{(r)*}\psi^{(r)}({\rm pr}_1^*{\cal F})
@>>>
\psi({\cal F}_\eta)
\end{CD}
\label{eqbc0}
\end{equation}
is an isomorphism.

2.
The base change map
{\rm (\ref{eqbc0})}
induces an isomorphism
\begin{equation}
\begin{CD}
s^{(r)*}\psi^{(r)0}({\rm pr}_1^*{\cal F})
@>>>
{\cal F}_{\bar \eta}^{G_{K,\log}^{r+}}
\subset 
{\cal F}_{\bar \eta}=
\psi({\cal F}_\eta)
\end{CD}
\label{eqbc00}
\end{equation}
to the 
$G_{K,\log}^{r+}$-fixed part.
\end{cor}
{\it Proof.}
1. We may assume 
${\cal F}={\cal F}^{(q)}$
for some rational number $0\le q\le r$.
First, we consider the case $r>0$.
We use the notation of
the proof of Proposition \ref{prgr}.
Since the inverse image $\pi^{-1}(0)$
of $0\in \Theta_{\log}^{(r)}
=\overline P^{(r)}_{S,\overline F}$
consists of
the image of the geometric closed point
by the section
$S\to P^{(r)}_S$,
the isomorphism
(\ref{eqbc0'}) shows that
the base change map
(\ref{eqbc0}) is an isomorphism.

Assume $r=q=0$.
Then, the smooth sheaf
${\rm pr}_1^*{\cal F}$
on $U\times \eta\subset P=P^{(0)}_S$
is tamely ramified
along $P\times_S{\rm Spec}\ F$.
Hence, \'etale locally on $P$, 
it is isomorphic
to the pull-back of
a sheaf on $\eta$.
Since $P$ is smooth over $S$,
the assertion follows.

2. We may assume 
${\cal F}={\cal F}^{(q)}$
for some rational number $q\ge 0$.
By 1,
it suffices to consider the case
$q>r$.
Since the base change map
$\psi^{(r)0}({\cal F})_{\bar s}
\to \psi^0({\cal F}_\eta)$
is injective,
it suffices to
show that 
the base change map is the 0-map.

Let $f_{rq}:
P_S^{(q)}\to P_S^{(r)}$
be the canonical map.
By Proposition \ref{prgr},
the sheaf $\psi^{(q)}({\cal F})$
has no non-trivial geometrically
constant subsheaf.
Since the image
$f_{rq}(P_{S,\overline F}^{(q)})$
is a point,
the base change map 
$f_{rq}^*
\psi^{(r)}({\cal F})
\to \psi^{(q)}({\cal F})$
is the 0-map.
Thus the composition
$\psi^{(r)}({\cal F})_{\bar s}
\to \psi^{(q)}({\cal F})_{\bar s}
\to \psi({\cal F}_\eta)$
is also the 0-map
as required.
\qed

We consider ${\cal H}=
{\cal H}om({\rm pr}_2^*{\cal F},
{\rm pr}_1^*{\cal F})$
on $P^{(r)}_{S,\eta}=
U\times \eta$ and
the base change map
with respect to the diagram
$$\begin{CD}
U\times \eta
@>{j^{(r)}}>> P^{(r)}_S\\
@AAA @AA{s^{(r)}}A\\
\eta @>j>> S.
\end{CD}$$

\begin{cor}\label{corgr0}
Let $r\ge 0$ be a rational number.

1. The following conditions
are equivalent:

{\rm (1)} The log
ramification of ${\cal F}$
is bounded by $r+$.

{\rm (2)} 
The base change map
$$\begin{CD}
s^{(r)*}j^{(r)}_*{\cal H}
@>>> j_*{\cal E}nd({\cal F}_\eta)
\end{CD}$$
is an isomorphism.

{\rm (3)}
The identity
$1\in {\rm End}_{G_K}({\cal F}_{\bar \eta})
=\Gamma(S,
j_*{\cal E}nd({\cal F}_\eta))$
is in the image of
the base change map
$$\begin{CD}
\Gamma(S,s^{(r)*}j^{(r)}_*{\cal H})
@>>> \Gamma(S,j_*{\cal E}nd({\cal F}_\eta)).
\end{CD}$$

2.
Assume that the 
$G_{K,\log}^{r+}$-fixed part
${\cal F}_{\bar \eta}^{G_{K,\log}^{r+}}$
is $0$.
Then we have
$i^*s^{(r)*}j^{(r)}_*{\cal H}=0$.
\end{cor}
{\it Proof.}
1. (1)$\Rightarrow$(2)
It suffices to show the isomorphism
for the geometric closed fiber
at $\bar s={\rm Spec}\ \overline F\to S$.
By the assumption (1),
we have 
${\cal F}_{\bar \eta}=
{\cal F}_{\bar \eta}^{G_{K,\log}^{r+}}$
and the base change map (\ref{eqbc00}) 
induces an isomorphism
$$\begin{CD}
s^{(r)*}\psi^{(r)0}({\cal H})_{\bar s}
={\rm Hom}(\psi({\cal F}_{\eta}),
s^{(r)*}\psi^{(r)0}({\rm pr}_1^*{\cal F}))
\to
{\rm Hom}(\psi({\cal F}_{\eta}),
\psi({\cal F}_{\eta}))=
{\rm End}(\psi({\cal F}_{\eta})).
\end{CD}$$
Taking the fixed parts by
the inertia subgroup $I\subset G_K$,
we obtain an isomorphism
$(s^{(r)*}j^{(r)}_*{\cal H})_{\bar s}=
s^{(r)*}\psi^{(r)0}({\cal H})^I
\to {\rm End}_I(\psi({\cal F}_{\eta}))=
(j_*{\cal E}nd({\cal F}_\eta))_{\bar s}$
as required.

(2)$\Rightarrow$(3)
Clear.

(3)$\Rightarrow$(1)
We consider the direct sum decomposition
${\cal F}_{\bar \eta}=
\bigoplus_q{\cal F}_{\bar \eta}^{(q)}$.
It suffices to show
that the identity is not
in the image
assuming ${\cal F}_{\bar \eta}^{(q)}\neq 0$
for some $q>r$.
Thus it is reduced to the assertion 2.

2.
Assume ${\cal F}_{\bar \eta}^{G_{K,\log}^{r+}}=0$.
Then, similarly as in the proof of
1 (1)$\Rightarrow$(2)
above, we have
$(s^{(r)*}j^{(r)}_*{\cal H})_{\bar s}=
s^{(r)*}\psi^{(r)0}({\cal H})^I=0$.
\qed

\begin{cor}\label{corgr1}
Assume that $q=r>0$ is an integer
and that the restriction to 
$G^r_{K,\log}$ of
the action on
${\cal F}_{\bar \eta}$
is by the multiplication
by a character
$\chi:{\rm Gr}^r_{\log}G_K
\to \Lambda^\times$.

1. 
There exists a canonical isomorphism
\begin{equation}
\begin{CD}
\psi^{(r)}({\rm pr}_1^*{\cal F})
@>>>
{\cal L}_\chi
\otimes p^{(r)*}\psi
({\cal F}_{\eta})
\end{CD}
\label{eqpsr}
\end{equation}
on $P^{(r)}_{S,\overline F}$.

2. 
There exists a canonical isomorphism
\begin{equation}
\begin{CD}
i^{(r)*}j^{(r)}_*{\cal H}
@>>>
{\cal L}_\chi
\otimes p^{(r)*}i^*j_*
{\cal E}nd({\cal F}_{\eta})
\end{CD}
\label{eqH}
\end{equation}
on $P^{(r)}_{S,F}$.
\end{cor}
{\it Proof.}
1. Clear from Proposition \ref{prgr}.1.

2. By 1,
we have an isomorphism
\begin{equation}
\begin{CD}
\!\!\!\!\!\!\!\!\!
\!\!\!\!\!\!\!\!\!
\!\!\!\!\!\!\!\!\!
\!\!\!\!\!\!\!\!\!
\!\!\!\!\!\!\!\!\!
\psi^{(r)}({\cal H})=
@.
\\
{\rm Hom}(p^{(r)*}\psi({\cal F}_{\eta}),
\psi^{(r)}({\rm pr}_1^*{\cal F}))
@>>>
{\rm Hom}(p^{(r)*}\psi({\cal F}_\eta),
{\cal L}_\chi
\otimes p^{(r)*}\psi({\cal F}_{\eta}))\\
@.
=
{\cal L}_\chi
\otimes p^{(r)*}\psi({\cal E}nd({\cal F}_{\eta})).
\end{CD}
\label{eqpsi}
\end{equation}
We have canonical isomorphisms
$R\Gamma(I,\psi^{(r)})\to i^{(r)*}Rj^{(r)}_*$
and
$R\Gamma(I,\psi)\to i^*Rj_*$
of functors.
Thus, we obtain
the isomorphism (\ref{eqH})
by taking the inertia fixed
parts in (\ref{eqpsi}).
\qed

The following geometric construction
is crucial in the proof
of Theorem \ref{thmgr}.
Let $(X\times X)'\to X\times X$
be the blow-up at $D\times D$
and let
$(X\times X)^\sim \subset
(X\times X)'$
be the complement of
the proper transforms
of $D\times X$ and of $X\times D$.
We call the immersion 
$\tilde \delta:X\to (X\times X)^\sim$
induced by
the diagonal
$\delta:X\to X\times X$
the log diagonal.
Let ${\cal J}_X\subset (X\times X)^\sim$
be the ideal defining
the log diagonal
and let
$\tilde j:U\times U\to (X\times X)^\sim$
be the open immersion.
For an integer $r\ge 0$,
we define
a scheme $(X\times X)^{(r)}$
affine over $(X\times X)^\sim$
by the quasi-coherent 
${\cal O}_{(X\times X)^\sim}$-algebra
$\sum_{l\ge 0}
{\cal I}_D^{-lr}\cdot {\cal J}_X^l
\subset \tilde j_*{\cal O}_{U\times U}$.

Similarly as Corollary \ref{cormn},
the fiber product
$(X\times X)^{(r)}
\times_XD$
is canonically identified
with the vector bundle
${\mathbf V}
(\Omega^1_X(\log D)(rD))
\times_XD$.
Hence the map
$(X\times S)^{(r)}
\to 
(X\times X)^{(r)}$
defined
by the canonical map
$S\to X$
induces an isomorphism
\begin{equation}
\begin{CD}
(X\times S)^{(r)}
\times_X{\rm Spec}\ \overline F
@>>> 
\!\!\!\!\!\!\!\!
\!\!\!\!\!\!\!\!
\!\!\!\!\!\!\!\!
\!\!\!\!\!\!\!\!
\!\!\!\!\!\!\!\!
(X\times X)^{(r)}
\times_X{\rm Spec}\ \overline F\\
\!\!\!\!\!\!\!\!
\!\!\!\!\!\!\!\!
\!\!\!\!\!\!\!\!
=\Theta_{\log}^{(r)}
@.
\qquad={\mathbf V}
(\Omega^1_X(\log D)(rD))
\times_X{\rm Spec}\ \overline F.
\end{CD}
\label{eqThD}
\end{equation}

\begin{lm}\label{lmmu}
Let $r>0$ be an integer.
Then there exists a unique map
$\mu:
(X\times S)^{(r)}\times_S
(X\times S)^{(r)}\to
(X\times X)^{(r)}$
that makes the diagram
\begin{equation}
\begin{CD}
(X\times S)^{(r)}\times_S
(X\times S)^{(r)}@>{\mu}>>
(X\times X)^{(r)}\\
@VVV @VVV\\
(X\times S)\times_S
(X\times S)=
X\times X\times S
@>{{\rm pr}_{12}}>>
X\times X
\end{CD}
\label{eqmu}
\end{equation}
commutative.

2.
Under the identification
{\rm (\ref{eqThD})}
$(X\times X)^{(r)}
\times_X{\rm Spec}\ \overline F=
\Theta_{\log}^{(r)}$,
the map
$\mu:
(X\times S)^{(r)}\times_S
(X\times S)^{(r)}
\to
(X\times X)^{(r)}$
induces the difference
$-:\Theta_{\log}^{(r)}
\times_{\overline F} \Theta_{\log}^{(r)}
\to
\Theta_{\log}^{(r)}$
on the fiber over ${\rm Spec}\ \overline F$.
\end{lm}
{\it Proof.}
We put
$P=(X\times S)^\sim\times_S
(X\times S)^\sim$.
Applying the basic construction
to the smooth scheme
$P$ and the diagonal section
$S\to P$,
we define
$q:P^{(r)}_S\to P$
and a section
$s^{(r)}:S\to P^{(r)}_S$.
The projections
$P\to (X\times S)^\sim$
induce
$P^{(r)}_S\to (X\times S)^{(r)}$.
We show that the product
\begin{equation}
P^{(r)}_S\to 
(X\times S)^{(r)}
\times_S(X\times S)^{(r)}
\label{eqisP}
\end{equation}
is an isomorphism.
Since the ideal
defining the closed
subscheme $S\subset P$
is generated by
the two pull-backs of
the ideal
defining the closed
subscheme $S\subset (X\times S)^\sim$.
Hence,
the map
(\ref{eqisP})
is a closed immersion.
Since both
$P^{(r)}_S$ and 
$(X\times S)^{(r)}
\times_S(X\times S)^{(r)}$
are smooth over $S$ 
of the same dimension,
the closed immersion
(\ref{eqisP})
is an open immersion.
Since the map
(\ref{eqisP}) is an isomorphism
on each fiber,
it is an isomorphism.

Let $D_{(X\times S)^\sim}
\subset (X\times S)^\sim$
be the pull-back 
${\rm pr}_1^*D=
{\rm pr}_2^*D_S$.
Since ${\rm pr}_1^*D_{(X\times S)^\sim}=
{\rm pr}_2^*D_{(X\times S)^\sim}$
on $P$,
there exists a unique map
$\lambda:P\to
(X\times X)^\sim$
that makes the diagram
(\ref{eqmu}) with
$\mu:P^{(r)}=
(X\times S)^{(r)}
\times_S(X\times S)^{(r)}
\to
(X\times X)^{(r)}$
replaced by
$\lambda:P\to
(X\times X)^\sim$
commutative.
By the commutative diagram
$$\begin{CD}
S@>>> X\\
@VVV @VVV\\
P@>{\lambda}>>(X\times X)^\sim,
\end{CD}$$
the pull-back
$\lambda^*({\cal I}_D^{-lr}\cdot {\cal I}_X^l)$
is contained in
${\mathfrak m}_K^{-lr}\cdot {\cal I}_S^l$.
Hence the assertion follows.

2. Let ${\cal J}_X\subset 
{\cal O}_{(X\times X)^\sim}$
and ${\cal J}_S\subset 
{\cal O}_{(X\times S)^\sim}$
be the ideals defining
the closed subschemes
$X\subset (X\times X)^\sim$
and 
$S\subset (X\times S)^\sim$
respectively.
By the identification
in Corollary \ref{cormn},
the map $\Theta_{\log}^{(r)}
\times_F \Theta_{\log}^{(r)}
\to
\Theta_{\log}^{(r)}
\subset
(X\times X)^{(r)}$
is defined by
${\cal I}_D^{-r}{\cal J}_X
\to {\mathfrak m}_K^{-r}\cdot {\cal J}_S
\oplus {\mathfrak m}_K^{-r}\cdot {\cal J}_S$.
Hence,
it is a linear map of vector bundles.
Thus it suffices to
show that 
the composition
with the injections
$i_1,i_2:\Theta_{\log}^{(r)}
\to
\Theta_{\log}^{(r)}
\times_F \Theta_{\log}^{(r)}$
of the two factors
are the identity of
$\Theta_{\log}^{(r)}$
and the multiplication 
by $-1$ respectively.

Let $s:S\to (X\times S)^{(r)}$
be the map induced by
the canonical map $S\to X$.
We consider the map
$\iota_1=
({\rm id}_{(X\times S)^{(r)}},
s\circ {\rm pr}_2)
:(X\times S)^{(r)}
\to
(X\times S)^{(r)}
\times_S
(X\times S)^{(r)}$.
Then, its restriction
$\Theta_{\log}^{(r)}
\to
\Theta_{\log}^{(r)}
\times_{\overline F} 
\Theta_{\log}^{(r)}$
to the closed fiber
is the injection
to the first component.
The composition
$\mu\circ \iota_1$
is the map
$(X\times S)^{(r)}
\to
(X\times X)^{(r)}$
induced by the canonical map
$S\to X$.
Hence the composition
$\mu\circ i_1$ 
is the identity
of $\Theta_{\log}^{(r)}$.
Similarly,
we consider the map
$\iota_2=
(s\circ {\rm pr}_2,
{\rm id}_{(X\times S)^{(r)}})
:(X\times S)^{(r)}
\to
(X\times S)^{(r)}
\times_S
(X\times S)^{(r)}$.
Then the composition
$\mu\circ \iota_2
:(X\times S)^{(r)}
\to
(X\times X)^{(r)}$
is the composition of
the canonical map
$(X\times S)^{(r)}
\to
(X\times X)^{(r)}$
and the map
$(X\times X)^{(r)}
\to
(X\times X)^{(r)}$
switching the two factors.
Hence the composition
$\mu\circ i_2$ 
is the multiplication
by $-1$
of $\Theta_{\log}^{(r)}$.
Hence the assertion 
is proved.
\qed

{\it Proof of Theorem \ref{thmgr}.}
We start with some reduction steps.
For each non-trivial
character $\chi:
{\rm Gr}_{\log}^rG_K
\to \Lambda^\times$,
the surjection (\ref{eqGr})
defines a 
locally constant
sheaf ${\cal L}_\chi$
of $\Lambda$-modules
of rank $1$
on $\Theta^{(r)}_{\log}$.
By Lemma \ref{lmisog},
in order to prove Theorem \ref{thmgr},
it suffices to show that,
for every character
$\chi:
{\rm Gr}_{\log}^rG_K
\to \Lambda^\times$,
there exists an isomorphism
$-^*{\cal L}_\chi
\to {\cal H}om(p_2^*
{\cal L}_\chi,p_1^*
{\cal L}_\chi)$
assuming $\Lambda$ is a finite field.

We reduce it to the
case where $r$ is an integer.
Let $e>0$ be an integer
such that $er$
is an integer and
let $K_1$ be an extension
of $K$ of ramification index $e$
as in Lemma \ref{lmr}.
Then the construction of ${\cal L}_\chi$
commutes with the base change
$K\to K_1$.
Hence, it is reduced
to the case where
$r$ is an integer.

We show that 
it is reduced to the case
where the restriction
of the action on ${\cal F}_{\bar \eta}$
to $G_{K,\log}^r$
is by the multiplication
by a character
$\chi:
{\rm Gr}_{\log}^rG_K
\to \Lambda^\times$.
By the same argument as
in the last paragraph,
we may replace $K$
by a tamely ramified
extension.
Hence we may assume
the restriction
to $G_{K,\log}^{0+}$
is irreducible.
By \cite{AS2}
Theorem 5.12.1,
${\rm Gr}_{\log}^rG_K$
is in the center of
$G_{K,\log}^{0+}/
G_{K,\log}^{r+}$.
Hence, 
the action on ${\cal F}_{\bar \eta}$
to $G_{K,\log}^r$
is by the multiplication
by a character
$\chi:
{\rm Gr}_{\log}^rG_K
\to \Lambda^\times$.

We assume that $r>0$
is an integer
and the restriction
of ${\cal F}_{\bar \eta}$
to $G_{K,\log}^r$
is the multiplication 
by a non-trivial character
$\chi:{\rm Gr}^r_{\log}G_K\to 
\Lambda^\times$.
We consider the commutative diagram
$$\begin{CD}
\Theta_{\log}^{(r)}
\times_F \Theta_{\log}^{(r)}
@>i>>
(X\times S)^{(r)}\times_S
(X\times S)^{(r)}
@<j<<
(U\times \eta)\times_\eta
(U\times \eta)\\
@V-VV @.
=U\times U\times \eta\\
\Theta_{\log}^{(r)}&& @V{\mu}VV
@VV{{\rm pr}_{12}}V\\
@VVV\\
(X\times X)^{(r)}\times_XD
@>{i'}>>
(X\times X)^{(r)}@<{j'}<<
U\times U.
\end{CD}$$
The left square is commutative by
Lemma \ref{lmmu}.2.
We consider the base change map
\begin{equation}
\begin{CD}
-^*((i^{\prime *}j'_* {\cal H})
|_{\Theta_{\log}^{(r)}})
@>>> i^*j_*{\rm pr}_{12}^*{\cal H}
\end{CD}\label{eqbcmu}
\end{equation}
for ${\cal H}=
{\cal H}om({\rm pr}_2^*{\cal F},
{\rm pr}_1^*{\cal F})$
on $U\times U$.

First, we compute
$i^*j_*{\rm pr}_{12}^*{\cal H}$.
We have $
\psi({\rm pr}_{12}^*{\cal H})=
\psi(
{\cal H}om({\rm pr}_2^*{\cal F},
{\rm pr}_1^*{\cal F}))$
where ${\rm pr}_i:
U\times U\times \eta
\to U$ denote the projections.
Further, we have
$\psi(
{\cal H}om({\rm pr}_2^*{\cal F},
{\rm pr}_1^*{\cal F}))
=
{\cal H}om
(\psi({\rm pr}_2^*{\cal F}),
\psi({\rm pr}_1^*{\cal F}))
=
{\cal H}om
({\rm pr}_2^*\psi^{(r)}{\cal F},
{\rm pr}_1^*\psi^{(r)}{\cal F})$
where ${\rm pr}_i:
\Theta_{\log}^{(r)}\times 
\Theta_{\log}^{(r)}
\to \Theta_{\log}^{(r)}$ 
denote the projections
in the right hand side.
By Proposition \ref{prgr},
it is further identified with
${\cal H}om(
{\rm pr}_2^*{\cal L}_\chi,
{\rm pr}_1^*{\cal L}_{\chi})
\otimes 
\psi{\cal E}nd(
{\cal F}_{\eta})$.
Thus, similarly as Corollary \ref{corgr1},
we obtain an isomorphism
$i^*j_*{\rm pr}_{12}^*{\cal H}
\to {\cal H}om(
p_2^*{\cal L}_\chi,
p_1^*{\cal L}_{\chi})
\otimes 
{\rm End}_I(
{\cal F}_{\eta})$
by taking the inertia fixed parts.

Next, we compute the restriction
$(i^{\prime *}j'_* {\cal H})
|_{\Theta_{\log}^{(r)}}$.
This is the same as $i^*j^{(r)}_*{\cal H}$
computed
in Corollary \ref{corgr1}.
Hence it is canonically
isomorphic to
${\cal L}_{\chi}
\otimes 
{\rm End}_I(
{\cal F}_{\eta})$.
Hence, the map (\ref{eqbcmu})
induces a map
$$-^*{\cal L}_{\chi}
\otimes 
{\rm End}_I(
{\cal F}_{\eta})
\to 
{\cal H}om(
p_2^*{\cal L}_\chi,
p_1^*{\cal L}_{\chi})
\otimes 
{\rm End}_I(
{\cal F}_{\eta})$$
of smooth sheaves.
Since, this is an isomorphism at
the origin,
it is an isomorphism on 
$\Theta_{\log}^{(r)}
\times_F \Theta_{\log}^{(r)}$.
By evaluating
at the identity of ${\cal F}_\eta$,
we obtain an isomorphism
$-^*{\cal L}_{\chi}
\to 
{\cal H}om(
p_2^*{\cal L}_\chi,
p_1^*{\cal L}_{\chi})$
as required.
\qed

\section{Ramification along a divisor}

We introduce the notion
of additive sheaves
on vector bundles
and its generalization.
in \S2.1.
In \S2.2, we study a global variant
of the basic construction in \S1.1.
After these preliminaries,
we study the ramification of
smooth sheaves on the complement
of a divisor with normal crossings
along the divisor
in \S2.3.

\subsection{Additive sheaves
on vector bundles
and generalizations}

We recall the Fourier-Deligne transform
\cite{KL}.
Let $X$ be a scheme over ${\mathbb F}_p$.
Let $E={\mathbf V}({\cal E})\to X$ be 
a vector bundle of rank $d$
and let
$E^\vee={\mathbf V}({\cal E}^\vee)\to X$ be 
the dual.
The canonical pairing defines a map
$\langle\ ,\ \rangle:
E \times_XE^\vee
\to {\mathbf A}^1$.
We consider the diagram
$$\begin{CD}
E@<{{\rm pr}_1}<<
E \times_XE^\vee
@>{\langle\ ,\ \rangle}>> {\mathbf A}^1\\
@. @V{{\rm pr}_2}VV @.\\
@. E^\vee @.
\end{CD}$$
where ${\rm pr}_i$
denote the projections.

We fix a non-trivial character
$\psi:{\mathbb F}_p\to \Lambda^\times$
and let ${\cal L}_\psi$
be the smooth rank 1 
Artin-Schreier sheaf on
${\mathbf A}^1={\rm Spec}\ k[t]$
defined by the ${\mathbb F}_p$-torsor
${\mathbf A}^1\to
{\mathbf A}^1: t\mapsto
t^p-t$.
For a sheaf ${\cal G}$ on 
the dual $E^\vee$ of
a vector bundle $E$,
we define the naive Fourier transform
$F_\psi({\cal G})$ on $E$
by
$$F_\psi({\cal G})
=R{\rm pr}_{1!}
({\rm pr}_2^*{\cal G}
\otimes 
\langle\ ,\ \rangle^*
{\cal L}_\psi)$$
For a sheaf
${\cal H}$
on $E$, we
define the inverse Fourier transform
$F'_{\psi'}({\cal H})$
by 
$$F'_{\psi'}({\cal H})
=R{\rm pr}_{2!}
({\rm pr}_1^*{\cal H}
\otimes 
\langle\ ,\ \rangle^*
{\cal L}_{\psi'})(d)[2d]$$
where $\psi'
:{\mathbb F}_p\to \Lambda^\times$
denotes the inverse of $\psi$.

We have a canonical isomorphism
\begin{equation}
{\cal H}
\to F_\psi F'_{\psi'}{\cal H}.
\label{eqFi}
\end{equation}
Let $f:E\to F$
be a linear morphism
of vector bundles
over $X$
and $f^\vee:F^\vee\to E^\vee$
be the dual.
Then, we have
a canonical isomorphism
\begin{equation}
f^*F_\psi{\cal G}
\to F_\psi Rf^\vee_!{\cal G}
\label{eqFf}
\end{equation}
for a sheaf ${\cal G}$ on $F^\vee$.
Dually,
we have a canonical isomorphism
\begin{equation}
Rf_*F_\psi{\cal G}
\to F_\psi Rf^{\vee !}{\cal G}
\label{eqFf*}
\end{equation}
for a sheaf ${\cal G}$ on $E^\vee$.

We introduce the
notion of additive sheaves
on vector bundles.

\begin{df}\label{dfadd}
Let $E=
{\mathbf V}({\cal E})$
be a vector bundle
over a scheme $X$
over $k$ and
let ${\cal H}$ be 
a constructible sheaf on $E$.
Let ${\cal G}=F'_{\psi'} {\cal H}$ 
be the inverse Fourier transform
and define a constructible
subset $S\subset E^\vee$
to be the support of ${\cal G}$.

We say ${\cal H}$ on $E$
is additive if,
for every point
$x$ of $X$, the fiber 
$S\times_Xx$ is finite.
For an additive constructible sheaf
${\cal H}$ on $E$,
we call the support
$S=S_{\cal H}\subset E^\vee$
of the inverse  Fourier transform
${\cal G}=F'_{\psi'} {\cal H}$
the dual support of ${\cal H}$.
We say an additive constuctible sheaf
is non-degenerate
if the intersection of
the {\em closure} of the dual support
$S_{\cal H}$ with
the $0$-section is empty.
\end{df}

A constructible sheaf ${\cal H}$
on a vector bundle $E$
is additive if and only if,
for every geometric point $\bar x\to X$,
the pull-back ${\cal H}|_{E_{\bar x}}$
is additive.
If $X={\rm Spec}\ \overline F$
is the spectrum of an algebraically closed 
field,
a constructible sheaf 
${\cal H}$ on a vector space
$E$ is additive if and only if
${\cal H}$ is a direct sum
of rank 1 Artin-Schreier sheaves 
defined by linear forms
by the isomorphism (\ref{eqFi}).
A constructible subsheaf ${\cal H}'$
of an additive constructible sheaf 
${\cal H}$ is additive
if it is fiberwisely smooth.

We have the following elementary properties
on additive sheaves.

\begin{lm}\label{lmadd2}
1. Let $f:E'\to E$ be a linear
map of vector bundles over $X$
and $f^\vee:E^\vee\to E^{\prime \vee}$ be the dual.
If ${\cal H}$ is additive,
then $f^*{\cal H}$ is additive
and we have $S_{f^*{\cal H}}=
f^\vee (S_{\cal H})$.

Assume $f:E'\to E$
is surjective and identify
$E^\vee$ with the image $f^\vee(E^\vee)$
by the closed immersion
$f^\vee:E^\vee\to E^{\prime \vee}$.
Then, conversely, if $f^*{\cal H}$
is additive,
then ${\cal H}$ is additive.

2. Let $f:E\to X$
be a vector bundle and
let ${\cal H}$
be an additive 
constructible sheaf.
If ${\cal H}$ is non-degenerate, 
we have
$Rf_*{\cal H}=Rf_!{\cal H}=0$.
\end{lm}
{\it Proof.}
1. Clear from (\ref{eqFf}).

2. Clear from (\ref{eqFf})
and (\ref{eqFf*}).
\qed

\begin{pr}\label{pradd}
Let $X$ be a scheme 
over $k$
and $E\to X$ be a vector bundle.
For a constructible
sheaf ${\cal H}$ on $E$,
the following conditions
are equivalent:

{\rm (1)}
${\cal H}$ is additive.

{\rm (2)}
For every geometric point $\bar x\in X$
and for every closed point $a\in E_{\bar x}$,
there exists an isomorphism
$(+a)^*({\cal H}|_{E_{\bar x}})
\to
{\cal H}|_{E_{\bar x}}$.
\end{pr}
{\it Proof.}
We may assume $k$ is
algebraically closed
and $X={\rm Spec}\ k$.
Let ${\cal G}=F_\psi {\cal H}$
be the Fourier transform
and $S\subset E^\vee$
be the support of ${\cal G}$.
A closed point $a\in E$
defines a linear form
$\langle a,\ \rangle:
E^\vee\to {\mathbf A}^1$.
The conditions (1) and (2)
are equivalent to the
following conditions respectively:

(1$'$)
For every closed point $a\in E$,
the image of $S$ by
the map
$\langle a,\ \rangle:
E^\vee\to {\mathbf A}^1$
is finite.

(2$'$)
For every closed point $a\in E$,
there exists an isomorphism
${\cal G}\otimes 
\langle a,\ \rangle^*{\cal L}_\psi
\to
{\cal G}$.

It is clear that
the condition (1$'$) implies (2$'$).
We show (2$'$) implies (1$'$).
Let $U\subset E^\vee$
be a normal integral 
locally closed subscheme
supported in $S$
such that the restriction
${\cal G}|_U$ is locally constant.
Let $\pi:V\to U$
be a connected finite
etale covering
such that
$\pi^*{\cal G}|_U$ is constant.
Then, by the condition (2$'$),
$\pi^*\langle ca,\ \rangle^*{\cal L}_\psi$
is constant on $V$
for every $a\in E$
and $c\in k$.
Hence,
the map
$(\langle a,\ \rangle\circ \pi)_*:
\pi_1(V)^{\rm ab}\to
\pi_1({\mathbf A}^1)^{\rm ab}$
has infinite cokernel.
Therefore
the image of
$\langle a,\ \rangle\circ \pi:
V\to {\mathbf A}^1$
collapses to a point.
Hence the condition
(2$'$) implies (1$'$).
\qed

\begin{pr}\label{pradd2}
Let $X$ be a scheme 
over $k$
and $E\to X$ be a vector bundle.
Let ${\cal H}$ be
an additive constructible
sheaf on $E$
and ${\cal K}$ be
a constructible
sheaf on $E$.
Let $e\in \Gamma(X,{\cal H}|_0)$
be a section of
the restriction ${\cal H}|_0$
on the $0$-section $X\subset E$
and
$u:{\cal H}\boxtimes
{\cal K}\to +^*{\cal K}$
be a map such that
the composition
\begin{equation}
\begin{CD}
u|_{0\times E}
\circ (e\otimes 1_{\cal K})
:{\cal K}
@>>>
{\cal H}|_0
\otimes {\cal K}
@>>>
{\cal K}
\end{CD}
\label{eque}
\end{equation}
is the identity of
${\cal K}$.
Then ${\cal K}$
is additive and
the support $S_{\cal M}\subset E^\vee$
of ${\cal M}=F_\psi {\cal K}$
is a subset of
the support $S_{\cal G}\subset E^\vee$
of ${\cal G}=F_\psi {\cal H}$.
\end{pr}
{\it Proof.}
We regard $e$
as a global section
$e\in \Gamma(E^\vee,{\cal G})
=\Gamma(X,{\cal H}|_0)$.
The map $u$
induces 
${\cal G}\boxtimes
{\cal M}
\to \delta_*{\cal M}$
on the Fourier transform
and hence a bilinear map
$v:{\cal G}\otimes
{\cal M}\to {\cal M}$
by adjunction.
We show that the composition
\begin{equation}
\begin{CD}
v\circ (e\otimes 1_{\cal M})
:{\cal M}
@>>>
{\cal G}
\otimes {\cal M}
@>>>
{\cal M}
\end{CD}
\label{eqve}
\end{equation}
is the identity of
${\cal M}$.
Let $\tilde e:
\Lambda_X(-d)[-2d]
\to {\cal H}$
be the cup-product
of $e:\Lambda
\to {\cal H}|_0$
with the map
$\Lambda_X(-d)[-2d]
\to \Lambda_E$
defined by the cycle
class of the 0-section $X
\subset E$.
We consider the map
\begin{equation}
\begin{CD}
u\circ (\tilde e\boxtimes 1_{\cal K})
:\Lambda_X(-d)[-2d]
\boxtimes
{\cal K}
@>>>
{\cal H}
\boxtimes {\cal K}
@>>>
+^*{\cal K}.
\end{CD}
\label{equte}
\end{equation}
By the assumption
that the composition
of (\ref{eque})
is the identity,
the induced map
$+_*(\Lambda_X(-d)[-2d]
\boxtimes
{\cal K})
={\cal K}(-d)[-2d]
\to 
+_*+^*{\cal K}
={\cal K}(-d)[-2d]$
is the identity map.
Therefore,
the Fourier transform
$$\begin{CD}
F_\psi u\circ (e\boxtimes 1_{\cal M})
:\Lambda
\boxtimes
{\cal M}
@>>>
{\cal G}
\boxtimes {\cal M}
@>>>
\delta_*{\cal M}.
\end{CD}$$
of (\ref{equte})
induces the identity
in (\ref{eqve}).

Since the composition 
in (\ref{eqve})
is the identity of ${\cal M}$,
the support $S_{\cal M}$
is a subset of
the support of $e
\in \Gamma(E^\vee,{\cal G})$.
Hence we have
$S_{\cal M}
\subset S_{\cal G}$
and ${\cal K}$ is additive.
\qed

\begin{lm}\label{lmadd0}
Let $X$ be a normal scheme 
over $k$
and $E\to X$ be a vector bundle.
Let ${\cal H}$
be an additive constructible
sheaf on $E$
satisfying the following
condition:
There exists a dense open subscheme
$U\subset X$ 
such that, if
$j:E_U=E\times_XU
\to E$
denote the open immersion,
the pull-back
${\cal H}_U=j^*{\cal H}$
is locally constant
and that the canonical
map ${\cal H}\to j_*j^*{\cal H}$
is injective.
Then, we have $S_{{\cal H}}
\subset 
\overline{S_{{\cal H}_U}}
\subset E^\vee$.
\end{lm}
{\it Proof.}
It suffices to show
that,
for each  $x\in X\setminus U$,
we have $S_{{\cal H}|_{E_x}}
\subset 
\overline{S_{{\cal H}_U}}$.
By replacing $X$
by the normalization
of the blowing-up of $X$
at the closure $\overline{\{x\}}$,
we may assume ${\cal O}_{X,x}$
is a discrete valuation ring.
Further replacing $X$
by the normalization
in a finite extension
of the function field,
we may assume
${\cal H}_U$
is a direct sum of
rank one sheaves defined
by the Artin-Schreier
equations
$T^p-T=f$
for linear forms $f$
on $E_U$.
Thus, it is reduced
to the following Lemma.

\begin{lm}\label{lmadd}
Let $K$ be a discrete valuation field
of characteristic $p>0$
and we consider the valuation $v_L$
of $L=K(t_1,\ldots,t_n)$
defined by the the prime ideal
${\mathfrak m}_K\cdot {\cal O}_K[t_1,\ldots,t_n]$
of the polynomial ring.
Then, for a linear form $f
\in Kt_1+\cdots+Kt_n
\subset L$,
the Artin-Schreier extension
of $L$
defined by
$T^p-T=f$ is unramified with respect to
$v_L$ if and only if
$f\in {\cal O}_Kt_1+\cdots+{\cal O}_Kt_n$.
\end{lm}
{\it Proof.}
It suffices to show that 
the Artin-Schreier extension
is ramified assuming $v_L(f)=-n<0$.
If $p\nmid n$,
it is a totally ramified extension.
If $p|n$,
the residue field extension
is the purely inseparable extension
generated by the $p$-th root of
the non-zero linear form
$\overline{\pi^nf}$.
\qed

We introduce a generalization of
vector bundles.
\begin{df}\label{dfbld}

Let $X$ be a scheme
and ${\cal L}$
and ${\cal E}$
be an invertible
${\cal O}_X$-module
and a locally free
${\cal O}_X$-module of
finite rank respectively
and $n\ge 1$ be
an integer.
We call the vector bundloid
of degree $n$ associated
to $({\cal E},{\cal L})$
the affine $X$-scheme
$$E={\mathbf V}_n
({\cal E},{\cal L})$$
defined by the quasi-coherent ${\cal O}_X$-algebra
$\bigoplus_{l\ge 0}
S^{nl}{\cal E}
\otimes 
{\cal L}^{\otimes l}$.
We call 
$E^\vee={\mathbf V}_n
({\cal E}^\vee,{\cal L}^\vee)$
the dual of
$E$.
\end{df}

The grading defines 
a natural action of
the multiplicative group
${\mathbf G}_m$
on ${\mathbf V}_n
({\cal E},{\cal L})$.
For $n=1$,
we have
${\mathbf V}_1
({\cal E},{\cal L})=
{\mathbf V}
({\cal E}\otimes {\cal L})$.
For $m=nr$,
the inclusion
$\bigoplus_{l\ge 0}
S^{nrl}{\cal E}
\otimes 
{\cal L}^{\otimes rl}
\subset
\bigoplus_{l\ge 0}
S^{nl}{\cal E}
\otimes 
{\cal L}^{\otimes l}$
defines a finite surjection
$$\begin{CD}
\pi_{mn}:
{\mathbf V}_n
({\cal E},{\cal L})
@>>> 
{\mathbf V}_m
({\cal E},{\cal L}^{\otimes r}).
\end{CD}$$
It induces
an isomorphism
${\mathbf V}_n
({\cal E},{\cal L})/\mu_r
\to
{\mathbf V}_m
({\cal E},{\cal L}^{\otimes r})$
with respect to
the action restricted
to the group $\mu_r
\subset {\mathbf G}_m$
of $r$-th roots of unity.
If $X$ is a scheme over
${\mathbb F}_p$
and if $r$ is a power of
$p$,
the map
$\pi_{mn}:
{\mathbf V}_n
({\cal E},{\cal L})
\to
{\mathbf V}_m
({\cal E},{\cal L}^{\otimes r})$
induces an isomorphism
on the \'etale site.
If ${\cal L}\to {\cal O}_X$
is an isomorphism,
the map $\pi_{n1}$
defines a finite surjection
${\mathbf V}({\cal E})=
{\mathbf V}_1({\cal E},{\cal O}_X)\to
{\mathbf V}_n({\cal E},{\cal O}_X)\to
{\mathbf V}_n({\cal E},{\cal L})$.
If ${\cal E}={\cal O}_X$,
we have
${\mathbf V}_n
({\cal O}_X,{\cal L})=
{\mathbf V}({\cal L})$.
For a vector bundloid $E=
{\mathbf V}_n
({\cal E},{\cal L})$,
we call
$E^\vee={\mathbf V}_n
({\cal E}^\vee,{\cal L}^\vee)$ 
the dual of $E$.

We call the section
$X\to E={\mathbf V}_n
({\cal E},{\cal L})$
defined by the augmentation
$\bigoplus_{l\ge 0}
S^{nl}{\cal E}
\otimes 
{\cal L}^{\otimes l}
\to {\cal O}_X$
the 0-section of $E$.
We identify $X$ with
a closed subscheme of $E$
by the 0-section.
On the complement
$E^0=E\setminus X$
of the 0-section,
we have a natural map
$$\begin{CD}
\varphi:
E^0
@>>> {\mathbf P}({\cal E})
={\cal P}roj(S^\bullet {\cal E})
\end{CD}$$
since ${\mathbf P}({\cal E})$
is canonically identified with
${\cal P}roj
(\bigoplus_{l\ge 0}
S^{nl}{\cal E}
\otimes 
{\cal L}^{\otimes l})$.
It induces an isomorphism
$E^0/{\mathbf G}_m
\to {\mathbf P}({\cal E}).$
The finite map
$\pi_{mn}:
{\mathbf V}_n
({\cal E},{\cal L})
\to 
{\mathbf V}_m
({\cal E},{\cal L}^{\otimes r})$
is compatible 
with the map
$\varphi:
E^0\to {\mathbf P}({\cal E}).$

\begin{lm}\label{lmOn}
Let ${\cal O}(n)$ 
be the tautological sheaf
on ${\mathbf P}({\cal E})$.
Then, there exists
a canonical isomorphism
$\varphi^*{\cal O}(n)
\to {\cal L}^\vee$
on $E^0$.
\end{lm}
{\it Proof.}
The invertible sheaf ${\cal O}(n)$
on ${\mathbf P}({\cal E})$
is the pull-back of
${\cal O}(1)$ on
${\mathbf P}(S^n{\cal E})$ by
the Veronese embedding
${\mathbf P}({\cal E})
\to
{\mathbf P}(S^n{\cal E})$.
On $E={\mathbf V}_n({\cal E},{\cal L})$,
we have a tautological map
$p^*(S^n{\cal E}\otimes {\cal L})
\to {\cal O}_E$.
On $E^0$,
this is a surjection
and defines a surjection
$p^*S^n{\cal E}\to
p^*{\cal L}^\vee$.
Since the composition
$E^0\to 
{\mathbf P}({\cal E})
\to
{\mathbf P}(S^n{\cal E})$
is defined by the surjection
$p^*S^n{\cal E}\to
p^*{\cal L}^\vee$,
the assertion follows.
\qed

\begin{lm}\label{lmpadd}
Let
$E={\mathbf V}_n({\cal E}, {\cal L})
\to X$
be a vector bundloid on 
a scheme $X$ over $k$.
Let ${\cal M}$
be an invertible ${\cal O}_X$-module
and ${\cal L}\to {\cal M}^{\otimes n}$
be an isomorphism
and let
$\pi:\widetilde E=
{\mathbf V}({\cal E}\otimes {\cal M})
={\mathbf V}_1({\cal E}, {\cal M})
\to E={\mathbf V}_n({\cal E}, {\cal L})$
be the induced map.
Let $g:E^0=E\setminus X\to E$
and $\tilde g:
\widetilde E^0=
\widetilde E\setminus X\to 
\widetilde E$ be the open immersions
of the complements of the $0$-section
and $\pi^0:\widetilde E^0\to E^0$
be the restriction of $\pi:
\widetilde E\to E$.
We consider the following condition
on a constructible sheaf ${\cal H}$
on $E$:

{\rm (P)} The canonical map
${\cal H}\to g_*g^*{\cal H}$ is 
an isomorphism
and the sheaf
$\tilde g_*\pi^{0*}g^*{\cal H}$
on $\widetilde E$
is additive.

1.
Let ${\cal M}'$
be another invertible ${\cal O}_X$-module
and ${\cal L}\to {\cal M}^{\prime\otimes n}$
be an isomorphism.
We define $\pi':\widetilde E'\to E$
etc.\ as above.
Then the condition {\rm (P)}
for ${\cal H}$
with respect to
$\pi:\widetilde E\to E$
is equivalent to 
that for
$\pi':\widetilde E'\to E$.

Assume ${\cal H}$
satisfies the equivalent conditions
and put
$\widetilde {\cal H}=
\tilde g_*\pi^{0*}g^*{\cal H}$
on $\widetilde E$
and
$\widetilde {\cal H}'=
\tilde g'_*\pi^{\prime 0*}g^{\prime *}{\cal H}$
on $\widetilde E'$.
Then 
$S_{\widetilde {\cal H}}\subset 
\widetilde E^\vee$
and
$S_{\widetilde {\cal H}'}\subset 
\widetilde E^{\prime \vee}$
have the same images in $E^\vee$.

2. Let $n'$ be the prime-to-$p$
part of $n$ and assume 
that $k$ contains
a primitive $n'$-th root of $1$.
We consider the natural action of
$G=\mu_{n'}$ on $\widetilde E$ over $E$.
Then, 
the condition {\rm (P)}
for ${\cal H}$
is equivalent to the following
condition:

{\rm (P$'$)} There exist
an additive constructible sheaf 
$\widetilde{\cal H}$
on $\widetilde E$ with an action of
$G$ and an isomorphism
${\cal H}\to 
(\pi_*\widetilde{\cal H})^G$.
\end{lm}
{\it Proof.}
1. The assertion is flat local on $X$.
Hence, we may assume
there exists an isomorphism
${\cal M}\to {\cal M}'$
compatible with
${\cal L}\to {\cal M}^{\otimes n}$
and ${\cal L}\to {\cal M}^{\prime\otimes n}$.
Then the assertion is clear.

2. On the restriction on $E^0$,
the canonical map
$g^*{\cal H}
\to (\pi^0_*\pi^{0*}g^*{\cal H})^G$
to the $G$-fixed part is an isomorphism.
Hence, it induces an isomorphism
\begin{equation}
g_*g^*{\cal H}
\to (g_*\pi^0_*\pi^{0*}g^*{\cal H})^G
\to
(\pi_*\tilde g_*\pi^{0*}g^*{\cal H})^G.
\label{eqpadd}
\end{equation}

(P)$\Rightarrow$(P$'$)
We put
$\widetilde {\cal H}=\tilde g_*\pi^{0*}g^*{\cal H}$.
Then, if the canonical map
${\cal H}\to g_*g^*{\cal H}$ is 
an isomorphism,
we obtain an isomorphism
${\cal H}
\to (\pi_*\widetilde {\cal H})^G$
by the isomorphism (\ref{eqpadd}).

(P$'$)$\Rightarrow$(P).
Let ${\cal H}\to 
(\pi_*\widetilde{\cal H})^G$
be an isomorphism.
Then, it induces an isomorphism
$\pi^{0*}g^*{\cal H}
\to \tilde g^*
\widetilde{\cal H}$
compatible with the $G$-action.
Since $\widetilde{\cal H}$
is additive,
the canonical map
$\widetilde{\cal H}
\to \tilde g_*
\tilde g^*
\widetilde{\cal H}$
is an isomorphism
by \cite{acyc} Proposition 3.2.
Hence
the isomorphism
$\pi^{0*}g^*{\cal H}
\to \tilde g^*
\widetilde{\cal H}$
is extended to an isomorphism
$\tilde g_*\pi^{0*}g^*{\cal H}
\to \widetilde{\cal H}$
and
$\tilde g_*\pi^{0*}g^*{\cal H}$
is additive.
By the isomorphism
(\ref{eqpadd}),
the isomorphism
${\cal H}\to 
(\pi_*\widetilde{\cal H})^G$
implies that 
the canonical map
${\cal H}\to g_*g^*{\cal H}$ is 
an isomorphism.
\qed

We generalize
the notion of
additive sheaves
on vector bundloids.

\begin{df}\label{dfpadd}
Let $E=
{\mathbf V}_n({\cal E},{\cal L})
\to X$
be a vector bundloid
of degree $n$
over a scheme $X$
over $k$.
We say a constructible
sheaf ${\cal H}$ on $E$
is potentially additive
if it satisfies
the condition {\rm (P)}
in Lemma {\rm \ref{lmpadd}}
Zariski locally on $X$.

Let ${\cal H}$ be
a potentially additive constructible
sheaf on $E$.
Then, we define
a constructible subset
$S_{\cal H}$ of the dual $E^\vee$
as the image
of $S_{\widetilde {\cal H}}$
in the notation of
Lemma {\rm \ref{lmpadd}.2}
Zariski locally on $X$
and call $S_{\cal H}$
the dual support of ${\cal H}$.
We say a potentially additive constructible
sheaf ${\cal H}$ on $E$
is non-degenerate
if the intersection of
$S_{\cal H}\subset E^\vee$ with
the $0$-section is empty.
\end{df}

\begin{lm}\label{lmnd}
Let $p:E=
{\mathbf V}_n({\cal E},{\cal L})
\to X$
be a vector bundloid
of degree $n$
over a scheme $X$
over $k$.
Let ${\cal H}$ be 
a potentially additive
constructible ${\mathbb Q}_\ell$-sheaf on $E$.
If it is non-degenerate,
then we have
$Rp_!{\cal H}=
Rp_*{\cal H}=0$.
\end{lm}
{\it Proof.}
Since the assertion is Zariski local
on $X$,
we may use the notation in
Lemma {\rm \ref{lmpadd}}.
Let $\tilde p:\widetilde E\to X$
denote the structural map.
Since ${\cal H}$ is assumed non-degenerate,
we have
$R\tilde p_!\widetilde {\cal H}=
R\tilde p_*\widetilde {\cal H}=0$
by Lemma \ref{lmadd2}.2.
Therefore,
$Rp_!{\cal H}=
(R\tilde p_!\widetilde {\cal H})^G$
and
$Rp_*{\cal H}=
(R\tilde p_*\widetilde {\cal H})^G$
are 0.
\qed

\subsection{Global basic construction}

We study the basic construction in \S1.1
in a global setting.
Let $X$ be a smooth scheme over $k$,
$D$ be a divisor with simple normal
crossings and $j:U=X\setminus D\to X$
be the open immersion of
the complement.
Let $p:P\to X$ be a smooth morphism
of relative dimension $d$
and $s:X\to P$ be a section.
By the section $s$,
we regard $X$ as a closed
subscheme of $P$.

Let $D_1,\ldots,D_m$ be
the irreducible components of $D$.
We consider an effective divisor
$R=r_1D_1+\cdots+r_mD_m$
with rational coefficients
$r_1,\ldots,r_m\ge0$.
For an integer $l\ge 0$,
let $[lR]$ denote
the integral part of $lR$ and
${\cal I}_{[lR]}\subset {\cal O}_X$ be the 
ideal sheaf of the 
effective divisor $[lR]$.
Let ${\cal I}_X\subset {\cal O}_P$
be the ideal sheaf of $X\subset P$
and $j_P:P_U=P\times_XU
\to P$ be the open immersion.
We define 
an affine $P$-scheme
$q:P^{(R)}\to P$
by the quasi-coherent
${\cal O}_P$-algebra
\begin{equation}
\sum_{l\ge 0}p^*{\cal I}_{[lR]}^{-1}\cdot {\cal I}_X^l
\subset j_{P*}{\cal O}_{P_U}.
\label{eqPR}
\end{equation}
Let $p^{(R)}:P^{(R)}\to X$
be the canonical map
and $s^{(R)}:X\to P^{(R)}$
be the section induced by
$s:X\to P$.
We also regard
$D\subset X$ as 
closed subschemes of $P^{(R)}$
by the section $s^{(R)}$.

Here is an alternative
construction of
$q:P^{(R)}\to P$.
Let $n>0$ be an integer
such that
$M=nR$ has integral coefficients.
Let $\bar q:P^{[M/n]}\to P$
be the blow-up by the ideal
$p^*{\cal I}_M+{\cal I}_X^n\subset {\cal O}_P$
and 
$P^{(M/n)}\subset P^{[M/n]}$
be the complement
of the support of
$\bar q^*(p^*{\cal I}_M+{\cal I}_X^n)/
\bar q^*p^*{\cal I}_M$.
The morphism $P^{(M/n)}\to P$
is affine and $P^{(M/n)}$ is defined
by the quasi-coherent
sub ${\cal O}_P$-algebra
${\cal O}_P[p^*{\cal I}_M^{-1}\cdot {\cal I}_X^n]
\subset j_{P*}{\cal O}_{P_U}$.
Then, 
$P^{(R)}$ is identified with the
normalization of
$P^{(M/n)}$.

We put $I^+=\{i|1\le i\le m,r_i>0\}$
and $D^+=\sum_{i\in I^+}D_i$.
We describe the structure
of the inverse image
$E^+=P^{(R)}\times_XD^+$
in terms of vector bundloids
introduced in the previous subsection.

\begin{lm}\label{lmR}
Let $D_1,\ldots,D_m$
be the irreducible components of $D$
and put $I^+=\{i|1\le i\le m,r_i>0\}$.

1. 
Let $I\subset I^+$
be a non-empty subset
and $n_I\ge 1$ be
the minimum integer
such that the coefficients
in $n_IR$ of $D_i$
are integers for all $i\in I$. 
Let $D_I$ be the intersection
$\bigcap_{i\in I}D_i$
and put $D_I^\circ
=D_I\setminus 
\bigcup_{i\in I^+\setminus I}
(D_i\cap D_I)$.

Then, there exists a canonical isomorphism
$$\begin{CD}
E_I^{\circ}
=(P^{(R)}\times_XD_I^\circ)_{\rm red}
@>>>
{\mathbf V}_{n_I}
({\cal N}_{X/P},{\cal O}(n_IR))\times_X
D_I^\circ
\end{CD}$$
over $D_I^\circ$.
The restriction 
$D_I^\circ\to E_I^{\circ}$
of the section
$s^{(R)}:X\to P^{(R)}$
corresponds
to the $0$-section of
the right hand side.

2. 
Let 
$R^*=X\times_P(P^{(R)}\setminus X)$ 
be the inverse image
of $X=s(X)\subset P$ 
by the restriction 
of the canonical map
$q:P^{(R)}\to P$
on the complement
$P^{(R)}\setminus X$
of the section $s^{(R)}$.
Then $R^*$ is a divisor of
$P^{(R)}\setminus X$
and satisfies
$R^*=p^{(R)*}R$.

3. Assume the coefficients of $R$
are integers.
Then, the map $p^{(R)}:P^{(R)}\to P$
is smooth.
The inverse image
$E^+=P^{(R)}\times_XD^+$
of $D^+=\sum_{i\in I^+}D_i$
is canonically isomorphic
to the vector bundle
${\mathbf V}
({\cal N}_{X/P}\otimes {\cal O}(R))\times_XD^+$.
\end{lm}
{\it Proof.}
1.
We may assume $I=I^+$.
By the definition (\ref{eqPR})
of $P^{(R)}$,
we have a surjection
$$\bigoplus_{l\ge 0,n_I|l}({\cal O}(lR)\otimes 
S^l {\cal N}_{X/P})\otimes_{{\cal O}_X} {\cal O}_{D_I^\circ}
\to {\cal O}_{E_I^\circ}.$$
Namely, we have a closed immersion
$E_I^\circ
\to 
{\mathbf V}_{n_I}
({\cal N}_{X/P},{\cal O}(n_IR))\times_X
D_I^\circ$.
We show this is an isomorphism.
Since the question is 
\'etale local on $P$,
we may assume
$P={\mathbf V}({\cal E})$
is a vector bundle defined by
a locally free ${\cal O}_X$-module ${\cal E}$
of rank $d$
and $s:X\to P$
is the 0-section.
Then $P^{(R)}$ is the affine scheme over
$X$ defined by the ${\cal O}_X$-algebra
$\bigoplus_{l\ge 0}
S^l{\cal E}\otimes {\cal O}([lR])$.
Hence the assertion follows.

2.
Let $n>0$ be an integer
such that $M=nR$ has 
integral coefficients.
Since the question is 
local on $P$,
we may assume the ideal
${\cal I}_X\subset {\cal O}_P$
is generated by $d$ sections
$e_1,\ldots,e_d$
and ${\cal I}_{nR}$ has a basis $l$.
Then, on the open subscheme
of $P^{(R)}$ where
$f_i=l^{-1}e_i^n$ is invertible,
the pull-back of the ideal
${\cal I}_X=(e_1,\ldots,e_d)$
is generated by $e_i$
since $e_j=e_i
\cdot l^{-1}e_je_i^{n-1}/f_i$.
Since the support of
the closed subscheme of $P^{(R)}$
defined by the ideal $(f_1,\ldots,
f_d)$ is $s^{(R)}(X)$,
the assertion follows.

3.
We show that the scheme
$P^{(R)}$ is smooth over $X$.
Since the question is
\'etale local on $P$,
we may assume
$P={\mathbf V}({\cal E})$
is a vector bundle defined by
a locally free ${\cal O}_X$-module ${\cal E}$
of rank $d$
and $s:X\to P$
is the 0-section as in the proof of 1.
Then $P^{(R)}$ is the vector bundle
${\mathbf V}({\cal E}\otimes {\cal O}(R))$
and the assertion follows.

Similarly as in the proof of 1,
we obtain a closed immersion
$E^+=E\times_XD^+
\to 
{\mathbf V}({\cal N}_{X/P}\otimes {\cal O}(R))
\times_XD^+$
and we see that
this is an isomorphism.
\qed

We have the following functoriality of
the construction of $P^{(R)}$.

\begin{lm}\label{lmPRf}
We consider a commutative
diagram 
$$\begin{CD}
Y@>t>> Q@>q>> Y\\
@VfVV @VgVV @VVfV\\
X@>s>>P@>p>>X
\end{CD}$$
of smooth schemes over $k$.
We assume that
$s:X\to P$ and
$t:Y\to Q$
are sections of
smooth maps $p:P\to X$
and $q:Q\to Y$
respectively.
Let $D$ be a divisor of $X$
with simple normal crossings.
Assume that
the divisor $D_Y=
(D\times_XY)_{\rm red}$
has simple normal crossings.
Let $R=\sum_ir_iD_i\ge 0$
be an effective divisor
with rational coefficients
$r_i\ge 0$
and let
$R_Y=f^*R$
be the pull-back.

1.
There exists
a unique map
$g^{(R)}:Q^{(R_Y)}\to P^{(R)}$
lifting $g:Q\to P$.

2. 
Suppose that the coefficients
of $R$ are integral.
Let $D^+$
and $D_Y^+$
be the supports
of $R$ and of $R_Y$
respectively.
We identify
$E^+=
P^{(R)}\times_XD^+$
with
${\mathbf V}(
{\cal N}_{X/P}\otimes {\cal O}(R))
\times_XD^+$
and
$E_Y^+=
Q^{(R_Y)}\times_YD_Y^+$
with
${\mathbf V}(
{\cal N}_{Y/Q}\otimes {\cal O}(R_Y))
\times_YD_Y^+$
as in Lemma {\rm \ref{lmR}.1}.
Then the restriction
$$\begin{CD}
E_Y^+=
{\mathbf V}(
{\cal N}_{Y/Q}\otimes {\cal O}(R_Y))
\times_YD_Y^+
@>>>
E^+=
{\mathbf V}(
{\cal N}_{X/P}\otimes {\cal O}(R))
\times_XD^+
\end{CD}$$
of 
$g^{(R)}:Q^{(R_Y)}\to P^{(R)}$
is the linear map
of vector bundles
induced by the canonical map
$f^*{\cal N}_{X/P}\to {\cal N}_{Y/Q}$.

3.
Suppose further that $f:Y\to X$
is the identity of $X$
and $g:Q\to P$
is smooth.
Then the induced map
$g^{(R)}:Q^{(R_Y)}\to P^{(R)}$
is smooth.
\end{lm}
{\it Proof.}
1. We have
$g^* {\cal I}_X\subset {\cal I}_Y$
since the left square is commutative.
By the inequalities
$g^*[lR]\le
[lR_Y]\le lg^*R$,
we have
$g^*{\cal I}_{[lR]}^{-1}
\subset {\cal I}_{[lR_Y]}^{-1}$.
Hence we have
$g^*({\cal I}_{[lR]}^{-1}\cdot
{\cal I}_X^l)\subset 
{\cal I}_{[lR_Y]}^{-1}
\cdot {\cal I}_Y^l$ and
the assertion follows
from the definition of $Q^{(R_Y)}$.

2.
The restriction
$E_Y^+\to E^+$
is induced by the linear map
$g^*:{\cal I}_R^{-1}\cdot
{\cal I}_X\to
{\cal I}_{R_Y}^{-1}
\cdot {\cal I}_Y$
and the assertion follows.

3.
On the complements of the inverse images
of $D^+$,
the maps
$P^{(R)}\to P$
and 
$Q^{(R)}\to Q$
are isomorphisms.
Hence, the assumption that
$f:Q\to P$
is smooth implies that
the restriction on
the complements of the inverse images
of $D^+$
is smooth.
Since $P^{(R)}\to X$
and 
$Q^{(R)}\to X$
are smooth,
it suffices to show
that the induced map
$E^{+\prime}=
Q^{(R)}\times_XD^+
\to E^+=
P^{(R)}\times_XD^+$
is smooth.

By 2,
it is identified
with the map
${\mathbf V}(
{\cal N}_{X/Q}\otimes {\cal O}(R))
\times_XD^+
\to
{\mathbf V}(
{\cal N}_{X/P}\otimes {\cal O}(R))
\times_XD^+$
of vector bundles
induced by the canonical map
$f^*{\cal N}_{X/P}\to {\cal N}_{X/Q}$.
Since $Q\to P$
is assumed smooth,
the map
$f^*{\cal N}_{X/P}\to {\cal N}_{X/Q}$
is a locally splitting injection
and the map
${\mathbf V}(
{\cal N}_{X/Q}\otimes {\cal O}(R))
\times_XD^+
\to
{\mathbf V}(
{\cal N}_{X/P}\otimes {\cal O}(R))
\times_XD^+$
is smooth.
\qed

\begin{cor}\label{corPRf}
Let $P$ and $Q$ be smooth schemes
over $X$ and
$s:X\to P$
and $t:X\to Q$
be sections.
Similarly as
$P^{(R)}$ and
$Q^{(R)}$,
we define
$(P\times_X Q)^{(R)}$
by the section
$(s,t):X\to P\times_XQ$.

Assume the coefficients
of $R$ are integers.
Then the maps
$(P\times_X Q)^{(R)}
\to P^{(R)}$ and
$(P\times_X Q)^{(R)}
\to Q^{(R)}$
induces an isomorphism
\begin{equation}
\begin{CD}
(P\times_X Q)^{(R)}
@>>>
P^{(R)}\times_X Q^{(R)}.
\end{CD}
\label{eqPQR}
\end{equation}
\end{cor}
{\it Proof.}
The ideal
defining the closed subscheme
$X\subset P\times_X Q$
is generated by
the pull-backs of
those defining
$X\subset P$
and
$X\subset Q$.
Hence the map
(\ref{eqPQR})
is a closed immersion.
Since the both 
schemes
$(P\times_X Q)^{(R)}$
and
$P^{(R)}\times_X Q^{(R)}$
are smooth of the same dimension
over $X$,
the closed immersion
(\ref{eqPQR})
is an open immersion.
By Lemma \ref{lmR}.3,
it induces an isomorphism
on the fibers over $D^+$.
Hence the assertion follows.
\qed

We establish some cohomological properties
of $P^{(R)}$.

\begin{pr}\label{prPr}
1. The cycle class defines
an isomorphism
\begin{equation}
\begin{CD}
{\mathbb Q}_\ell(d)[2d]
@>>>
 Rp^{(R)!}{\mathbb Q}_\ell.
\end{CD}
\label{eqPrK}
\end{equation}

2. Define the cycle class 
$[X]\in H^{2d}_X(
P^{(R)},{\mathbb Q}_\ell(d))$
to be the inverse image of
$1\in H^0(X,{\mathbb Q}_\ell)
=H^0_X(P^{(R)},
Rp^{(R)!}{\mathbb Q}_\ell)$
by the isomorphism
{\rm (\ref{eqPrK})}.
Then, for 
the pull-back $s^{(R)*}[X]
=(X,X)_{P^{(R)}}
\in H^{2d}(X,{\mathbb Q}_\ell(d))$,
we have
\begin{eqnarray}
(X,X)_{P^{(R)}}
&=&(X,X)_P
-(c({\cal N}_{X/P})^*\cap (1+R)^{-1}\cap [R])_{\deg d}
\label{eqPr0}
\\
&=&
(-1)^d
(c_d({\cal N}_{X/P})+
(c({\cal N}_{X/P})\cap (1-R)^{-1}\cap [R])_{\deg d}.
\label{eqPr1}
\end{eqnarray}
\end{pr}
{\it Proof.}
1.
Since the question is \'etale local,
we may assume there exists 
a smooth map $X\to {\mathbf A}^m_k$
such that
$D$ is the inverse image
of the union of
coordinate hyperplanes.
Let $n_1,\ldots,n_m>0$
be integers
and let $\pi:\widetilde X=
X\times_{{\mathbf A}^m_k}
{\mathbf A}^m_k\to X$
be the base change by
the map ${\mathbf A}^m_k
\to {\mathbf A}^m_k$
defined by
$t_i\mapsto t_i^{n_i}$.
We put $\widetilde P
=P\times_X\widetilde X$
and
$\widetilde R=
\pi^*R$.
We consider the commutative diagram
\begin{equation}
\begin{CD}
P^{(R)}@<{\pi}<< \widetilde P^{(\widetilde R)}\\
@V{p^{(R)}}VV @VV{\tilde p^{(\widetilde R)}}V \\
X@<<< \widetilde X.
\end{CD}
\label{eqPRsim}
\end{equation}
The map 
$\tilde p^{(\widetilde R)}:
\widetilde P^{(\widetilde R)}
\to \widetilde X$
is smooth if
$n_1r_1,\ldots, n_mr_m$
are integers
by Lemma \ref{lmR}.
Let $n'_1,\ldots, n'_m$
be the prime-to-$p$ parts of
$n_1,\ldots, n_m$
and put $G=
\mu_{n'_1}\times \cdots\times \mu_{n'_m}$.
Then, we have a natural action of
$G$ on 
$\widetilde P^{(\widetilde R)}$.
The induced map
$\widetilde P^{(\widetilde R)}/G\to
P^{(R)}$ defines an isomorphism
on the \'etale sites.

We put ${\cal K}_{P^{(R)}}=
Rp^{(R)!}{\mathbb Q}_\ell$
and ${\cal K}_{\widetilde P^{
(\widetilde R)}}=
R(p^{(R)}\circ \pi)^!
{\mathbb Q}_\ell$.
The trace map
$\pi_!\pi^*{\cal K}_{P^{(R)}}
=\pi_*\pi^*{\cal K}_{P^{(R)}}
\to {\cal K}_{P^{(R)}}$
defines its adjoints
$\pi^*{\cal K}_{P^{(R)}}\to 
\pi^!{\cal K}_{P^{(R)}}={\cal K}_{\widetilde P^{
(\widetilde R)}}$
and
${\cal K}_{P^{(R)}}\to \pi_*{\cal K}_{\widetilde P^{
(\widetilde R)}}$.
We also have the adjunction map
$\pi_*{\cal K}_{\widetilde P^{
(\widetilde R)}}=
\pi_!\pi^!{\cal K}_{P^{(R)}}\to {\cal K}_{P^{(R)}}$.
The composition
${\cal K}_{P^{(R)}}\to 
\pi_*{\cal K}_{\widetilde P^{
(\widetilde R)}}\to {\cal K}_{P^{(R)}}$
is the multipication 
by the degree $[\widetilde P^{
(\widetilde R)}:P^{(R)}]$.
Hence ${\cal K}_{P^{(R)}}$
is a direct summand of
the $G$-fixed part
$(\pi_*{\cal K}_{\widetilde P^{
(\widetilde R)}})^G.$

We consider the commutative diagram
$$\begin{CD}
{\mathbb Q}_\ell(d)[2d]
@>>>
{\cal K}_{P^{(R)}}\\
@VVV @VVV\\
(\pi_*{\mathbb Q}_\ell(d)[2d])^G
@>>>
(\pi_*{\cal K}_{\widetilde P^{
(\widetilde R)}})^G
\end{CD}
$$
where the horizontal arrows 
are defined by the cycle classes.
Since $\widetilde P^{
(\widetilde R)}$ is smooth
over $k$,
the lower horizontal arrow
is an isomorphism.
Since the left vertical arrow
is an isomorphism,
$(\pi_*{\cal K}_{\widetilde P^{
(\widetilde R)}})^G$
is a direct summand of
${\cal K}_{P^{(R)}}$. 
Thus the assertion is proved.

2.
First, we reduce it to the case
where $P$ is a vector bundle over $X$
and $s:X\to P$ is the $0$-section,
by the deformation to the normal bundle.
We put $\widetilde X=X\times {\mathbf A}^1$
and $\widetilde D=D\times {\mathbf A}^1$.
Let $\overline P$
be the blow-up of $P\times {\mathbf A}^1$
at $X\times \{0\}$
and $\widetilde P\subset
\overline P$
be the complement of the
proper transform of $P\times \{0\}$.
Then, the map
$\tilde p:
\widetilde P\to 
\widetilde X$ is smooth.
We consider the cartesian diagram
$$\begin{CD}
V@>>>\widetilde P
@<<<P\times {\mathbb G}_m\\
@VVV @V{\tilde p}VV @VV{p\times {\rm id}}V\\
X @>>> \widetilde X@<<<
X\times {\mathbb G}_m
\\
@VVV @VVV @VVV\\
\{0\}@>>> {\mathbf A}^1
@<<<{\mathbb G}_m
\end{CD}$$
where $V={\mathbf V}
({\cal N}_{X/P})$ denotes
the normal bundle.

The section $s:X\to P$
induces a section
$\tilde s:\widetilde X
\to \widetilde P$.
By applying the construction,
we define
$\tilde p^{(R)}:
\widetilde P^{(R)}
\to \widetilde X$
and its section
$\tilde s^{(R)}:
\widetilde X
\to \widetilde P^{(R)}$.
By the standard argument,
the assertion for $(X,V)$
is equivalent to
that for $(\widetilde X,
\widetilde P)$
and implies that for
$(X,P)$.
Thus we may assume 
$P$ is a vector bundle over $X$
and $s:X\to P$ is the $0$-section.

Let $q:P^{(R)}\to P$ denote
the canonical map.
It suffices to
show the equality
$$[s^{(R)}(X)]
=q^*[s(X)]
-p^{(R)*}((c({\cal N}_{X/P})^*\cap 
(1+R)^{-1}
\cap 
[R])_{\deg d}$$
in 
$H^{2d}_{q^{-1}(s(X))}
(P^{(R)},{\mathbb Q}_\ell(d))$.
For a closed subscheme
$F$ of $P^{(R)}$,
let $F^\circ$ denote
the complement
$F\setminus (F\cap D^+)$.
Let $R^*$ be the divisor
of $P^{(R)}\setminus X$
defined in Lemma \ref{lmR}.2.
Then, we have
$q^{-1}(s(X))
=s^{(R)}(X)\cup
\overline R^*$ and
$s^{(R)}(X)\cap
\overline R^*=D^+$.
Hence, by 1, we have
$H^{2d}_{s^{(R)}(X)\cap
\overline R^*}
(P^{(R)},{\mathbb Q}_\ell(d))=
H^{2d}_{D^+}
(X,s^{(R)!}{\mathbb Q}_\ell(d))=
H^0_{D^+}
(X,{\mathbb Q}_\ell)=0$.
Thus, the restriction map
$$H^{2d}_{q^{-1}(s(X))}
(P^{(R)},{\mathbb Q}_\ell(d))
\to
H^{2d}_{X^\circ}
(P^{(R)\circ},
{\mathbb Q}_\ell(d))\oplus
H^{2d}_{R^*}
(P^{(R)\circ },
{\mathbb Q}_\ell(d))
$$
is an injection.
Therefore, it suffices to show that
the components of the
restriction of
$q^*[s(X)]$ are
$[s^{(R)}(X)^\circ]$
and 
$p^{(R)*}(c({\cal N}_{X/P})^*\cap (1+R)^{-1}\cap 
[R])_{\deg d}$
respectively.

This is clear for the first component
$[s^{(R)}(X)^\circ]$.
By the excess intersection formula,
the second component is
$(c({\cal N}_{X/P})^*\cap 
c({\cal N}_{R^*/P^{(R)\circ}})^{*-1}
\cap 
[R^*])_{\deg d}.$
Hence,
the assertion follows
by Lemma \ref{lmR}.2.
\qed

\subsection{Ramification along 
a divisor}

We globalize the constructions
in \S\S1.1 and 1.2
and the computations in \S1.4.
Let $X$ be a
smooth scheme of dimension $d$
over $k$
and $D$ be a divisor with
simple normal crossings.
Let $D_1,\ldots,D_m$
be the irreducible components
of $D$.
We put $U=X\setminus D$
and let $j:U\to X$
denote the open immersion.

We define the log blow up
$(X\times X)'
\to X\times X$
to be the blow-up at
$D_1\times D_1,
D_2\times D_2,
\ldots,
D_m\times D_m$.
Namely the blow-up by the product
${\cal I}_{D_1\times D_1}\cdot
{\cal I}_{D_2\times D_2}\cdots
{\cal I}_{D_m\times D_m}
\subset {\cal O}_{X\times X}$
of ideal sheaves.
We define the log product
$(X\times X)^\sim\subset
(X\times X)'$
to be the complement
of the proper transforms
of $D\times X$ and $X\times D$.
The diagonal map
$X\to X\times X$
induces a closed immersion
$\tilde \delta:X\to (X\times X)^\sim
\subset (X\times X)'$
called the log diagonal map.
The scheme $(X\times X)^\sim$
is affine over $X\times X$
and is defined
by the quasi-coherent
${\cal O}_{X\times X}$-algebra
$${\cal O}_{X\times X}[
{\rm pr}_1^*{\cal I}_{D_i}^{-1}\cdot
{\rm pr}_2^*{\cal I}_{D_i},\
{\rm pr}_1^*{\cal I}_{D_i}\cdot
{\rm pr}_2^*{\cal I}_{D_i}^{-1};
i=1,\ldots,m]
\subset
j^\times_*
{\cal O}_{U\times U}$$
where 
$j^\times:U\times U
\to X\times X$
is the open immersion.
The projections
$p_1,p_2:(X\times X)^\sim\to X$ are smooth.
The conormal sheaf
${\cal N}_{X/(X\times X)^\sim}$
is canonically identified with
the locally free ${\cal O}_X$-module 
$\Omega^1_X(\log D)$ of rank $d$.

Let
$R=r_1D_1+\cdots+r_mD_m$
be  an effective divisor
with rational coefficients
$r_1,\ldots,r_m\ge0$.
We apply the construction
of \S2.1
to the smooth map
$p_2:P=(X\times X)^\sim\to X$
and its section
$\tilde \delta:
X\to (X\times X)^\sim$.
Then, we obtain
$P^{(R)}=(X\times X)^{(R)}\to X$
and its section
$\delta^{(R)}:X\to (X\times X)^{(R)}$.
Thus, we have constructed a diagram
$$\begin{CD}
X\times X
@<<< 
(X\times X)'
@<<<
(X\times X)^{[R]}\\
@. @AAA @AAA\\
@.
(X\times X)^\sim
@<<<
(X\times X)^{(R)},
\end{CD}$$
where the vertical arrows are open immersions.
For $R=0$,
we have
$(X\times X)^{(R)}
=(X\times X)^\sim$.

We consider the cartesian diagram
$$\begin{CD}
U\times U
@>{j^{(R)}}>>
(X\times X)^{(R)}\\
@A{\delta_U}AA
@AA{\delta^{(R)}}A\\
U@>j>>X
\end{CD}$$
where the horizontal arrows
are open immersions
and the vertical arrows
are the diagonal immersions.

\begin{df}\label{dfbdR}
Let ${\cal F}$
be a smooth sheaf on $U=
X\setminus D$.
We define a smooth sheaf
${\cal H}$ on $U\times U$
by ${\cal H}=
{\cal H}om(
{\rm pr}_2^*{\cal F},
{\rm pr}_1^*{\cal F})$.
Let $R=\sum_ir_iD_i\ge 0$
be an effective
divisor with rational coefficients
and we consider
the open immersion 
$j^{(R)}:U\times U
\to (X\times X)^{(R)}$.
We identify
$\delta_U^*{\cal H}
={\cal E}nd({\cal F})$
and regard
the identity
${\rm id}_{\cal F}
\in 
{\rm End}_U({\cal F})$
as a section of
$\Gamma(U,{\cal E}nd({\cal F}))=
\Gamma(X,j_*{\cal E}nd({\cal F}))=
\Gamma(X,j_*\delta_U^*{\cal H}).$

We say that the log ramification
of ${\cal F}$ along $D$ is bounded
by $R+$ if
the identity
${\rm id}_{\cal F}\in 
{\rm End}_U({\cal F})=
\Gamma(X,j_*\delta_U^*{\cal H})$
is in the image of the base change
map
\begin{equation}
\begin{CD}
\Gamma(X,\delta^{(R)*}j^{(R)}_*{\cal H})@>>>
\Gamma(X,j_*\delta_U^*{\cal H})
={\rm End}_U({\cal F}).
\end{CD}
\label{eqcbc}
\end{equation}
\end{df}

We compare Definition \ref{dfbdR} 
with Definition \ref{dfFr}.

\begin{lm}\label{lmbc}
Let ${\cal F}$
be a smooth sheaf on $U=
X\setminus D$
and let $R=\sum_ir_iD_i\ge 0$
be an effective
divisor with rational coefficients.
We consider the smooth sheaf
${\cal H}=
{\cal H}om(
{\rm pr}_2^*{\cal F},
{\rm pr}_1^*{\cal F})$ on $U\times U
\subset (X\times X)^{(R)}$.

1. We consider the following conditions:

{\rm (1)}
The log ramification
of ${\cal F}$ along $D$ is bounded
by $R+$.

{\rm (2)}
For every irreducible componenent
$D_i$ of $D$, 
the log ramification
of ${\cal F}$ along $D$ is bounded
by $r_i+$ at the generic point
$\xi_i$ of $D_i$.

{\rm (3)}
There exists an open
subscheme $X'\subset X$
such that $X'\supset U$,
that $D'=X'\cap D$ is dense in $D$
and that
the base change map
$$\begin{CD}
\delta^{(R)*}j^{(R)}_*{\cal H}
@>>>
j_*\delta_U^*{\cal H}
\end{CD}$$
is an isomorphism on
$D'$.

{\rm (4)}
The base change map
$$\begin{CD}
\delta^{(R)*}j^{(R)}_*{\cal H}
@>>>
j_*\delta_U^*{\cal H}
\end{CD}$$
is an isomorphism.

\noindent
Then, we have implications
$(4)\Rightarrow
(1)\Rightarrow
(2)\Rightarrow(3)$.

2.
Let $D_i$
be a component of $D$ 
satisfying $r_i>0$
and let $E_i=(X\times X)^{(R)}
\times_XD_i$
be the inverse image.
Then, the vanishing 
${\cal F}_{\bar \eta_i}^
{G_{K,\log}^{r_i+}}=0$
implies
$j_*{\cal H}|_{E_i}=0$.
\end{lm}
{\it Proof.}
1.
The implication
$(1)\Rightarrow(2)$
follows from
Corollary \ref{corgr0}
$(3)\Rightarrow(1)$.
The implication
$(2)\Rightarrow(3)$
follows from
Corollary \ref{corgr0}
$(1)\Rightarrow(2)$.
The implication
$(4)\Rightarrow(1)$
is obvious.

2. 
Let $\xi_i$ be the generic point of $D_i$.
It suffices to show
$j_*{\cal H}|_{E_{i,\xi_i}}=0$.
Hence, it follows from
Corollary \ref{corgr0}.
\qed

The author does not know a counterexample
for the implication
$(1)\Rightarrow(4)$.
The conditions (1) to (4)
are equivalent,
if the rank of ${\cal F}$ is 1.

In the tamely ramified case,
we have the following equivalence
for $R=0$.

\begin{cor}\label{corbd0}
The following conditions
are equivalent:

{\rm (1)}
The log ramification
of ${\cal F}$ along $D$ is bounded
by $0+$.

{\rm (2)}
${\cal F}$ is tamely ramified along $D$.

{\rm (3)}
The base change map
$$\begin{CD}
\tilde \delta^*\tilde j_*{\cal H}
@>>>
j_*\delta_U^*{\cal H}
\end{CD}$$
is an isomorphism on
$D$.
\end{cor}
{\it Proof.}
By Lemma \ref{lmbc}
(1)$\Rightarrow$(2),
the condition
(1) implies (2).

Assume ${\cal F}$
is tamely ramified along $D$.
Then ${\cal H}$ on $U\times U$
is tamely ramified along 
$(X\times X)^\sim 
\setminus (U\times U)$.
Hence, \'etale locally on
$(X\times X)^\sim$,
it is isomorphic to 
the pull-back of
a sheaf on $U$
with respect to the projection
$(X\times X)^\sim\to X$.
Since the projection
is smooth,
the condition (3) is satisfied.

It is clear that
(3) implies (1).
\qed

We have the following
stability
under the pull-back.

\begin{lm}\label{lmcutY}
Let $Y$ be a smooth scheme over $k$
and $f:Y\to X$
be a morphism over $k$.
Assume that the reduced inverse image
$D_Y=(D\times_XY)_{\rm red}$
is a divisor with simple normal crossings
and let $R_Y$ be the pull-back
$f^*R$.

Let ${\cal F}$
be a smooth sheaf on
$U=X\setminus D$
and ${\cal F}_V$
be the pull-back to
$V=U\times_XY=Y\setminus D_Y$.
If the log ramification of
${\cal F}$ is bounded by $R+$,
then the log ramification of
${\cal F}_V$ is bounded by $R_Y+$.
\end{lm}

{\it Proof.}
We show that
the map $f\times f:
Y\times Y\to X\times X$
is lifted to
$(f\times f)^{(R)}:
(Y\times Y)^{(R_Y)}\to 
(X\times X)^{(R)}$.
For each irreducible component
$D_i$ of $D$,
the pull-backs
of ${\rm pr}_1D_i$
and ${\rm pr}_2D_i$
are equal on
the log product
$(Y\times Y)^\sim$.
Hence, the map $f\times f:
Y\times Y\to X\times X$
is uniquely lifted to
$(f\times f)^\sim:
(Y\times Y)^\sim\to 
(X\times X)^\sim$.
This is uniquely lifted to
$(Y\times Y)^{(R_Y)}\to 
(X\times X)^{(R)}$
by Lemma \ref{lmPRf}.1.

Let $j^{(R)}:U\times U
\to (X\times X)^{(R)}$
and $j^{(R_Y)}:V\times V
\to (Y\times Y)^{(R_Y)}$
be the open immersions
and $f_U:V\to U$
be the restriction of 
$f:Y\to X$.
We put
${\cal F}_V=f_U^*{\cal F}$,
${\cal H}=
{\cal H}om({\rm pr}_2^*{\cal F},
{\rm pr}_1^*{\cal F})$
and 
${\cal H}'=
{\cal H}om({\rm pr}_2^*{\cal F}_V,
{\rm pr}_1^*{\cal F}_V)$.
Then, the base change map
$(f\times f)^{(R)^*}j^{(R)}_*
\to j^{(R_Y)}_*(f_U\times f_U)^*$
defines a commutative diagram
$$\begin{CD}
\Gamma(X,\delta^{(R)*}j^{(R)}_*{\cal H})@>>> 
{\rm End}_U({\cal F})\\
@VVV @VVV\\
\Gamma(Y,\delta^{(R_Y)*}j^{(R_Y)}_*{\cal H}')@>>> 
{\rm End}_V({\cal F}_V)
\end{CD}$$
By the assumption
that the log ramification
of ${\cal F}$
is bounded by $R+$,
the identity of ${\cal F}$
is in the image of the upper horizontal
arrow.
Hence, the identity of ${\cal F}_V$
is in the image of the lower horizontal
arrow
and the log ramification
of ${\cal F}_V$
is bounded by $R_Y+$.
\qed

We consider the restrictions
of ${\cal F}$ on 
smooth curves in $X$
and compare them.

\begin{pr}\label{prjet}
Let ${\cal F}$
be a smooth sheaf on
$U=X\setminus D$
such that
the log ramification of
${\cal F}$ is bounded
by $R+$.
Let $C$ and $C'$ be
smooth curves in $X$
and $x$ be a closed
point in $C\cap C'\cap D$.
We assume that
$C\cap U$ and $C'\cap U$
are not empty
and let
${\cal F}_C$
and 
${\cal F}_{C'}$
denote the restrictions of
${\cal F}$ on
$C\cap U$ and $C'\cap U$
respectively.
Assume that the following
conditions are satisfied:

{\rm (1)}
For every irreducible component
$D_i$ of $D$,
we have
$(C,D_i)_x=
(C',D_i)_x$.

{\rm (2)}
${\rm length}_x
{\cal O}_{C\cap C',x}
\ge (C,R+D)_x$.

\noindent
Then, \'etale locally at $x$, 
there exist an isomorphism
$f:C\to C'$
and an isomorphism
$f^*{\cal F}|_{C'}\to
{\cal F}|_C$.
\end{pr}
{\it Proof.}
It suffices to consider
the case $C\neq C'$.
Since the assertion is \'etale local,
we may assume $C\cap D=
C'\cap D=C\cap C'=\{x\}$
set theoretically
and the residue field
of $x$ is $k$.
We put $n=
{\rm length}_x
{\cal O}_{C\cap C',x}$.
Take an isomorphism
$k[t]/(t^n)
\to {\cal O}_{C\cap C',x}$
and lift it to
\'etale morphisms
$C\to {\mathbf A}^1_k$
and $C'\to {\mathbf A}^1_k$.
Since the assertion is \'etale local,
we may assume
there exists an isomorphism
$f:C\to C'$
inducing the identity
on $C\cap C'$.

We consider the graph of $f$
$$
\begin{CD}
g=(1,f):C@>>>
C\times C'\subset X\times X
\end{CD}.$$
The intersection with 
the diagonal defines
an isomorphism
$C\times_{X\times X}X
\to 
(C\times C')
_{X\times X}X$
$=C\cap C'$
since
$f:C\to C'$
induces the identity
on $C\cap C'$.
By the assumption (1),
the immersion
$g:C\to X\times X$ 
is uniquely lifted 
to an immersion
$\tilde g:C\to (X\times X)^{\sim}$
to the log product.
We put
$C\cap^{\log} C'=
C\times_{
(X\times X)^{\sim}}X
\subset
C\cap C'$.
We show
\begin{equation}
{\rm length}_x
{\cal O}_{C\cap^{\log} C',x}=
{\rm length}_x
{\cal O}_{C\cap C',x}-
(C,D)_x.
\label{eqlg}
\end{equation}
Let ${\cal I}_X\subset 
{\cal O}_{X\times X}$
and ${\cal J}_X\subset 
{\cal O}_{(X\times X)^\sim}$
be the ideal sheaves
of $X\subset X\times X$
and of $X\subset (X\times X)'$
respectively
and let
${\cal I}_E
\subset {\cal O}_{(X\times X)^\sim}$
be the ideal sheaves
of $E=p^*D$.
Then, we have
${\cal I}_X{\cal O}_{(X\times X)^\sim}=
{\cal J}_X\cdot {\cal I}_E$.
By pulling back it by $\tilde g$,
we obtain the equality
(\ref{eqlg}).

By the assumption (2)
and by (\ref{eqlg}),
we have
${\rm length}_x
{\cal O}_{C\times_{(X\times X)^\sim}X,x}
\ge (C,R)_x$.
In other word,
we have
inclusions
${\cal J}_X^l{\cal O}_C
\subset {\cal I}_{[lR]}{\cal O}_C$
for every integer $l\ge 0$.
Hence the immersion
$\tilde g:C\to
(X\times X)^\sim$
is further lifted to
$h:C\to (X\times X)^{(R)}$.
We consider a cartesian
diagram
$$\begin{CD}
C\cap U @>{h_U}>> U\times U\\
@V{j_C}VV @VV{j^{(R)}}V\\
C @>h>> (X\times X)^{(R)}
\end{CD}$$
where the vertical arrows are
open immersions.
We also consider the base change maps
\begin{equation}
\begin{CD}
h^*j^{(R)}_*{\cal H}
@>>>
j_{C*}h_U^*{\cal H}&
=j_{C*}{\cal H}om
(f^*{\cal F}_{C'}
,{\cal F}_C)
,\\
\delta^{(R)*}j^{(R)}_*{\cal H}
@>>>
j_*\delta_U^*{\cal H}&
\!\!\!\!
\!\!\!\!
\!\!\!\!
\!\!\!\!
\!\!\!\!
\!\!\!\!
\!\!\!\!
=j_*{\cal E}nd({\cal F}).
\end{CD}
\label{eqjet}
\end{equation}

Let $K$ denote the fraction
field of the henselization
${\cal O}_{C,x}$
and let $\bar \eta$
denote the geometric point
of $C$ defined
by an algebraic closure
$\overline K$ of $K$.
Let $G_K$ be the absolute Galois
group ${\rm Gal}(\overline K/K)$.
By the assumption
that the log ramification
is bounded by $R+$,
we have a unique element
$e$ in
$\Gamma(x,(\delta^{(R)*}j^{(R)}_*{\cal H})|_x)=
\Gamma(x,(h^*j^{(R)}_*{\cal H})|_x)$
whose image in
$\Gamma(x,(j_*{\cal E}nd({\cal F}))|_x)$
is the identity
of ${\cal F}_{\bar \eta}$.
The image of $e$ in
$\Gamma(x,(j_{C*}
{\cal H}om(f^*{\cal F}_{C'}
,{\cal F}_C))|_x)$
defines a $G_K$-homomorphism
$\varphi:{\cal F}_{f(\bar \eta)}
\to {\cal F}_{\bar \eta}$.
Switching the two factors,
we obtain a $G_K$-homomorphism
$\psi:{\cal F}_{\bar \eta}
\to {\cal F}_{f(\bar \eta)}$.
Since the construction
is compatible
with the composition,
the maps $\varphi$ and $\psi$
are the inverse of each other.
\qed

We study the higher direct image
$R^qj^{(R)}_*{\cal H}$.
We put
$I^+=\{i|1\le i\le m,r_i>0\}$
and 
$D^+=\bigcup_{i\in I^+}D_i$.
First, we consider the case
where the coefficients
of $R=\sum_ir_iD_i$ are integers.
If $R$ is integral,
the inverse image
$E^+=(X\times X)^{(R)}\times
_XD^+$
is identified with
a vector bundle
${\mathbf V}
(\Omega^1_X(\log D)(R))\times_X
D^+$ over $D^+$
by Lemma \ref{lmR}.3.
We prepare a global analogue
of Lemma \ref{lmmu}.

\begin{lm}\label{lmD}
Assume $R$ is integral.

1. 
We identify
the fiber product
$(X\times X)
\times_X
(X\times X)$
with respect to
the second and the
first projections
$X\times X\to X$
with
$X\times 
X\times X$.
Then, there exists 
a smooth map
$\mu:(X\times X)^{(R)}
\times_X
(X\times X)^{(R)}
\to (X\times X)^{(R)}$
that makes the diagram
$$\begin{CD}
(X\times X)^{(R)}
\times_X
(X\times X)^{(R)}
@>{\mu}>> (X\times X)^{(R)}\\
@VVV@VVV\\
(X\times X)
\times_X
(X\times X)=
X\times 
X\times X
@>{{\rm pr}_{13}}>>
X\times X
\end{CD}$$
commutative.

2. Let $D^+$ be the support
of $R$
and we identify
$E^+=(X\times X)^{(R)}
\times_XD^+$
with the vector bundle
${\mathbf V}
(\Omega^1_X(\log D)\otimes
{\cal O}(R))\times_XD^+$
as above.
The restriction of
$\mu$ defines
the addition
$E^+\times_{D^+}E^+\to E^+$
of the vector bundle 
$E^+={\mathbf V}(\Omega^1_X(\log D)(R))
\times_XD^+$.
\end{lm}
{\it Proof.}
1.
Let
$P=(X\times X)^\sim
\times_X
(X\times X)^\sim$
be the fiber product
with respect to
the second and the
first projections
$(X\times X)^\sim\to X$.
We define $P^{(R)}\to P$
by applying the constuction
in \S2.2 to
the smooth map
$P=(X\times X)^\sim
\times_X
(X\times X)^\sim
\to X$
and the diagonal section
$X\to P$.
The projections
$P\to (X\times X)^\sim$
induce an isomorphism
$P^{(R)}\to
(X\times X)^{(R)}
\times_X
(X\times X)^{(R)}$
by Corollary \ref{corPRf}.

On $P=(X\times X)^\sim
\times_X
(X\times X)^\sim$,
the pull-backs of
${\rm pr}_1^*D_i$
and
${\rm pr}_3^*D_i$
are equal
for each component $D_i$
of $D$.
Hence
the map
${\rm pr}_{13}
:(X\times X)
\times_X
(X\times X)
\to
X\times X$
is lifted to
$P=(X\times X)^\sim
\times_X
(X\times X)^\sim\to
(X\times X)^\sim$.
This is uniquely lifted
to a smooth map $P^{(R)}
\to (X\times X)^{(R)}$
by Lemma \ref{lmPRf}.

2.
The restriction
$E^+\times_{D^+}E^+\to E^+$
is a linear map
of vector bundles
by Lemma \ref{lmPRf}.2.
Hence, it suffices
to show that
the compositions
with the injections
$i_1,i_2:E^+\to E^+\times_{D^+}E^+$
of the two factors
are the identity of $E^+$.
We consider the map
$\iota_1:(X\times X)^{(R)}
\to
(X\times X)^{(R)}
\times_X
(X\times X)^{(R)}$
defined by
the identity
of $(X\times X)^{(R)}$
and $\delta^{(R)}\circ {\rm pr}_2$.
Then, its restriction
$E^+\to E^+\times_{D^+}E^+$
is the injection
of the first factor.
Since
the composition
$\mu\circ \iota_1$
is the identity,
the composition
$\mu\circ i_1:
E^+\to E^+\times_{D^+}E^+
\to E^+$ 
is the identity.
Similarly,
by considering the map
$\iota_2:(X\times X)^{(R)}
\to
(X\times X)^{(R)}
\times_X
(X\times X)^{(R)}$
defined by
$\delta^{(R)}\circ {\rm pr}_1$
and the identity
of $(X\times X)^{(R)}$,
we see that
$\mu\circ i_2$
is the identity.
Hence the assertion follows.
\qed

\begin{pr}\label{prD}
Let $X$ be a smooth scheme
over $k$ and
${\cal F}$ be a smooth sheaf
on the complement
$U=X\setminus D$
of a divisor
with simple normal crossings.
Let $R=\sum_ir_iD_i\ge 0$
be an effective divisor
with integral coefficients
$r_i\ge 0$.
Assume that the log
ramification of ${\cal F}$
is bounded by $R+$.
We put
$D^+=\bigcup_{i:r_i>0}D_i$
and 
$E^+=
{\mathbf V}
(\Omega^1_X(\log D)\otimes
{\cal O}(R))\times_X
D^+$.
Let $j^{(R)}:U\times U\to
(X\times X)^{(R)}$
be the open immersion.

1. For every integer $q\ge 0$,
the restriction of
$R^qj^{(R)}_*{\cal H}$
on $E^+$
is additive.

2. 
Let $S^q\subset E^{+\vee}=
{\mathbf V}
(\Omega^1_X(\log D)^\vee \otimes
{\cal O}(-R))\times_X
D^+$
be the dual support of
$R^qj^{(R)}_*{\cal H}|_{E^+}$.
Then, we have
$S^q\subset S^0$.

3.
Let $D_i$ be an irreducible
component $D_i$ of $D^+$
and $\xi_i$ be the generic point.
Then,
the intersection
$S^0\cap E_i^\vee$
is a subset
of the closure
$\overline{S^0_{\xi_i}}$
of the generic fiber.
\end{pr}

{\it Proof.}
Since $\mu:
(X\times X)^{(R)}\times_X
(X\times X)^{(R)}\to
(X\times X)^{(R)}$
is smooth,
the base change map
$\mu^*Rj^{(R)}_*{\cal H}
\to Rj_{3*}{\rm pr}_{13}^*{\cal H}$
is an isomorphism,
where $j_3:U\times U\times U
\to 
(X\times X)^{(R)}
\times_X
(X\times X)^{(R)}$
denotes the open immersion.
Hence, the composition
${\cal H}\boxtimes{\cal H}
={\cal H}om(
{\rm pr}_2^*{\cal F},
{\rm pr}_1^*{\cal F})
\otimes
{\cal H}om(
{\rm pr}_3^*{\cal F},
{\rm pr}_2^*{\cal F})
\to
{\cal H}om(
{\rm pr}_3^*{\cal F},
{\rm pr}_1^*{\cal F})
={\rm pr}_{13}^*{\cal H}$
induces
\begin{equation}
\begin{CD}
Rj^{(R)}_*{\cal H}\boxtimes Rj^{(R)}_*{\cal H}
@>>>
Rj_*{\rm pr}_{13}^*{\cal H}=
\mu^*Rj^{(R)}_*{\cal H}.
\end{CD}\label{eqpa}
\end{equation}

Let $\bar x$ be
an arbitrary geometric point
of $D^+$.
We show that
the restriction of
$R^qj^{(R)}_*{\cal H}$
on the fiber $E^+_{\bar x}$
satisfies
the condition
(2) in Proposition \ref{pradd}.
By the assumption that the log
ramification of ${\cal F}$
is bounded by $R+$,
we have a unique section
$e\in \Gamma(X,\delta^{(R)*}
j^{(R)}_*{\cal H})$
lifting the identity
${\rm id}_{\cal F}
\in \Gamma(X,j_*\delta_U^*{\cal H})$.
Take an \'etale neighborhood 
$V\to (X\times X)^{(R)}$
of $\bar x$
and a section
$\tilde e\in \Gamma(V,
j^{(R)}_*{\cal H})$
whose stalk in $
(j^{(R)}_*{\cal H})_{\bar x}=
(\delta^{(R)*}j^{(R)}_*{\cal H})_{\bar x}$
is the stalk of $e$ above.

Since $e$ is
a lifting of the identity,
the pairing (\ref{eqpa}) with the
restriction of $\tilde e$
is an isomorphism
${\rm pr}_1^*{\cal H}
\to {\rm pr}_{13}^*{\cal H}$
on $(U\times U)\times_X 
((U\times U)\times_{(X\times X)^{(R)}}V)$.
Hence, 
the pairing (\ref{eqpa}) with $\tilde e$
defines an isomorphism
${\rm pr}_1^*Rj^{(R)}_*{\cal H}
\to \mu^*Rj^{(R)}_*{\cal H}$
on $(X\times X)^{(R)}\times_X V$.
Thus, for a closed
point $a\in E_{\bar x}$
in the image of
$V\times_X\bar x
\to (X\times X)^{(R)}\times_X\bar x=E_{\bar x}$,
the restriction
$Rj^{(R)}_*{\cal H}|_{E_{\bar x}}$
is isomorphic
to the translate
$(+a)^*
(Rj^{(R)}_*{\cal H}|_{E_{\bar x}})$.
Since $E_{\bar x}$
is generated by the image of
$V\times_X\bar x$,
there is an isomorphism
$Rj^{(R)}_*{\cal H}|_{E_{\bar x}}
\to (+a)^*
(Rj^{(R)}_*{\cal H}|_{E_{\bar x}})$
for every closed
point $a\in E_{\bar x}$.
Thus the sheaf
$R^qj^{(R)}_*{\cal H}|_{E_{\bar x}}$
satisfies the condition (2)
in Proposition \ref{pradd}
and hence 
is additive for every $q\ge 0$.

2. It suffices to apply
Proposition \ref{pradd2} to
$j^{(R)}_*{\cal H}\boxtimes R^qj^{(R)}_*{\cal H}
\to
\mu^*R^qj^{(R)}_*{\cal H}$.

3. 
It follows immediately from
Lemma \ref{lmadd0}.
\qed

We consider the general case.
For a non-empty
subset $I\subset I^+$,
we put $D_I=
\bigcap_{i\in I}D_i$
and $D_I^\circ
=D_I\setminus 
\bigcup_{i\in I^+\setminus I}
(D_i\cap D_I)$.
Recall that
$n_i$ denotes
the denominator
of $r_i=m_i/n_i$
and $n_I$
is the least common multiple
of $n_i$ for $i\in I$.
The inverse image
$E_I^{\circ}
=((X\times X)^{(R)}\times_XD_I^\circ)_{\rm red}$
is identified with
${\mathbf V}_{n_I}
(\Omega^1_X(\log D),{\cal O}(n_IR))\times_X
D_I^\circ$
by Lemma \ref{lmR}.1.

\begin{pr}\label{prR}
Let the notation be
as in Proposition {\rm \ref{prD}}
except that we
do not assume the coefficients
of $R$ are integers.
Assume that the log
ramification of ${\cal F}$
is bounded by $R+$.
Let $I\subset I^+
=\{i|1\le i\le m,r_i>0\}$
be a non-empty subset
and 
$E_I^{\circ}
={\mathbf V}_{n_I}
(\Omega^1_X(\log D),{\cal O}(n_IR))\times_X
D_I^\circ$
be the reduced inverse image.

1. For every integer $q\ge 0$,
the restriction of
$R^qj^{(R)}_*{\cal H}$
on $E_I^{\circ}$
is potentially additive.

2. Let $S^q_I\subset E_I^{\circ\vee}
={\mathbf V}_{n_I}
(\Omega^1_X(\log D)^\vee,{\cal O}(-n_IR))\times_X
D_I^\circ$
be the dual support of
$(R^qj^{(R)}_*{\cal H})|_{E_I^{\circ}}$.
Then, we have
$S^q_I\subset S^0_I$.

3. For $i\in I$,
let $\xi_i$ be the 
generic point of
the irreducible component $D_i$
and $F_i=\kappa(\xi_i)$
be the function field
of $D_i$.
We consider the canonical map
$E_i^\vee\times_{D_i}D_I^\circ
={\mathbf V}_{n_i}
(\Omega^1_X(\log D)^\vee,{\cal O}(-n_iR))\times_X
D_I^\circ\to
E_I^{\circ\vee}
={\mathbf V}_{n_I}
(\Omega^1_X(\log D)^\vee,{\cal O}(-n_IR))\times_X
D_I^\circ$.
Then, 
$S^0_I$ is a subset
of the image of the
intersection
$\overline{S^0_{i,\xi_i}}\times_{D_i}D_I^\circ
\subset
E_i^\vee\times_{D_i}D_I^\circ$
of the closure of the generic fiber.
\end{pr}
{\it Proof.}
1. By replacing
$X$ by $X\setminus
\bigcup_{i\in I^+\setminus I}D_i$,
we may assume $I^+=I$,
$E_I=E_I^\circ$
and $n=n_I$.
Since the assertion is Zariski local,
we may take a smooth map
$X\to A={\mathbf A}^m
={\rm Spec}\ k[T_1,\ldots,T_m]$
such that
$D$ is the inverse image
of the union of
coordinate hyperplanes.
We put
$A'={\mathbf A}^m\times {\mathbf G}_m^m
={\rm Spec}\ k[T_1,\ldots,T_m,
U_1^{\pm1},\ldots,U_m^{\pm1}]$
and 
$\widetilde A
={\rm Spec}\ k[S_1,\ldots,S_m,
U_1^{\pm1},\ldots,U_m^{\pm1}]$
and define a map
$\widetilde A
\to A'$ by
$T_i\mapsto U_iS_i^{n_i}$.
We consider the base change
$X\gets X'\gets \widetilde X$ of
$A\gets A'\gets \widetilde A$.

We define schemes
$(X\times X')^{\sim}$ and 
$(X\times \widetilde X)^{\sim}$
by the cartesian diagram
$$\begin{CD}
(X\times X)^\sim
@<<<
(X\times X')^\sim
@<<<
(X\times \widetilde X)^\sim
\\
@V{p_2}VV @VVV @VVV \\
X@<<< X'@<<<\widetilde X
\end{CD}$$
and consider the sections
$X'\to (X\times X')^{\sim}$
and
$\widetilde X
\to (X\times \widetilde X)^{\sim}$
induced by the log diagonal
$X\to (X\times X)^{\sim}$.
The map
$\widetilde X\to X$
induces 
$(\widetilde X\times \widetilde X)^\sim
\to 
(X\times \widetilde X)^\sim$.
By applying the construction
in \S2.2 to the pull-backs
$R'$ and $\widetilde R$ of $R$
to $X'$ and
to $\widetilde X$,
we obtain the commutative diagram
\begin{equation}
\begin{CD}
(X\times X)^{(R)}
@<f<<
(X\times X')^{(R')}
@<g<<
(X\times \widetilde X)^{(\widetilde R)}
@<h<<
(\widetilde X\times \widetilde X)^{(\widetilde R)}
\\
@VVV @VVV @VV{\tilde p}V @.\\
X@<<< X'@<<<\widetilde X.
\end{CD}\label{eqls}
\end{equation}

Since $X'\to X$
is smooth,
the left square of
(\ref{eqls}) is cartesian
and 
the horizontal arrow $f
:(X\times X')^{(R')}
\to (X\times X)^{(R)}$ 
is smooth.
Since $\widetilde R$
has integral coefficients,
the right vertical arrow
$\tilde p:
(X\times \widetilde X)^{(\widetilde R)}
\to \widetilde X$ is smooth.
We show that the map
$h:(\widetilde X\times 
\widetilde X)^{(\widetilde R)}
\to
(X\times \widetilde X)^{(\widetilde R)}$
is smooth.
By Lemma \ref{lmPRf}.3,
it suffices to show that the map
$(\widetilde X\times 
\widetilde X)^\sim
\to
(X\times \widetilde X)^\sim$
is smooth.
Thus, it is reduced to showing that
the map
$(\widetilde A\times 
\widetilde A)^\sim
\to
(A\times \widetilde A)^\sim$
is smooth.
Since the map
$(\widetilde A\times 
\widetilde A)^\sim
={\rm Spec}\ k
[S_i,U_i^{\pm1},S'_i,U_i^{\prime\pm1},
V_i^{\pm1}\ (i=1,\ldots,m)]/
(S'_i-V_iS_i\ (i=1,\ldots,m))
={\rm Spec}\ k
[S_i,U_i^{\pm1},U_i^{\prime\pm1},
V_i^{\pm1}\ (i=1,\ldots,m)]
\to
(A\times \widetilde A)^\sim
={\rm Spec}\ k
[T_i,S'_i,U_i^{\prime\pm1},W_i^{\pm1}\ (i=1,\ldots,m)]/
(U'_iS_i^{\prime n_i}-W_iT_i\ (i=1,\ldots,m))
={\rm Spec}\ k
[S'_i,U_i^{\prime\pm1},W_i^{\pm1}\ (i=1,\ldots,m)]$
is defined by
$W_i\mapsto V_i^{n_i}U'_i/U_i$,
it is smooth.

We put $U'=U\times_XX'$ and $
\widetilde U=U\times_X
\widetilde X$
and consider the diagram
$$\begin{CD}
(X\times X)^{(R)}
@<f<<
(X\times X')^{(R')}
@<g<<
(X\times \widetilde X)^{(\widetilde R)}
@<h<<
(\widetilde X\times \widetilde X)^{(\widetilde R)}
\\
@A{j^{(R)}}AA @A{j^{(R)\prime}}AA @A{j^{(R)\sim}}AA @A{\tilde j^{(R)}}AA\\
U\times U@<<< U\times U'
@<<<U\times \widetilde U
@<<<\widetilde U\times \widetilde U.
\end{CD}$$
where the vertical arrows are open immersions.
We consider the pull-backs
${\cal H}'$,
${\cal H}^\sim$ 
and $\widetilde {\cal H}$
of ${\cal H}$ respectively
on 
$U\times U'$,
$U\times \widetilde U$
and on 
$\widetilde U\times \widetilde U$.

Since 
$\widetilde R$ is integral,
the restriction
of $R^q\tilde j^{(R)}_*\widetilde{\cal H}$
on $\widetilde E^+$
is additive
by Proposition \ref{prD}
for every $q\ge0$.
Since $h$ is smooth,
the base change map
$h^*R^qj^{(R)\sim}_*{\cal H}^\sim
\to 
R^q\tilde j^{(R)}_*\widetilde {\cal H}$
is an isomorphism.
The conormal sheaves
${\cal N}_{\widetilde X/
(X\times \widetilde X)^\sim}$
and
${\cal N}_{\widetilde X/
(\widetilde X\times \widetilde X)^\sim}$
are canonically identified
with $\Omega^1_X(\log D)
\otimes {\cal O}_{\widetilde X}$
and with $\Omega^1_{\widetilde X}
(\log \widetilde D)$
respectively.
Since the map
$\Omega^1_X(\log D)
\otimes {\cal O}_{\widetilde X}
\to \Omega^1_{\widetilde X}
(\log \widetilde D)$
is a locally splitting injection,
the restriction
of $R^qj^{(R)\sim}_*{\cal H}^\sim$
is additive
by Lemma \ref{lmadd2}.

To study the restriction of
$R^qj^{(R)\prime}_*{\cal H}'$,
we introduce some notations.
For $i\in I$,
let $D'_i\subset X'$
be the inverse image of $D_i$
and $\widetilde D_i$
be the divisor defined by $S_i$.
We put $D'_I=\bigcap_{i\in I}D'_i$
and
$\widetilde D_I=\bigcap_{i\in I}
\widetilde D_i$.
The natural map
$\widetilde D_I
\to D'_I$
is an isomorphism.
Let $E'_I\subset 
(X\times X')^{(R')}$ and
$\widetilde E_I\subset
(X\times \widetilde X)^{(\widetilde R)}$
be the reduced inverse images of
$D'_I$ and of 
$\widetilde D_I$.
Recall $n=n_I$
is the least common multiple
of the denominators $n_i$ for $i\in I$.
By Lemma \ref{lmR}.1, 
we have a canonical isomorphism
$E'_I\to
{\mathbf V}_n(\Omega^1_X(\log D),
{\cal O}(-nR'))_{D'_I}$
and
$\widetilde E_I\to
{\mathbf V}_1(\Omega^1_X(\log D),
{\cal O}(-\widetilde R))_{\widetilde D_I}$.
The natural map
$\pi_I:\widetilde E_I\to
E'_I$
is then identified with the
canonical map
$\pi_n:
{\mathbf V}_1(\Omega^1_X(\log D),
{\cal O}(-\widetilde R))_{D'_I}
\to
{\mathbf V}_n(\Omega^1_X(\log D),
{\cal O}(-nR'))_{D'_I}$.

Let $n'_i$ be the prime-to-$p$
part of $n_i$ and
$n'$ be the prime-to-$p$ part
of $n=n_I$.
We consider the natural action
of $G=\prod_{i\in I}\mu_{n'_i}$
on $\widetilde X$ over $X'$.
Since the map
$\widetilde D_I
\to D'_I$
is an isomorphism,
the action of $G$ on
$\widetilde D_I$
is trivial.
The action of $G$
on $\widetilde E_I$
factors through
the product map
$G\to \mu_{n'}$
and the action of $\mu_{n'}$
on $\widetilde E_I$
is by the multiplication.

We show that the restriction 
of $R^qj^{(R)\prime}_*{\cal H}'$
on $E'_I$
is potentially additive.
The canonical map
$R^qj^{(R)\prime}_*{\cal H}'
\to g_*(R^qj^{(R)\sim}_*{\cal H}^\sim)^G$
is an isomorphism.
Let
$G_1$ be the kernel
of 
$G\to \mu_{m'}$.
Then since
the restriction of
$R^qj^{(R)\sim}_*{\cal H}^\sim$
on $\widetilde E_I$
is additive,
its $G_1$-fixed part
$(R^qj^{(R)\sim}_*{\cal H}^\sim)|_{
\widetilde E_I}^{G_1}$
is also additive.
Hence by Lemma \ref{lmpadd},
the $\mu_{n'}$-fixed part
$\pi_{n*}\left((R^qj^{(R)^\sim}_*{\cal H}^\sim)|_{
\widetilde E_I}^{G_1}\right)^{\mu_{n'}}=
g_*(R^qj^{(R)\sim}_*{\cal H}^\sim)^G|_{E'_I}$
is potentially additive.
Thus the restriction
$(R^qj^{(R)\prime}_*{\cal H}')|_{E'_I}$
is potentially additive.

Since the map $X'\to X$ is smooth,
the base change map
$f^*R^qj^{(R)}_*{\cal H}
\to 
R^qj^{(R)\prime}_*{\cal H}'$
is an isomorphism.
Since the map $X'\to X$ 
admits a section,
the restriction of
$R^qj^{(R)}_*{\cal H}$
on $E_I$ is 
also potentially additive.

2.
Similarly as in the proof of 1,
we may assume $I=I^+$
and $D_I=D_I^\circ$.
We show the inclusion
$S^q_I\subset S^0_I$.
Let $S^{q\sim}_I\subset E_I^{\sim\vee}$
be the dual support
of the additive sheaf
$(R^qj^{(R)^\sim}_*{\cal H}^\sim)^{G_1}$.
We apply Proposition \ref{pradd2}
to the map
$(j^{(R)^\sim}_*{\cal H}^\sim)^{G_1}
\boxtimes
(R^qj^{(R)^\sim}_*{\cal H}^\sim)^{G_1}
\to 
\mu^*
(R^qj^{(R)^\sim}_*{\cal H}^\sim)^{G_1}$
and the pull-back to
$\Gamma(\widetilde X,
\tilde \delta^{(R)*}j^{(R)^\sim}_*{\cal H}^\sim)$
of the section
$e\in \Gamma(X,
\delta^{(R)*}j^{(R)}_*{\cal H})$
lifting the identity of ${\cal F}$.
Then, we obtain the inclusion
$S^{q\sim}_I\subset
S^{0\sim}_I$.

Since the dual support
$S^{q\prime}_I\subset E_I^{\prime\vee}$
of the potentially additive sheaf
$R^qj^{(R)\prime}_*{\cal H}'$
is the image of
$S^{q\sim}_I$
by the canonical map
$E_I^{\sim*}
={\mathbf V}_1(\Omega^1_X(\log D)^\vee,
{\cal O}_{\widetilde X}(-\widetilde R))
\times_{\widetilde X}\widetilde D_I
\to 
E_I^{\prime\vee}=
{\mathbf V}_{n_I}(\Omega^1_X(\log D)^\vee\otimes 
{\cal O}_{X'}, {\cal O}(-n_IR'))
\times_{X'}\widetilde D_I$,
we obtain the inclusion
$S^{q\prime}_I\subset
S^{0\prime}_I$.
Thus, we deduce
$S^q_I\subset
S^0_I$ by pull-back.

3. 
We have the inclusion
$S^{0\sim}_I\subset
\overline{S^{0\sim}_{i,\xi_i}}
\times_{\widetilde D_i}\widetilde D_I^\circ$
by Lemma \ref{lmadd0}.
Hence the assertion follows
as in the proof of 2.
\qed

For an integer $n>0$
such that $nR$
has integral coefficients,
we define the dual support
$$S_{\cal F}^{(n\cdot R)}
\subset
E_n^\vee={\mathbf V}_n(
\Omega^1_X(\log D)^\vee,
{\cal O}(-nR))_{D^+}$$
as a constructible subset
as follows.
Let $I$
be a non-empty subset
of $I^+=
\{i|1\le i\le m,r_i>0\}$.
Then,
the restriction
${\cal H}_I^\circ=
j^{(R)}_*{\cal H}|_{E_I^\circ}$
on $E_I^\circ
={\mathbf V}_{n_I}(
\Omega^1_X(\log D),
{\cal O}(n_IR))\times_XD_I^{\circ}$
is potentially additive
by Proposition
\ref{prD}.
Hence the dual support
$S_{{\cal H}_I^\circ}$
is defined
as a constructible subset
of the dual
${\mathbf V}_{n_I}(
\Omega^1_X(\log D)^\vee,
{\cal O}(-n_IR))\times_XD_I^{\circ}$.
Since $n$ is
divisible by $n_I$,
the canonical map
$\pi_{nn_I}:
{\mathbf V}_{n_I}(
\Omega^1_X(\log D)^\vee,
{\cal O}(-n_IR))\times_XD_I^{\circ}
\to 
{\mathbf V}_n(
\Omega^1_X(\log D)^\vee,
{\cal O}(-nR))\times_XD_I^{\circ}
=E_n^\vee\times_XD_I^{\circ}$
is defined.

\begin{df}\label{dfdlsp}
Let the notation be as above.
We define the dual support
$S_{\cal F}^{(n\cdot R)}
\subset E_n^\vee$
with respect to $R$
as the union
$$S_{\cal F}^{(n\cdot R)}
=\bigcup_{I\subset I^+,
I\neq \emptyset}
\pi_{nn_I}(S_{{\cal H}_I^\circ}).$$
We say that the log ramification of ${\cal F}$
is non-degenerate with respect to $R$
if the intersection
of dual support
$S_{\cal F}^{(n\cdot R)}$
with the $0$-section of
$E_n^\vee$ is empty.
\end{df}

For $n|m$,
we have
$S_{\cal F}^{(m\cdot R)}
=\pi_{mn}(S_{\cal F}^{(n\cdot R)})$.

\begin{cor}\label{corSF}
Assume that the log ramification of
${\cal F}$ is bounded by $R+$.

1. For an irreducible component
$D_i$ of $D$,
let $\xi_i$ be the generic point
of $D_i$.
Then, we have
$$S^{(n\cdot R)}_{\cal F}
\subset \bigcup_{i\in I^+}
\overline{S^{(n\cdot R)}_{{\cal F},\xi_i}}
.$$

2.
Assume that the log ramification of ${\cal F}$
is non-degenerate with respect to $R$
and that $\Lambda$
is a finite extension of ${\mathbb Q}_\ell$.
Then
$Rp_*(R^qj_*{\cal H}|_{D^+})$
and 
$Rp_!(R^qj_*{\cal H}|_{D^+})$
are $0$
for every $q\ge 0$.
\end{cor}
{\it Proof.}
1.
Clear from Proposition \ref{prR}.3.

2.
Clear from Lemma \ref{lmnd}.
\qed

We make explicit
the relation between
the dual support
and the refined Swan character.
Let $D_i$ be an irreducible component
of $D$, $\xi_i$ be the
generic point of $D_i$
and $K_i$ be the
fraction field of
the henselization
${\cal O}_{X,\xi_i}^h$ of the local ring.
The residue field $F_i$
of $K_i$ is the function field of $D_i$.
Let $\chi$ be a character of
${\rm Gr}^{r_i}_{\log}G_{K_i}$.
Recall that ${\mathfrak m}_{\overline K}^r
\supset {\mathfrak m}_{\overline K}^{r+}$
denote
$\{a\in \overline K|v(a)\ge r\}
\supset \{a\in \overline K|v(a)> r\}$.
Then, by Corollary \ref{corgr},
the refined Swan character
of $\chi$ defines 
an $\overline{F_i}$-valued
point
$${\rm rsw}\ \chi\in 
\Omega^1_{F_i}(\log)
\otimes {\mathfrak m}_{\overline K_i}^{(-r_i)}/
{\mathfrak m}_{\overline K_i}^{(-r_i)+}
=
{\mathbf V}_1(
\Omega_X(\log D)^\vee
\otimes \overline{F_i},
{\mathfrak m}_{\overline K_i}^{r_i}/
{\mathfrak m}_{\overline K_i}^{r_i+})
(\overline{F_i})
.$$

\begin{lm}\label{lmSH}
Assume the log ramification
of ${\cal F}$ is bounded by $R+$.
Let $D_i$ be
an irreducible component of $D$
such that $r_i>0$.
We consider the stalk
${\cal F}_{\bar \eta_i}$
as a representation
of $G_{K_i}$
and 
the direct sum decomposition
${\cal F}_{\bar \eta_i}=
\bigoplus_\chi 
\chi^{\oplus n_\chi}$
of the restriction
to $G_{K_i,\log}^{r_i}$
by characters of
${\rm Gr}^{r_i}_{\log}G_{K_i}$.
Let 
$\pi_n:{\mathbf V}_1(
\Omega_X(\log D)^\vee
\otimes \overline{F_i},
{\mathfrak m}_{\overline K_i}^{r_i}/
{\mathfrak m}_{\overline K_i}^{r_i+})
\to
{\mathbf V}_n(
\Omega_X(\log D)^\vee,
{\cal O}(-nR))_{\xi_i}
=E_{n,\xi_i}^\vee$
be the canonical map.

Then, the generic fiber 
$S^{(n\cdot R)}_{{\cal F},\xi_i}
\subset E_{n,\xi_i}^\vee=
{\mathbf V}_n(
\Omega_X(\log D)^\vee,
{\cal O}(-nR))_{\xi_i}$
of the dual support
consists of the images 
$\pi_n({\rm rsw}\ \chi)$ of
the refined Swan characters
of $\chi$ appearing
in the direct sum decomposition
${\cal F}_{\bar \eta_i}=
\bigoplus_\chi 
\chi^{\oplus n_\chi}$.
\end{lm}
{\it Proof.}
It is reduced to the case
where $r_i$ is an integer,
by Lemma \ref{lmr}
and by the proof of Proposition \ref{prR}.1.
In the case where
$r_i$ is an integer,
it follows from
Corollary \ref{corgr1}.
\qed

\begin{cor}\label{corSH}
Assume the log ramification
of ${\cal F}$ is bounded by $R+$.
The equality
${\cal F}_{\bar \eta_i}=
{\cal F}^{(r_i)}_{\bar \eta_i}$
is equivalent to
that the generic fiber
$S^{(n\cdot R)}_{{\cal F},\xi_i}$
does not contain $0$.
Also 
the vanishing
${\cal F}^{(r_i)}_{\bar \eta_i}=0$
is equivalent to
that the generic fiber
$S^{(n\cdot R)}_{{\cal F},\xi_i}$
is a subset of $\{0\}$.
\end{cor}
{\it Proof.}
By Corollary \ref{corgr},
the map 
${\rm rsw}:{\rm Hom}(
{\rm Gr}_{\log}^rG_K,
{\mathbb F}_p)
\to
\Omega^1_F(\log)\otimes_F
{\mathfrak m}^{(-r)}_{\overline K}/
{\mathfrak m}^{(-r)+}_{\overline K}$
is injective.
Hence it follows from
Lemma \ref{lmSH}.
\qed

We study the functoriality of
the dual support
$S_{\cal F}^{(n\cdot R)}$.
Let $Y$ be a smooth scheme over $k$
and $f:Y\to X$
be a morphism over $k$.
Assume that the reduced inverse image
$D_Y=(D\times_XY)_{\rm red}$
is a divisor with simple normal crossings
and let $R_Y$ be the pull-back
$f^*R$.
Let $n>0$ be
an integer such that
$nR$ is integral
and we put
$E_n^\vee={\mathbf V}_n(
\Omega^1_X(\log D)^\vee,
{\cal O}(-nR))$
and
$E_n^{\prime \vee}={\mathbf V}_n(
\Omega^1_Y(\log D_Y)^\vee,
{\cal O}(-nR_Y))$.
Then, the canonical map
$f^*\Omega^1_X(\log D)
\to \Omega^1_Y(\log D_Y)$
induces a map
$\varphi:E_n^\vee\times_XY
\to E_n^{\prime \vee}$.

\begin{lm}\label{lmSHf}
Let $Y$ be a smooth scheme over $k$
and $f:Y\to X$
be a morphism over $k$.
Assume that the reduced inverse image
$D_Y=(D\times_XY)_{\rm red}$
is a divisor with simple normal crossings
and let $R_Y$ be the pull-back
$f^*R$.
Let $n>0$ be
an integer such that
$nR$ is integral
and let
$\varphi:E_n^\vee\times_XY
\to E_n^{\prime \vee}$
be the map defined above.

Let ${\cal F}$
be a smooth sheaf on
$U=X\setminus D$
and ${\cal F}_V$
be the pull-back to
$V=U\times_XY=Y\setminus D_Y$.
Assume the log ramification of
${\cal F}$ is bounded by $R+$.
Then, we have
$$\varphi(f^*
S_{\cal F}^{(n\cdot R)})
\subset 
S_{{\cal F}_V}^{(n\cdot R_Y)}.$$
\end{lm}
{\it Proof.}
Let $D_1,\ldots,D_m$
and 
$D'_1,\ldots,D'_{m'}$
be the components
of $D$ and of $D_Y$
respectively.
We put $R=
\sum_{i=1}^mr_iD_i$ and
$R_Y=
\sum_{j=1}^{m'}r'_jD'_j$.
Let $J$ be a non-empty subset
of $J^+=\{j|r'_j>0,j=1,\ldots, m'\}$
and put
$I=\{i|f^{-1}(D_i)\supset D'_J,
r_i>0, i=1,\ldots,m\}$.
The map
$(f\times f)^{(R)}:
(Y\times Y)^{(R_Y)}
\to (X\times X)^{(R)}$
defined in the
proof of Lemma
\ref{lmcutY}
induces 
$E^{\prime \circ}_J
={\mathbf V}_{n_J}
(\Omega^1_Y(\log D_Y),
{\cal O}(n_JR_Y))
\to E_I^{\circ}
={\mathbf V}_{n_I}
(\Omega^1_X(\log D),
{\cal O}(n_IR))$.
Since the base change map
$(f\times f)^{(R)}j_*{\cal H}
\to j'_*(f_U\times f_U)^*{\cal H}$
is injective,
the assertion follows.
\qed

For an irreducible component
$D_i$ of $D$,
the residue map 
$\Omega^1_X(\log D)
\otimes_{{\cal O}_X}{\cal O}_{D_i}
\to {\cal O}_{D_i}$
defines a map
\begin{eqnarray*}
{\rm res}_i:
E_n^\vee\times_XD_i={\mathbf V}_n
(\Omega^1_X(\log D)^\vee,
{\cal O}(-nR))\times_XD_i
&\to&
{\mathbf V}_n
({\cal O}_{D_i},{\cal O}_{D_i}(-nR))\\
&&={\mathbf V}({\cal O}_{D_i}(-nR)).
\end{eqnarray*}

\begin{cor}\label{corSHf}
Let the notation be
as in Lemma {\rm \ref{lmSHf}}.
Let $D_i$ 
be an irreducible component
of $D$ and $D'_j$
be an irreducible component
of $D_Y$ such that
the multiplicity $e$
of $D'_j$ in the pull-back
$f^{-1}(D_i)$
is non-zero.
We consider the map
${\mathbf V}({\cal O}_{D'_j}(-nR))=
{\mathbf V}({\cal O}_{D_i}(-nR))\times_{D_i}D'_j
\to 
{\mathbf V}({\cal O}_{D_i}(-nR))$
induced by $f$.

1. 
We have
$e\cdot f^*{\rm res}_i(
S^{(n\cdot R)}_{\cal F})
\subset
{\rm res}_j(
S^{(n\cdot R_Y)}_{{\cal F}_V})$.

2. Assume 
${\rm res}_i(
S^{(n\cdot R)}_{\cal F})\setminus D_i
\to D_i$ is surjective
and $e$ is prime to $p$.
Then, the log ramification of
${\cal F}_V$ along $D'_j$
is not bounded by $r_j$.
\end{cor}

{\it Proof.}
1. By the commutative
diagram
$$\begin{CD}
E_n^\vee\times_XD'_j
@>{{\rm res}_i}>>
{\mathbf V}({\cal O}_{D'_j}(-nR))\\
@V{\varphi}VV @VVeV\\
E_n^{\prime\vee}\times_YD'_j
@>{{\rm res}_j}>>
{\mathbf V}({\cal O}_{D'_j}(-nR_Y)),
\end{CD}$$
it follows from Lemma \ref{lmSHf}.

2. By 1,
${\rm res}_j(
S^{(n\cdot R_Y)}_{{\cal F}_V})$
is not contained in
the 0-section.
\qed

\section{Characteristic cycle}

We recall in \S3.1 the definition of
the characteristic class
and compute it under
a certain assumption.
We propose a definition of
the characteristic class in some case
and prove
that it computes the characteristic class
in \S3.2. 

\subsection{Characteristic class}
We recall the definition of
the characteristic class.
For more detail on the construction,
we refer to \cite{cc} \S\S 1 and 2
and \cite{SGA5} \S\S 1-3.
Let $X$ be a scheme over $k$
and ${\cal F}$ be a constructible
sheaf of flat $\Lambda$-modules.
We put ${\cal K}_X=Ra^!\Lambda$ 
for the structural map
$a:X\to {\rm Spec}\ k$
and
$D_X{\cal F}
=R{\cal H}om({\cal F},{\cal K}_X)$.
We consider 
${\cal H}=
R{\cal H}om({\rm pr}_2^*{\cal F},
R{\rm pr}_1^!{\cal F})$
on $X\times X$.

The canonical pairing
$$
{\cal F}
\otimes^L
R\delta^!
{\cal H}
=\delta^*
{\rm pr}_2^*{\cal F}
\otimes^L
R\delta^!
R{\cal H}om({\rm pr}_2^*{\cal F},
R{\rm pr}_1^!{\cal F})
\longrightarrow 
R\delta^!R{\rm pr}_1^!{\cal F}={\cal F}
$$
induces
an isomorphism
\begin{equation}
\begin{CD}
H^0_X(X\times X,
{\cal H})
@>>>
{\rm End}_X({\cal F}).
\end{CD}
\label{eqcc1}
\end{equation}
Alternatively, one can apply the
canonical isomorphism \cite{SGA5} (3.2.1).
The inverse of the 
canonical isomorphism
${\cal F}\boxtimes
D_X{\cal F}
\to
R{\cal H}om({\rm pr}_2^*{\cal F},
R{\rm pr}_1^!{\cal F})={\cal H}$
and the canonical map
$R\delta^!\to \delta^*$
induce a map
\begin{equation}
\begin{CD}
H^0_X(X\times X,
{\cal H})
@>>>
 H^0(X,
{\cal F}\otimes^L
D_X{\cal F}).
\end{CD}
\label{eqcc2}
\end{equation}
The evaluation map
${\cal F}\otimes^L
D_X{\cal F}
\to {\cal K}_X$
induces a map
\begin{equation}
\begin{CD}
H^0(X,
{\cal F}\otimes^L
D_X{\cal F})
@>>>
H^0(X,{\cal K}_X).
\end{CD}
\label{eqcc3}
\end{equation}
We define the characteristic
class $C({\cal F})
\in H^0(X,{\cal K}_X)$
to be the image
of
$1\in {\rm End}_X({\cal F})$
by the composition
$$\begin{CD}
{\rm End}_X({\cal F})
@>{\rm (\ref{eqcc1})^{-1}}>>
H^0_X(X\times X,
{\cal H})
@>{\rm (\ref{eqcc2})}>>
 H^0(X,
{\cal F}\otimes^L
D_X{\cal F})
@>{\rm (\ref{eqcc3})}>>
H^0(X,{\cal K}_X).
\end{CD}
$$
If $X$ is proper,
we have an index formula
\begin{equation}
\chi(X_{\bar k},{\cal F})
=
{\rm Tr}\ C({\cal F})
\label{eqidx}
\end{equation}
for the Euler number
$\chi(X_{\bar k},{\cal F})
=\sum_{q=0}^{2\dim X}
(-1)^q\dim H^q(X_{\bar k},
{\cal F})$.

Assume that
$X$ is smooth of dimension $d$
and that ${\cal F}$
is a smooth sheaf
of free $\Lambda$-modules
of finite rank.
Then, 
the isomorphism
${\rm End}_X({\cal F})
\to
H^0_X(X\times X,{\cal H})$
(\ref{eqcc1})
is described as follows.
We put
${\cal H}_0=
{\cal H}om({\rm pr}_2^*{\cal F},
{\rm pr}_1^*{\cal F})$.
By the assumptions
on $X$ and on ${\cal F}$,
we have a canonical isomorphism
${\cal H}_0(d)[2d]\to {\cal H}=
R{\cal H}om({\rm pr}_2^*{\cal F},
R{\rm pr}_1^!{\cal F})$.
We identify
$\delta^*{\cal H}_0
=
{\cal E}nd({\cal F})$
and
$H^0(X,\delta^*{\cal H}_0)
=
{\rm End}_X({\cal F})$.
Then, the isomorphism
(\ref{eqcc1})
is the inverse of the cup-product
\begin{equation}
\begin{CD}
{\rm End}_X({\cal F})
=
H^0(X,\delta^*{\cal H}_0)
@>{\cup [X]}>>
H^0_X(X\times X,{\cal H})
\end{CD}
\label{eqid}
\end{equation}
with the cycle class
$[X]\in H^0_X(X\times X,\Lambda(d)[2d])$.
Further, in this case,
the evaluation map
$\delta^*{\cal H}\to {\cal K}_X$
is the tensor product of
the trace map
${\rm Tr}:\delta^*{\cal H}_0
=
{\cal E}nd({\cal F})
\to \Lambda$
with the isomorphism
$\Lambda(d)[2d]\to {\cal K}_X$
defined by the cycle class.
Thus, in this case,
we have
$$C({\cal F})=
{\rm rank}\ {\cal F}\cdot
(X,X)_{X\times X}$$
in $H^{2d}(X,\Lambda(d))$
where 
$(X,X)_{X\times X}=
(-1)^dc_d(\Omega^1_X)$
is the self intersection class.

We will compute 
the characteristic class
in some cases.
First we consider the tamely ramified case.
Let $X$ be a smooth scheme 
of dimension $d$ over $k$
and $U=X\setminus D$ be 
the complement
of a divisor $D$
with simple normal crossings.
We consider the
diagram
\begin{equation}
\xymatrix{
X\times X&
(X\times X)^\sim
\ar[l]_f&
U\times U
\ar[l]_{\ \ \tilde j}\\
&X\ar[lu]^{\delta}
\ar[u]_{\tilde \delta}
&U\ar[u]_{\delta_U}
\ar[l]^{j}
}\label{eqtame}
\end{equation}
where
$(X\times X)^\sim$
is the log product
and $f:(X\times X)^\sim
\to X\times X$
is the canonical map.
The diagonal maps
for $X$ and $U$ are
denoted by $\delta$
and $\delta_U$
respectively and
$\tilde \delta$
is the log diagonal map.
The map 
$\tilde j:U\times U
\to (X\times X)^\sim$
is the open immersion.

\begin{pr}\label{prtame}
Let the notation be as 
in the diagram {\rm (\ref{eqtame})}
above and let
${\cal F}$ be a smooth sheaf
of free $\Lambda$-modules
of finite rank on
$U=X\setminus D$.
Assume that
${\cal F}$ is tamely ramified along $D$.

We put ${\cal H}_0=
{\cal H}om({\rm pr}_2^*{\cal F},
{\rm pr}_1^*{\cal F})$
on $U\times U$
and $\overline {\cal H}=R{\cal H}om(
{\rm pr}_2^*j_!{\cal F},
{\rm pr}_1^!j_!{\cal F})$
on $X\times X$.
We also put
$\widetilde {\cal H}_0=
\tilde j_*{\cal H}_0$
and $\widetilde {\cal H}=
\widetilde {\cal H}_0(d)[2d]$
on $(X\times X)^\sim$.
Let $e\in \Gamma(X,
\tilde \delta^*\widetilde {\cal H}_0)$
be the unique element
that maps to the identity
${\rm id}_{\cal F}\in {\rm End}_U({\cal F})
=H^0(U,\delta_U^*{\cal H}_0)$
and let
$e\cup [X]
\in H^0_X((X\times X)^\sim,
\widetilde {\cal H})$
be the cup product
with the cycle class
$[X]\in H^0_X((X\times X)^\sim,$
$\Lambda(d)[2d])$.

1. There exists a unique map
$f^*\overline {\cal H}
\to \widetilde {\cal H}$
inducing the canonical isomorphism
${\cal H}=R{\cal H}om({\rm pr}_2^*{\cal F},
R{\rm pr}_1^!{\cal F})\to
{\cal H}om({\rm pr}_2^*{\cal F},
{\rm pr}_1^*{\cal F})(d)[2d]={\cal H}_0(d)[2d]$
on $U\times U$.

2. We consider the pull-back 
$f^*({\rm id}_{j_!\cal F})
\in H^0_{f^{-1}(X)}(
(X\times X)^\sim,\widetilde {\cal H})$
of the identity
${\rm id}_{j_!\cal F}\in 
{\rm End}_X(j_!{\cal F})=
H^0_X(X\times X,\overline{\cal H})$.
Then, we have
\begin{equation}
f^*({\rm id}_{j_!\cal F})
=e\cup [X]
\label{eqeX}
\end{equation}
in $H^0_{f^{-1}(X)}(
(X\times X)^\sim,\widetilde {\cal H})$.
Consequently,
we have
$$C(j_!{\cal F})=
{\rm rank}\ {\cal F}\cdot
(X,X)_{(X\times X)^\sim}$$
in $H^0(X,{\cal K}_X)$
where 
$(X,X)_{(X\times X)^\sim}=
(-1)^dc_d(\Omega^1_X(\log D))$
is the self intersection class.
\end{pr}
{\it Proof.}
1.
Let $\bar f:(X\times X)'\to X\times X$
be the log blow-up.
Let $(U\times X)'\subset 
(X\times X)'$
be the complement
of the proper transform
of $D\times X$
and we consider
the open immersions
$$\begin{CD}
(X\times X)^\sim
@>{j_2}>>
(U\times X)'
@>{j_1}>>
(X\times X)'.
\end{CD}$$
We put 
${\cal H}'=j_{1!}Rj_{2*}
\widetilde {\cal H}$
on $(X\times X)'$.
The log blow-up
$\bar f:(X\times X)'\to X\times X$
is an isomorphism on
the complement 
$U\times X$ of
$D\times X$.
The restriction 
of $\overline{\cal H}$
on $D\times X$ is 0
and the restriction
$\overline{\cal H}|_{
U\times X}
=R{\cal H}om({\rm pr}_2^*j_!{\cal F},
R{\rm pr}_1^!{\cal F})$
is canonically identified with
the restriction
${\cal H}'|_{
U\times X}
=R(1\times j)_*
{\cal H}_0(d)[2d]$.
Hence, 
there exists a unique map
$\bar f^*\overline {\cal H}
\to {\cal H}'$
inducing the canonical isomorphism
${\cal H}\to
{\cal H}_0(d)[2d]$
on $U\times U$.

2.
To prove the equality
$f^*({\rm id}_{j_!\cal F})
=e\cup [X]$
(\ref{eqeX}),
it suffices to show
the equality
$\bar f^*({\rm id}_{j_!\cal F})
=e\cup [X]$
in $H^0_{\bar f^{-1}(X)}(
(X\times X)',{\cal H}')$.
By \cite{cc} Lemma 2.2.4,
the adjunction
$\overline {\cal H}\to Rf_*
{\cal H}'$
of the canonical map
$f^*\overline {\cal H}\to 
{\cal H}'$
is an isomorphism.
Hence,
the pull-back
$\bar f^*:
H^0_X(X\times X,\overline {\cal H})
\to
H^0_{\bar f^{-1}(X)}
((X\times X)',
{\cal H}')$
is an isomorphism.
Since the restriction map
$H^0_X(X\times X,\overline {\cal H})
={\rm End}_X(j_!{\cal F})
\to
H^0_U
(U\times U,{\cal H})
={\rm End}_U({\cal F})$
is an isomorphism,
the arrows in the commutative diagram
$$\xymatrix{
H^0_X(X\times X,\overline {\cal H})
\ar[r]\ar[d]&
H^0_U
(U\times U,
{\cal H})\\
H^0_{\bar f^{-1}(X)}
((X\times X)',
{\cal H}')
\ar[ru]}$$
are isomorphisms.
Thus, it suffices to show the equality
in 
$H^0_U(U\times U,
{\cal H})$.
Hence the assertion
follows from the description
(\ref{eqid})
of the isomorphism (\ref{eqcc1})
in the smooth case.

We consider the pull-back
to $H^0(X,\tilde \delta^*\widetilde {\cal H})$
of the equality
$f^*({\rm id}_{j_!\cal F})
=e\cup [X]$
by the log diagonal map
$\tilde\delta:X\to
(X\times X)^\sim$.
Then, since $\delta=f\circ \tilde \delta$,
we obtain
$\delta^*({\rm id}_{j_!{\cal F}})
=e\cup (X,X)_{(X\times X)^\sim}$
in $H^0(X,\tilde \delta^*
\widetilde {\cal H})
=H^0(X,j_*{\cal E}nd_U{\cal F}(d)[2d])$.
Since the evaluation map
$\tilde \delta^*
\widetilde {\cal H}
\to {\cal K}_X$
is induced by the trace map
$j_*{\cal E}nd_U{\cal F}
\to \Lambda$,
we obtain
$C(j_!{\cal F})
={\rm rank}\ {\cal F}\cdot
(X,X)_{(X\times X)^\sim}$.
\qed

\begin{cor}\label{cortame}
Let $X$ be a smooth scheme 
of dimension $d$ over $k$
and $U=X\setminus D$ be 
the complement
of a divisor $D$
with simple normal crossings.
We keep the notation
in the diagram
{\rm(\ref{eqtame})}.

Let $D^+\subset D$
be the union of some 
irreducible components
and put $U^+=X
\setminus D^+\supset U$.
Let $g:{\mathbb X}\to (X\times X)^\sim$
be a morphism of schemes over $k$
that is an isomorphism
on $(U^+\times U^+)^\sim
=(U^+\times U^+)
\times_{X\times X}
(X\times X)^\sim
\subset
(X\times X)^\sim$
and let $j^\natural:U\times U
\to {\mathbb X}$
be the open immersion.
Let $\delta^\natural:
X\to {\mathbb X}$
be a closed immersion
satisfying
$\tilde \delta=
g\circ \delta^\natural$.
We assume that 
the cycle map
\begin{equation}
\begin{CD}
\Lambda(d)[2d]
@>>> {\cal K}_{{\mathbb X}}=
R(p_2\circ f\circ g)^!\Lambda
\end{CD}\label{eqU3}
\end{equation}
is an isomorphism.
Define the cycle class $[X]
\in H^0_X({\mathbb X},
\Lambda(d)[2d])$
to be the inverse image
of $1$ by
the isomorphism
$H^0_X({\mathbb X},
\Lambda(d)[2d])
\to
H^0_X({\mathbb X},
{\cal K}_{{\mathbb X}})=H^0(X,\Lambda)$
induced by
the isomorphism {\rm (\ref{eqU3})}.

Let ${\cal F}$ be a smooth sheaf on $U$ 
of free $\Lambda$-modules of finite rank.
We assume that 
${\cal F}$ is tamely ramified
along $D\cap U^+=
U^+\setminus U$.
We put ${\cal H}_0=
{\cal H}om({\rm pr}_2^*{\cal F},
{\rm pr}_1^*{\cal F})$ on $U\times U$
and 
$\overline {\cal H}=R{\cal H}om(
{\rm pr}_2^*j_!{\cal F},
{\rm pr}_1^!j_!{\cal F})$
on $X\times X$.
We also put
${\cal H}_0^\natural=
j^\natural_*{\cal H}_0$
on ${\mathbb X}$.
Let $e\in \Gamma(X,
\delta^{\natural *}{\cal H}_0^\natural)$
be a section
such that
the restriction $e|_U
\in \Gamma(U,
\delta_U^*{\cal H}_0)=
{\rm End}_U({\cal F})$
is the identity
of ${\cal F}$.
We put $E={\mathbb X}
\setminus (U^+\times U^+)^\sim$
and assume
\begin{equation}
H^{2d}_E({\mathbb X},
{\cal H}^\natural_0(d))=0.
\label{eqHE0}
\end{equation}

1. There exists a unique map
$(f\circ g)^*\overline {\cal H}
\to {\cal H}^\natural=
{\cal H}_0^\natural(d)[2d]$
inducing the canonical isomorphism
${\cal H}=R{\cal H}om({\rm pr}_2^*{\cal F},
R{\rm pr}_1^!{\cal F})\to
{\cal H}om({\rm pr}_2^*{\cal F},
{\rm pr}_1^*{\cal F})(d)[2d]={\cal H}_0(d)[2d]$
on $U\times U$.

2. We consider the pull-back 
$(f\circ g)^*({\rm id}_{j_!\cal F})
\in H^0_{(f\circ g)^{-1}(X)}(
{\mathbb X},{\cal H}^\natural)$
of the identity
${\rm id}_{j_!\cal F}\in 
{\rm End}_X(j_!{\cal F})=
H^0_X(X\times X,\overline{\cal H})$.
Then, we have
\begin{equation}
(f\circ g)^*({\rm id}_{j_!\cal F})
=e\cup [X]
\label{eqXn0}
\end{equation}
in $H^0_{(f\circ g)^{-1}(X)\cup E}(
{\mathbb X},{\cal H}^\natural)$.
Consequently,
we have
\begin{equation}
C(j_!{\cal F})=
{\rm rank}\ {\cal F}\cdot
(X,X)_{\mathbb X}
\label{eqXn1}
\end{equation}
in $H^0(X,{\cal K}_X)$
where 
$(X,X)_{\mathbb X}$
denotes pull-back of the cycle class
$[X]\in H^{2d}_X({\mathbb X},\Lambda(d))$.
\end{cor}
{\it Proof.}
1.
Since the image of
$(X\times X)^\sim
\setminus
(U\times U)$
in $X\times X$
is a subset of $D\times X$,
we have
$(f\circ g)^*\overline {\cal H}
=j^\natural_!{\cal H}$
as in the proof of Proposition \ref{prtame}.
Hence the assertion follows.

2.
Since $g^{-1}(U^+)
=(f\circ g)^{-1}(X)
\cap
(U^+\times U^+)^\sim$
is the complement
of $E\cap (f\circ g)^{-1}(X)$
in $(f\circ g)^{-1}(X)$,
we have an exact sequence
$$\begin{CD}
H^0_{E}
({\mathbb X},{\cal H}^\natural)
@>>>
H^0_{(f\circ g)^{-1}(X)\cup E}
({\mathbb X},{\cal H}^\natural)
@>>>
H^0_{g^{-1}(U^+)}
((U^+\times U^+)^\sim,{\cal H}^\natural).
\end{CD}$$
By the assumption
(\ref{eqHE0}),
we have
$H^0_{E}
({\mathbb X},{\cal H}^\natural)=0$.
Hence the restriction map

\noindent
$H^0_{(f\circ g)^{-1}(X)\cup E}
({\mathbb X},{\cal H}^\natural)
\to
H^0_{g^{-1}(U^+)}
((U^+\times U^+)^\sim,{\cal H}^\natural)$
is an injection.
By Proposition \ref{prtame},
the equality 
$(f\circ g)^*({\rm id}_{j_!\cal F})
=e\cup [X]$
holds in $H^0_{g^{-1}(U^+)}
((U^+\times U^+)^\sim,{\cal H}^\natural)$.
Thus we obtain
$(f\circ g)^*({\rm id}_{j_!\cal F})
=e\cup [X]$
(\ref{eqXn0})
in $H^0_{(f\circ g)^{-1}(X)\cup E}(
{\mathbb X},{\cal H}^\natural)$.

The equality
(\ref{eqXn1})
is deduced from
(\ref{eqXn0})
as in the proof of
Proposition \ref{prtame}.
\qed

\begin{thm}\label{thmcc}
Let $X$ be a smooth scheme 
of dimension $d$ over $k$
and $U=X\setminus D$ be 
the complement
of a divisor $D$
with simple normal crossings.
Let $R=\sum_ir_iD_i\ge 0$
be an effective
divisor with rational coefficients.
Let $g:(X\times X)^{(R)}\to 
(X\times X)^\sim$
and $\delta^{(R)}:X
\to (X\times X)^{(R)}$
be as in \S {\rm 2.3}
and let
$j^{(R)}:U\times U
\to (X\times X)^{(R)}$
be the open immersion.

Let ${\cal F}$ be a smooth sheaf on $U$ 
of free $\Lambda$-modules of finite rank.
We put ${\cal H}_0=
{\cal H}om({\rm pr}_2^*{\cal F},
{\rm pr}_1^*{\cal F})$ on $U\times U$.
We assume that the log ramification of
${\cal F}$ is bounded by $R+$
and let
$e\in \Gamma(X,
\delta^{(R)*}j^{(R)}_*
{\cal H}_0)$
be the unique section
whose image
by the base change map
in 
$\Gamma(X,
j_*\delta_U^*{\cal H}_0)
={\rm End}_U({\cal F})$
is the identity of ${\cal F}$.
We further assume that
the log ramification of ${\cal F}$
along $D$ is non-degenerate
with respect to $R$
(cf.\ Definition {\rm \ref{dfdlsp}}).

Then, we have
\begin{eqnarray}
C(j_!{\cal F})
&=&
{\rm rank}\ {\cal F}\cdot
(X,X)_{(X\times X)^{(R)}}
\label{eqXX3}
\\
&=&
(-1)^d\cdot
{\rm rank}\ {\cal F}\times
\label{eqXX4}
\\
&&\qquad
(c_d(\Omega^1_X(\log D))
+c(\Omega^1_X(\log D))\cap
(1-R)^{-1}\cap [R]))_{\dim d}
\nonumber
\end{eqnarray}
in $H^0(X,{\cal K}_X)$.
\end{thm}
{\it Proof.}
We put
$D^+=\sum_{i:r_i>0}D_i$.
We verify
that
$g:{\mathbb X}=
(X\times X)^{(R)}
\to (X\times X)^\sim$
satisfies
the assumptions
in Corollary \ref{cortame}.
By the construction,
the map
$g:{\mathbb X}=
(X\times X)^{(R)}
\to (X\times X)^\sim$
is an isomorphism on the complement of
$D^+\subset X
\subset (X\times X)^\sim$.
The log diagonal map
$\tilde \delta:X\to 
(X\times X)^\sim$
is lifted to a closed immersion
$\delta^{(R)}:X\to (X\times X)^{(R)}$.
The cycle map
$\Lambda(d)[2d]\to {\cal K}_{(X\times X)^{(R)}}$
is an isomorphism
by Proposition \ref{prPr}.1.

By the definition of $D^+$
and by the assumption
that the log ramification
of ${\cal F}$ is bounded by $R+$,
it follows that
${\cal F}$ is tamely ramified
along $D\setminus D^+=
U^+\setminus U$
by Corollary \ref{corbd0}.
The complement
$(X\times X)^{(R)}
\setminus (U^+\times U^+)^\sim$
equals the inverse image
$E^+$ of $D^+$.
We show that 
${\cal H}_0^{(R)}=
j^{(R)}_*{\cal H}_0$
satisfies the assumption
$$H^{2d}_{E^+}((X\times X)^{(R)},
{\cal H}^{(R)}_0(d))=0$$
(\ref{eqHE0}).
Let $i:E^+\to 
(X\times X)^{(R)}$
be the closed immersion
and $p:E^+\to D^+$
be the projection.
Since
$H^{2d}_{E^+}((X\times X)^{(R)},
{\cal H}^{(R)}_0(d))=
H^{2d}(D^+,
Rp_*Ri^!j^{(R)}_*{\cal H}_0(d))$,
it suffices to
show
$Rp_*R^qi^!j^{(R)}_*{\cal H}_0=0$
for $q\ge 0$.
Since
$R^qi^!j^{(R)}_*{\cal H}_0$ is 0
for $q=0,1$ and
is isomorphic to
$R^{q-1}j^{(R)}_*{\cal H}_0$
for $q>1$,
it follows from 
$Rp_*R^qj^{(R)}_*{\cal H}_0=0$
proved in
Lemma \ref{lmSH}.2.
Thus, the assumptions in
Corollary \ref{cortame}
are satisfied and we obtain
the equality (\ref{eqXX3}).

The equality (\ref{eqXX4})
follows
from the equalities
(\ref{eqXX3})
and 
(\ref{eqPr1})
and the isomorphism
${\cal N}_{X/(X\times X)^\sim}
\to
\Omega^1_X(\log D)$.
\qed

\subsection{Characteristic cycle}

Let $X$ be a smooth scheme
of dimension $d$ over $k$
and let $D$ be
a divisor with simple normal crossings.
Let 
$$T^*X(\log D)={\mathbf V}
(\Omega^1_X(\log D)^\vee)$$
be the logarithmic cotangent
bundle.
We regard
$X$ as a closed subscheme of
$T^*X(\log D)$
by the $0$-section.
Let $D_i$ be an irreducible component
of $D$, $\xi_i$ be the
generic point of $D_i$
and $K_i$ be the
fraction field of
the henselization
${\cal O}_{X,\xi_i}^h$ of the local ring.
The residue field $F_i$
of $K_i$ is the function field of $D_i$.

Let $r>0$ be a rational number
and $\chi:{\rm Gr}^r_{\log}G_{K_i}
\to {\mathbb F}_p$
be a non-trivial character.
The refined Swan character
${\rm rsw}\ \chi\in 
\Omega^1_{F_i}(\log)
\otimes 
{\mathfrak m}_{{\overline K}_i}^{(-r)}
/{\mathfrak m}_{{\overline K}_i}^{(-r)+}$
regarded as an injection
${\mathfrak m}_{{\overline K}_i}^r
/{\mathfrak m}_{{\overline K}_i}^{r+}
\to
\Omega^1_X(\log D)_{\xi_i}
\otimes \overline F_i$
defines a line
in the $\overline F_i$-vector space
$\Omega^1_X(\log D)_{\xi_i}
\otimes \overline F_i$
and hence an $\overline F_i$-valued point
$[{\rm rsw}\ \chi]:
{\rm Spec}\ \overline F_i
\to 
{\mathbf P}(\Omega^1_X(\log D)^\vee)$.
We define a reduced closed subscheme
$T_{\chi}
\subset
{\mathbf P}(\Omega^1_X(\log D)^\vee)$ 
to be the Zariski closure
$\overline{\{[{\rm rsw}\ \chi]
({\rm Spec}\ \overline F_i)\}}$
of the image
and let $L_\chi={\mathbf V}({\cal O}_{T_\chi}(1))$
be the pull-back to $T_\chi$
of the tautological
sub line bundle 
$L\subset T^*X(\log D)\times_X
{\mathbf P}(\Omega^1_X(\log D)^\vee)$.
We have a commutative diagram
\begin{equation}
\begin{CD}
L_\chi@>>> 
T^*X(\log D)\times_XD_i
@>>>
T^*X(\log D)&={\mathbf V}(\Omega^1_X(\log D)^\vee)\\
@VVV @VVV@VVV\\
T_\chi@>{\pi_\chi}>>D_i
@>>> X.
\end{CD}
\label{eqSS}
\end{equation}
The natural map
$\pi_\chi:T_\chi\to D_i$
is generically finite.

Let ${\cal F}$ be a 
smooth $\ell$-adic sheaf on
$U=X\setminus D$
and $R=\sum_ir_iD_i$
be an effective divisor
with rational coefficients
$r_i\ge0$.
In the rest of the paper,
we assume that 
${\cal F}$ satisfies
the following conditions:
\begin{itemize}
\item[(R)]
The log ramification of ${\cal F}$
along $D$ is bounded by $R+$.

\item[(C)]
For each irreducible component
$D_i$ of $D$,
the closure 
$\overline{S^{(n\cdot R)}_{{\cal F},\xi_i}}$
of the generic fiber of the dual support
is finite over $D_i$
and its intersection
$\overline{S^{(n\cdot R)}_{{\cal F},\xi_i}}
\cap D_i$
with the 0-section of
${\mathbf V}_n(
\Omega^1_X(\log D)^\vee,
{\cal O}(-nR))_{D_i}$
is empty.
\end{itemize}
By Lemma \ref{lmSH}, the conditions
imply ${\cal F}_{\bar \eta_i}=
{\cal F}_{\bar \eta_i}^{(r_i)}$
for every irreducible component
$D_i$ of $D$.
They also imply
that the log ramification of
${\cal F}$ is non-degenerate
with respect to $R$ by
Proposition \ref{prR}.3.

\begin{lm}\label{lmSS}
Assume that the log ramification
of ${\cal F}$ 
is bounded by $R+$
and that ${\cal F}$
satisfies the conditions
{\rm (R)}
and {\rm (C)}.
Let $D_i$
be an irreducible component of
$D$
and $\chi$
be a character of
${\rm Gr}^{r_i}_{\log}G_{K_i}$
appearing in 
the direct sum decomposition
${\cal F}_{\bar \eta_i}=
\sum_{\chi}n_{\chi}\chi$.

1.
The scheme $T_\chi$
is finite over $D_i$.

2.
We put
\begin{equation}
SS_\chi=\frac1{[T_{\chi}:D_i]}\pi_{\chi*}[L_{\chi}]
\label{eqSSc}
\end{equation}
in $Z_d(T^*X(\log D)\times_XD_i)
_{\mathbb Q}$
in the notation in
{\rm (\ref{eqSS})}.
Then, we have
$$SS_{\chi}
=\left(c(\Omega^1_X(\log D))\cap
(1-R)^{-1}\cap 
[T^*X(\log D)\times_X{D_i}])\right)_{\dim d}.$$
\end{lm}
{\it Proof.}
1. 
By Lemma \ref{lmSH},
the generic fiber
$S^{(n\cdot R)}_{{\cal F},\xi_i}$
of the dual support
consists of
the refined Swan characters
${\rm rsw}\ \chi$
of the characters $\chi$
appearing in 
the direct sum decomposition
${\cal F}_{\bar \eta_i}=
\sum_{\chi}n_{\chi}\chi$.
Hence, 
by the condition (C),
the closure
$\overline{S^{(n\cdot R)}_{{\cal F},\xi_i}}$
is a closed subscheme of
${\mathbf V}_n(
\Omega^1_X(\log D)^\vee,
{\cal O}(-nR))_{D_i}
\setminus D_i$
finite over $D_i$.
Hence, 
the union 
$\bigcup_\chi T_{\chi}$
is the image of
$\overline{S^{(n\cdot R)}_{{\cal F},\xi_i}}$
by the canonical map
$\varphi:
{\mathbf V}_n(
\Omega^1_X(\log D)^\vee,
{\cal O}(-nR))_{D_i}
\setminus D_i
\to 
{\mathbf P}(
\Omega^1_X(\log D)^\vee)$.

2.
Since the conormal sheaf of
$L_\chi
\subset 
T^*X(\log D)\times_XT_\chi$
is the pull-back of
the locally free sheaf
${\rm Ker}(\Omega^1_X(\log D)^\vee\to {\cal O}(1))$
of rank $d-1$,
we have
$$[L_\chi]
=(-1)^{d-1}
c_{d-1}({\rm Ker}(\Omega^1_X(\log D)^\vee\to {\cal O}(1)))
\cap [T^*X(\log D)\times_XT_\chi].$$
Hence we have
$$SS_\chi=
(c(\Omega^1_X(\log D))\cap
c({\cal O}(-1))^{-1}\cap 
[T^*X(\log D)\times_X{D_i}])_{\dim d}.$$
By Lemma \ref{lmOn},
the pull-back of ${\cal O}(n)$
on ${\mathbf V}_n(
\Omega^1_X(\log D)^\vee,
{\cal O}(-nR))_{D_i}
\setminus D_i \supset
\overline{S_{{\cal F},\xi_i}^{(n\cdot R)}}$ 
is 
canonically isomorphic to
${\cal O}(nR)$.
Since 
the union 
$\bigcup_\chi T_{\chi}$
is the image of
$\overline{S^{(n\cdot R)}_{{\cal F},\xi_i}}$,
the assertion follows.
\qed

Assuming the conditions (R) and (C),
we define the characteristic cycle
$CC({\cal F})$
as a rational $d$-cycle
on $T^*X(\log D)$.

\begin{df}\label{dfCC}
Let ${\cal F}$
be a smooth $\Lambda$-sheaf on
$U=X\setminus D$
satisfying the conditions
{\rm (R)}
and {\rm (C)}.
For an irreducible component $D_i$
of $D$ with $r_i>0$,
let ${\cal F}_{\bar \eta_i}=
\sum_{\chi}n_{\chi}\chi$
be the direct sum decomposition
of the representation
induced on ${\rm Gr}^{r_i}_{\log}G_{K_i}$.
We define the
characteristic cycle by 
\begin{equation}
CC({\cal F})=
(-1)^d\left(
{\rm rank}\ {\cal F}
\cdot [X]
+\sum_{i,r_i>0}
r_i\cdot \sum_{\chi}
n_{\chi}\cdot
[SS_{\chi}]
\right)
\label{CC}
\end{equation}
in $Z^d(T^*X(\log D))_{\mathbb Q}$.
\end{df}

If $\dim X=1$,
we put
${\rm Sw}\ {\cal F}
=
\sum_{x\in D}{\rm Sw}_x\ {\cal F}
\cdot [x]\in Z_0(X)$
and let
$p:T^*X(\log D)
\to X$ be the projection.
Then, we have
$$
CC({\cal F})
=-\left(
{\rm rank}\ {\cal F}
\cdot [X]
+
p^*[{\rm Sw}\ {\cal F}]\right).
$$

\begin{thm}\label{thmCC}
Let $X$ be a smooth scheme
of dimension $d$ over
$k$ and $D$ be a divisor
with simple normal crossings.
Let ${\cal F}$
be a smooth $\ell$-adic sheaf on
$U=X\setminus D$
satisfying the conditions
{\rm (R)}
and {\rm (C)}.
Then we have
$$C(j_!{\cal F})=[CC({\cal F})]$$
in $H^{2d}(X,\Lambda(d))=
H^{2d}(T^*X(\log D),\Lambda(d))$.
In other words,
we have
$$C(j_!{\cal F})=(CC({\cal F}),
X)_{T^*X(\log D)}.$$
\end{thm}
{\it Proof.}
By the assumption (C)
and by Lemma \ref{lmSH},
the assumption in Theorem \ref{thmcc}
is satisfied.
Hence the left hand side is
equal to
$${\rm rank}\ {\cal F}\cdot
(-1)^d\cdot
(c_d(\Omega^1_X(\log D))
+c(\Omega^1_X(\log D))\cap
(1-R)^{-1}\cap [R]))_{\dim 0}).$$
By Lemma \ref{lmSS},
the right hand side is
also equal to this.
\qed

By the index formula
(\ref{eqidx}),
Theorem \ref{thmCC}
implies the following.

\begin{cor}
Further if $X$ is proper,
we have
$$\chi_c(U_{\bar k},{\cal F})=
\deg (CC({\cal F}),
X)_{T^*X(\log D)}.$$
\end{cor}

\begin{tabular}{l@{\qquad}l}
Takeshi Saito\\
Department of Mathematical Sciences,\\
University of Tokyo,\\
Tokyo 153-8914 Japan\\
t-saito@ms.u-tokyo.ac.jp
\end{tabular}

\end{document}